\documentclass{gtart}

\def\ifplaintex{\expandafter\ifx\csname documentclass\endcsname\relax}

\ifplaintex 
\else
\fi

\expandafter\ifx\csname beginpicture\endcsname\relax
\expandafter\ifx\csname documentclass\endcsname\relax
\input pictex \else
\input prepictex \input pictex \input postpictex \fi\fi

\def\gt{{\mathsurround=0pt\it $\cal G\mskip-2mu$eometry \&\ 
$\cal T\!\!$opology}}        

\def\gtp{{\mathsurround=0pt\it $\cal G\mskip-2mu$eometry \&\ 
$\cal T\!\!$opology $\cal P\!$ublications}}  

\def\lognumber#1{\def\thelognumber{#1}}
\def\volumenumber#1{\def\thevolumenumber{#1}}
\def\papernumber#1{\def\thepapernumber{#1}}
\def\volumeyear#1{\def\thevolumeyear{#1}}

\def\pagenumbers#1#2{\def\startpage{#1}\def\finishpage{#2}}
\def\published#1{\def\publishdate{#1}}
\def\proposed#1{\def\theproposer{#1}}
\def\seconded#1{\def\theseconders{#1}}
\def\received#1{\def\receiveddate{#1}}

\def\accepted#1{\def\accepteddate{#1}}
\def\asciititle#1{\def\theasciititle{#1}}

\long\def\asciiabstract#1{\long\def\theasciiabstract{#1}}
\def\asciikeywords#1{\def\theasciikeywords{#1}}

\let\\\par\let\thelognumber\relax
\let\thevolumenumber\relax\let\thepapernumber\relax
\let\thevolumeyear\relax\let\thesamplenumber\relax\let\startpage\relax
\let\finishpage\relax\let\publishdate\relax\let\receiveddate\relax
\let\reviseddate\relax\let\accepteddate\relax\let\theasciititle\relax
\let\theasciiauthors\relax
\let\theasciiabstract\relax\let\theasciikeywords\relax
\let\theasciiemail\relax\let\theshortauthors\relax\let\theshorttitle\relax

\long\def\maketitlep{   

\count0=\startpage

\gt\hfill      
\beginpicture
\setcoordinatesystem units <0.33truein, 0.33truein> point at 2.2 0.9
\setplotsymbol ({$\cal G$})
\plotsymbolspacing=9truept
\circulararc 315 degrees from 0 1 center at 0 0
\setplotsymbol ({$\cal T$})
\circulararc 315 degrees from 1 -1 center at 1 0
\endpicture
\break
{\small\ifx\thesamplenumber\relax 
Volume \else Sample
\fi\thevolumenumber\ (\thevolumeyear)
\startpage--\finishpage\nl
Published: \publishdate}
\vglue 0.5truein plus 0.4fil minus 0.1truein

{\parskip=0pt\leftskip 0pt plus 1fil\def\\{\par\smallskip}{\ifplaintex\large
\else\Large\fi\bf\thetitle}\par\medskip}   

\vglue 0pt plus 0.1fil 

{\parskip=0pt\leftskip 0pt plus 1fil\def\\{\par}{\sc\theauthors}
\par\medskip}

\vglue 0pt plus 0.1fil 

{\small\parskip=0pt\let\newline\\
{\leftskip 0pt plus 1fil\def\\{\par}{\sl\theaddress}\par}
\expandafter\ifx\theemail\relax    
\relax\else\vglue 5pt plus 0.02fil minus 2pt\def\\{\stdspace{\rm 
and}\stdspace} 
\cl{Email:\stdspace\tt\theemail}\fi
\ifx\theurl\relax                  
\relax\else\vglue 5pt plus 0.02fil minus 2pt\def\\{\stdspace{\rm 
and}\stdspace}
\cl{URL:\stdspace\tt\theurl}\fi\par}

\vglue 7pt plus 0.3fil minus 3pt

{\bf Abstract}
\vglue 5pt plus 0.1fil minus 2pt

\theabstract

\vglue 7pt plus 0.3fil minus 3pt

{\bf AMS Classification numbers}\quad Primary:\quad \theprimaryclass

Secondary:\quad \thesecondaryclass

\vglue 5pt plus 0.3fil minus 2pt

{\bf Keywords}\quad \thekeywords

\vglue 10pt plus 0.5fil minus 5pt

{\small  Proposed: \theproposer\hfill Received: \receiveddate\nl
Seconded: \theseconders\hfill 
\ifx\reviseddate\relax                         
Accepted: \accepteddate                        
\else
Revised: \reviseddate                          
\fi}
\eject
}       

\let\maketitlepage\maketitlep
\let\maketitle\maketitlepage

\font\phead=cmsl9 scaled 950
\font=cmsl9 scaled 1050
\font\pnum=cmbx10 scaled 913
\font\lnum=cmbx10 
\font\pfoot=cmsl9 scaled 950
\font=cmsl9 scaled 1050
\ifplaintex
\headline{\vbox to 0pt{\vskip -4.5mm\line{\small\phead\ifnum
\count0=\startpage ISSN 1364-0380 (on line)
1465-3060 (printed) \hfill {\pnum\folio}\else\ifodd\count0\def\\{ }%
\ifx\theshorttitle\relax\thetitle\else\theshorttitle\fi\hfill{\pnum\folio}
\else\def\\{ and }{\pnum\folio}\hfill\ifx\theshortauthors\relax\theauthors
\else\theshortauthors\fi\fi\fi}\vss}}
\footline{\vbox to 0pt{\vglue 0mm\line{\small\pfoot\ifnum\count0=\startpage
\copyright\ \gtp\hfill\else
\gt, Volume \thevolumenumber\ (\thevolumeyear)\hfill\fi}\vss
}}
\else
\makeatletter
\def\@oddhead{{\small\lhead\ifnum\count0=\startpage ISSN 1364-0380 (on line)
1465-3060 (printed) \hfill {\lnum\number\count0}\else\ifodd\count0
\def\\{ }\ifx\theshorttitle\relax \thetitle \else\theshorttitle\fi\hfill
{\lnum\number\count0}\else\def\\{ and }{\lnum\number\count0}
\hfill\ifx\theshortauthors\relax 
\theauthors\else\theshortauthors\fi\fi\fi}}\def\@evenhead{\@oddhead}
\def\@oddfoot{\small\lfoot\ifnum\count0=\startpage\copyright\ \gtp\hfill\else
\gt, Volume \thevolumenumber\ (\thevolumeyear)\hfill\fi}
\def\@evenfoot{\@oddfoot}
\makeatother
\fi

\newwrite\gtoutfile
\long\gdef\makeheadfile{  
{\def\\{, }\def\s{ }
\immediate\openout\gtoutfile head.xxx
\immediate\write\gtoutfile{To: math@arxiv.org}
\immediate\write\gtoutfile{Subject: put or rep NNNNN:pppp}
\immediate\write\gtoutfile{--text follows this line--}
\immediate\write\gtoutfile{Proxy-for: \ifx\theasciiauthors\relax
\theauthors\else\theasciiauthors\fi\s<\ifx\theasciiemail\relax\theemail\else\theasciiemail\fi>}
\immediate\write\gtoutfile{\noexpand\\}
\immediate\write\gtoutfile{Authors: \ifx\theasciiauthors\relax
\theauthors\else\theasciiauthors\fi}
{\def\\{ }\immediate\write\gtoutfile{Title: \ifx\theasciititle\relax
\thetitle\else\theasciititle\fi}}
\immediate\write\gtoutfile{Subj-class: GT or SG or MG etc}
\immediate\write\gtoutfile{MSC-class: \theprimaryclass\ifx\thesecondaryclass\relax\else, \thesecondaryclass\fi}
\immediate\write\gtoutfile{Journal-ref: Geom. Topol. \thevolumenumber
(\thevolumeyear) \startpage-\finishpage}
\immediate\write\gtoutfile{Comments: Published by Geometry and Topology at}
\immediate\write\gtoutfile{\s\s http://www.maths.warwick.ac.uk/gt/GTVol\thevolumenumber/paper\thepapernumber.abs.html}
\immediate\write\gtoutfile{\noexpand\\}
\immediate\write\gtoutfile{}
\ifx\theasciiabstract\relax
\immediate\write\gtoutfile{\theabstract}\else
\immediate\write\gtoutfile{\theasciiabstract}\fi
\immediate\write\gtoutfile{}
\immediate\write\gtoutfile{\noexpand\\}
\immediate\write\gtoutfile{}
\immediate\closeout\gtoutfile}}  

\def\maketitlepage{\maketitlep\makeheadfile}
\let\maketitle\maketitlepage

\def\ifplaintex{\expandafter\ifx\csname documentclass\endcsname\relax}

\ifplaintex 
\else
\fi

\expandafter\ifx\csname beginpicture\endcsname\relax
\expandafter\ifx\csname documentclass\endcsname\relax
\input pictex \else
\input prepictex \input pictex \input postpictex \fi\fi

\def\gt{{\mathsurround=0pt\it $\cal G\mskip-2mu$eometry \&\ 
$\cal T\!\!$opology}}        

\def\gtp{{\mathsurround=0pt\it $\cal G\mskip-2mu$eometry \&\ 
$\cal T\!\!$opology $\cal P\!$ublications}}  

\def\lognumber#1{\def\thelognumber{#1}}
\def\volumenumber#1{\def\thevolumenumber{#1}}
\def\papernumber#1{\def\thepapernumber{#1}}
\def\volumeyear#1{\def\thevolumeyear{#1}}

\def\pagenumbers#1#2{\def\startpage{#1}\def\finishpage{#2}}
\def\published#1{\def\publishdate{#1}}
\def\proposed#1{\def\theproposer{#1}}
\def\seconded#1{\def\theseconders{#1}}
\def\received#1{\def\receiveddate{#1}}

\def\accepted#1{\def\accepteddate{#1}}
\def\asciititle#1{\def\theasciititle{#1}}

\long\def\asciiabstract#1{\long\def\theasciiabstract{#1}}
\def\asciikeywords#1{\def\theasciikeywords{#1}}

\let\\\par\let\thelognumber\relax
\let\thevolumenumber\relax\let\thepapernumber\relax
\let\thevolumeyear\relax\let\thesamplenumber\relax\let\startpage\relax
\let\finishpage\relax\let\publishdate\relax\let\receiveddate\relax
\let\reviseddate\relax\let\accepteddate\relax\let\theasciititle\relax
\let\theasciiauthors\relax
\let\theasciiabstract\relax\let\theasciikeywords\relax
\let\theasciiemail\relax\let\theshortauthors\relax\let\theshorttitle\relax

\long\def\maketitlep{   

\count0=\startpage

\gt\hfill      
\beginpicture
\setcoordinatesystem units <0.33truein, 0.33truein> point at 2.2 0.9
\setplotsymbol ({$\cal G$})
\plotsymbolspacing=9truept
\circulararc 315 degrees from 0 1 center at 0 0
\setplotsymbol ({$\cal T$})
\circulararc 315 degrees from 1 -1 center at 1 0
\endpicture
\break
{\small\ifx\thesamplenumber\relax 
Volume \else Sample
\fi\thevolumenumber\ (\thevolumeyear)
\startpage--\finishpage\nl
Published: \publishdate}
\vglue 0.5truein plus 0.4fil minus 0.1truein

{\parskip=0pt\leftskip 0pt plus 1fil\def\\{\par\smallskip}{\ifplaintex\large
\else\Large\fi\bf\thetitle}\par\medskip}   

\vglue 0pt plus 0.1fil 

{\parskip=0pt\leftskip 0pt plus 1fil\def\\{\par}{\sc\theauthors}
\par\medskip}

\vglue 0pt plus 0.1fil 

{\small\parskip=0pt\let\newline\\
{\leftskip 0pt plus 1fil\def\\{\par}{\sl\theaddress}\par}
\expandafter\ifx\theemail\relax    
\relax\else\vglue 5pt plus 0.02fil minus 2pt\def\\{\stdspace{\rm 
and}\stdspace} 
\cl{Email:\stdspace\tt\theemail}\fi
\ifx\theurl\relax                  
\relax\else\vglue 5pt plus 0.02fil minus 2pt\def\\{\stdspace{\rm 
and}\stdspace}
\cl{URL:\stdspace\tt\theurl}\fi\par}

\vglue 7pt plus 0.3fil minus 3pt

{\bf Abstract}
\vglue 5pt plus 0.1fil minus 2pt

\theabstract

\vglue 7pt plus 0.3fil minus 3pt

{\bf AMS Classification numbers}\quad Primary:\quad \theprimaryclass

Secondary:\quad \thesecondaryclass

\vglue 5pt plus 0.3fil minus 2pt

{\bf Keywords}\quad \thekeywords

\vglue 10pt plus 0.5fil minus 5pt

{\small  Proposed: \theproposer\hfill Received: \receiveddate\nl
Seconded: \theseconders\hfill 
\ifx\reviseddate\relax                         
Accepted: \accepteddate                        
\else
Revised: \reviseddate                          
\fi}
\eject
}       

\let\maketitlepage\maketitlep
\let\maketitle\maketitlepage

\font\phead=cmsl9 scaled 950
\font=cmsl9 scaled 1050
\font\pnum=cmbx10 scaled 913
\font\lnum=cmbx10 
\font\pfoot=cmsl9 scaled 950
\font=cmsl9 scaled 1050
\ifplaintex
\headline{\vbox to 0pt{\vskip -4.5mm\line{\small\phead\ifnum
\count0=\startpage ISSN 1364-0380 (on line)
1465-3060 (printed) \hfill {\pnum\folio}\else\ifodd\count0\def\\{ }%
\ifx\theshorttitle\relax\thetitle\else\theshorttitle\fi\hfill{\pnum\folio}
\else\def\\{ and }{\pnum\folio}\hfill\ifx\theshortauthors\relax\theauthors
\else\theshortauthors\fi\fi\fi}\vss}}
\footline{\vbox to 0pt{\vglue 0mm\line{\small\pfoot\ifnum\count0=\startpage
\copyright\ \gtp\hfill\else
\gt, Volume \thevolumenumber\ (\thevolumeyear)\hfill\fi}\vss
}}
\else
\makeatletter
\def\@oddhead{{\small\lhead\ifnum\count0=\startpage ISSN 1364-0380 (on line)
1465-3060 (printed) \hfill {\lnum\number\count0}\else\ifodd\count0
\def\\{ }\ifx\theshorttitle\relax \thetitle \else\theshorttitle\fi\hfill
{\lnum\number\count0}\else\def\\{ and }{\lnum\number\count0}
\hfill\ifx\theshortauthors\relax 
\theauthors\else\theshortauthors\fi\fi\fi}}\def\@evenhead{\@oddhead}
\def\@oddfoot{\small\lfoot\ifnum\count0=\startpage\copyright\ \gtp\hfill\else
\gt, Volume \thevolumenumber\ (\thevolumeyear)\hfill\fi}
\def\@evenfoot{\@oddfoot}
\makeatother
\fi

\newwrite\gtoutfile
\long\gdef\makeheadfile{  
{\def\\{, }\def\s{ }
\immediate\openout\gtoutfile head.xxx
\immediate\write\gtoutfile{To: math@arxiv.org}
\immediate\write\gtoutfile{Subject: put or rep NNNNN:pppp}
\immediate\write\gtoutfile{--text follows this line--}
\immediate\write\gtoutfile{Proxy-for: \ifx\theasciiauthors\relax
\theauthors\else\theasciiauthors\fi\s<\ifx\theasciiemail\relax\theemail\else\theasciiemail\fi>}
\immediate\write\gtoutfile{\noexpand\\}
\immediate\write\gtoutfile{Authors: \ifx\theasciiauthors\relax
\theauthors\else\theasciiauthors\fi}
{\def\\{ }\immediate\write\gtoutfile{Title: \ifx\theasciititle\relax
\thetitle\else\theasciititle\fi}}
\immediate\write\gtoutfile{Subj-class: GT or SG or MG etc}
\immediate\write\gtoutfile{MSC-class: \theprimaryclass\ifx\thesecondaryclass\relax\else, \thesecondaryclass\fi}
\immediate\write\gtoutfile{Journal-ref: Geom. Topol. \thevolumenumber
(\thevolumeyear) \startpage-\finishpage}
\immediate\write\gtoutfile{Comments: Published by Geometry and Topology at}
\immediate\write\gtoutfile{\s\s http://www.maths.warwick.ac.uk/gt/GTVol\thevolumenumber/paper\thepapernumber.abs.html}
\immediate\write\gtoutfile{\noexpand\\}
\immediate\write\gtoutfile{}
\ifx\theasciiabstract\relax
\immediate\write\gtoutfile{\theabstract}\else
\immediate\write\gtoutfile{\theasciiabstract}\fi
\immediate\write\gtoutfile{}
\immediate\write\gtoutfile{\noexpand\\}
\immediate\write\gtoutfile{}
\immediate\closeout\gtoutfile}}  

\def\maketitlepage{\maketitlep\makeheadfile}
\let\maketitle\maketitlepage

\lognumber{273}

\volumenumber{6}
\papernumber{16}
\volumeyear{2002}
\pagenumbers{425}{494}
\received{4 September 2002}
\accepted{2 October 2002}
\published{23 October 2002}
\proposed{Robion Kirby}
\seconded{Ronald Stern, Yasha Eliashberg}

\usepackage{graphicx,latexsym}
\nocolon

\theoremstyle{plain}
\newtheorem{Thm}{Theorem}
\newtheorem{Lem}[Thm]{Lemma}
\newtheorem{Cor}[Thm]{Corollary}
\newtheorem{Prop}[Thm]{Proposition}

\theoremstyle{remark}
\newtheorem{Rm}[Thm]{Remark}

\def\<{\langle}
\def\>{\rangle}

\begin{document}

\title{Cappell--Shaneson's 4--dimensional  $s$--cobordism}
\asciititle{Cappell-Shaneson's 4-dimensional  s-cobordism}
\author{Selman Akbulut} 
\asciikeywords{s-cobordism, quaternionic space}
\address{Department of Mathematics\\Michigan State University\\MI, 48824, USA}
\email{akbulut@math.msu.edu}

\begin{abstract}
In 1987 S Cappell and J Shaneson constructed an $s$--cobordism $H$ from the
quaternionic 3--manifold $Q$ to itself, and asked whether $H$ or any of its
covers are trivial product cobordism? In this paper we study $H$, and
in particular show that its 8--fold cover is the product cobordism from
$S^3$ to itself. We reduce the triviality of $H$ to a question about the
3--twist spun trefoil knot in $S^4$, and also relate this to a question
about a Fintushel--Stern knot surgery.
\end{abstract}

\asciiabstract{In 1987 S Cappell and J Shaneson constructed an
s-cobordism H from the quaternionic 3-manifold Q to itself, and asked
whether H or any of its covers are trivial product cobordism? In this
paper we study H, and in particular show that its 8-fold cover is the
product cobordism from S^3 to itself. We reduce the triviality of H to
a question about the 3-twist spun trefoil knot in S^4, and also relate
this to a question about a Fintushel-Stern knot surgery.}

\maketitlepage

\setcounter{section}{-1}
\section{Introduction}

Let $Q^{3}=S^{3}/{\bf Q_{8}}$ be the quaternionic $3$--manifold, obtained as the
quotient of the $3$--sphere by the free action of the quaternionic group ${\bf
Q_{8}}$ of order eight, which can be presented by
${\bf Q_{8}}=\<i,j,k\;|\; i^{2}=j^{2}=k^{2}=-1, \; ij=k, jk=i, ki=j\> $. Also $Q$ is
the $2$--fold branched covering space of $ S^3 $ branched over the three Hopf circles;
combining this with the Hopf map $ S^3\to S^2 $ one sees that $Q$ is a Seifert
Fibered space with three singular fibers. $Q$ is also the $3$--fold branched covering space of
$ S^3 $ branched over the trefoil knot.
$Q$ can also be identified with the boundaries of the
$4$--manifolds of Figure 1 (one can easily
check that the above three definitions are equivalent to this one by drawing framed link
pictures). The second manifold $W$ of Figure 1, consisting
of a $1$-- and $2$--handle pair, is a Stein surface by
\cite{g1}. It is easily seen that $ W $ is a disk bundle over ${\bf RP}^{2}$ obtained as the
tubular neighborhood of an imbedded
${\bf RP}^{2}$ in $S^{4}$. The complement of this imbedding is also a copy of $W$,
decomposing $S^4 =W\smile_{\partial} W$.

\begin{figure}[ht!]
  \begin{center}
     \includegraphics{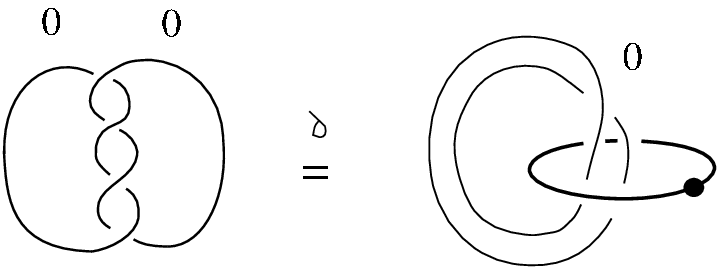}
   \caption{} \label{B1}
    \end{center}
  \end{figure}

In \cite{cs1}, \cite{cs2} Cappell and Shaneson  constructed an $s$--cobordism 
$H$ from $Q$ to itself as follows:  $Q$ is the union  of an $ I$--bundle over a Klein
bottle ${\bf K}$ and the solid torus $ S^{1}\times D^{2} $, glued along their
boundaries. Let $N$  be the $D^{2}$--bundle over ${\bf K}$ obtained as the open
tubular neighborhood of
${\bf K} \subset Q\times \{ 1/2 \} $ in the interior of 
$ Q\times [0,1] $. Then they constructed a certain punctured torus bundle $M$ over
${\bf K}$, with $ \partial M=\partial N $, and replaced  $N$  with $M$:  
$$ H=  M \smile( Q\times [0,1] - \mbox{interior } N) $$
They asked whether $H$ or any of its covers are trivial product
cobordisms? Evidently the $2$--fold cover of $H$ is an $s$--cobordism $\widetilde{H}$
from the lens space $L(4,1)$ to itself, and the further $4$--fold cyclic cover 
$\widetilde{\widetilde{H}}$ of
$ \widetilde{H} $ gives an $s$--cobordism from
$S^3$ to itself. For the past 15 years the hope was that this universal cover
$\widetilde{\widetilde{H}}$  might be a non-standard $s$--cobordism, inducing a fake smooth
structure on $S^4$. In this paper among other things  we will prove that this is
not the case by demonstrating the following smooth identification:

{\Thm  $\widetilde{\widetilde{H}}= S{^3} \times [0,1]$} 

\medskip

We will first describe a handlebody picture of $H$ (Figure 47).
Let $Q_{\pm}$ be the two boundary components of
$H$ each of which is homemorphic to
$Q$: $$\partial H = Q_{-}\cup Q_{+}$$
We can cap either ends of $H$ with $W$, by taking the union with $W$
along $ Q_{\pm} \approx \partial W $: 
 $$W_{\pm} = H\smile_{Q_{\pm}}W $$
There is more than one way of capping $H$ since $Q$ has nontrivial self diffeomorphisms, but it turns out
from the construction that there is a `natural' way of capping.  The reason for bringing the rational
ball $W$ into the picture while studying
$Q$ is  that philosophically  the relation of $W$ is to
$Q$ is similar to the relation of
$B^4$ to $S^3$. Unable to prove that $H$ itself is a product cobordism, we prove the
next best thing:

{\Thm $W_{-} =W$} 

\medskip

Unfortunately we are not able to find a similar proof for $W_{+}$. This is because the handlebody picture of
$H$ is highly non-symmetric (with respect to its two ends) which prevents us
adapting the above theorem to $W_{+}$. Even though, there is a  way of capping
$H$ with $W_{+}$ which gives back the standard $W$, it does not correspond to our `natural' 
way of capping (see the last paragraph of Section 1).  

The story for $W_{+}$ evolves in a completely
different way: Let $\widetilde{W}$ and $\widetilde{W}_{+}$ denote the
$2$--fold covers of $W$ and $W_{+}$ respectively (note that $ \pi_{1}(W)={\bf Z}_{2}$ and $\widetilde{W}$ 
is the Euler class
$-4$ disk bundle over $S^{2}$). We will manege to prove $\widetilde{W}_{+}$ is 
standard, by first showing that  
it splits as  $W\# \Sigma $, where $\Sigma $ a certain homotopy $4$--sphere, 
and then by proving
 $\Sigma $ is in fact diffeomorphic to $S^{4}$.

{\Thm  $\widetilde{W}_{+} =\widetilde{W}$}

\medskip

It turns out that the homotopy sphere $\Sigma $ is obtained from $S^{4}$ by the {\it
Gluck construction} along a certain remarkable
$2$--knot $A\subset S^{4}$ (ie, there is an imbedding 
$ F\co S^2 \hookrightarrow S^4 $ with $ F(S^2)=A $). Furthermore
 $A$ is the fibered knot in $S^4$ with fiber consisting of the
punctured quaternionic
$3$--manifold $Q_{0}$, with monodromy  $ \phi\co  Q_{0} \to Q_{0} $ coming from the
restriction of the order $ 3 $ diffemorphism of $Q$, which cyclically permutes
the three singular fibers of $Q$ (as Seifert Fibered space). Recall that, in
\cite{p} the Mapping Class group  $\pi_{0}(\mbox{Diff}Q$) of $Q$ was computed to
be $ {\it S}_{3} $, the symmetric group on three letters (transpositions
correspond to the diffeomorphisms pairwise permuting the three singular fibers).
In fact in a peculiar way this $2$--knot
$A$ completely determines the $s$--cobordism $H$ (this is explained in the next
paragraph). Then Theorem 3 follows by showing that the Gluck construction to 
$S^4$ along $A$ yields $S^4$, and Theorem 1 follows by showing  $\widetilde{\widetilde{H}}=
(S{^3} \times [0,1]) \# \Sigma \# \Sigma $.

Figure 77 describes a nice handlebody description of
$ W_{+}$. One remarkable thing about this figure is that it explicitly
demonstrates $H$ as the complement of the imbedded $ \;W\subset
W_{+}$ (we refer this as the {\it vertical handlebody of $H$}). Note that
$W$ is clearly visible in Figure 77. 
It turns out that $W_{+}$ is obtained by attaching a
$2$--handle to the complement of the tubular
neighborhood of $A$ in $S^4$ and the $2$--handle $ H_{1} $  is attached the
simplest possible way along the twice the meridional circle of
$A$ ! (Figures 77 and 82). Equivalently, $W_{+}$ is the complement of the tubular
neighborhood  of of an knotted ${\bf RP}^2$ in
$S^4$, which is obtained  from the standardly imbedded ${\bf RP}^2$ by 
connected summing operation ${\bf
RP}^2 \# A \subset S^4\#S^4 =S^4$. 

Put another way, if $f\co  B^2 \hookrightarrow B^4$ is the proper imbedding (with
standard boundary) induced from $A\subset S^{4}$ by deleting a small ball $B^4$
from $S^4$, then (up to $3$--handles) $W_{+}$  is obtained by removing
a tubular neighborhood of $f(B^2)$ from $B^4$  and attaching a
$2$--handle along the circle in $ S^3 $, which links the unknot $f(\partial B^2)$
twice as in Figure 1.  Note that in Figure 1 $ f(\partial B^2) $ corresponds to the
the circle with dot. So $W_{+}$ is obtained by removing a tubular neighborhood of a 
properly imbedded knotted $2$--disc from
$S^2\times B^2$, while $W$ is obtained by removing a tubular neighborhood of the unknotted
$2$--disc with the same boundary. This is very similar to the structure of the fake
fishtail of \cite{a1}. Also, it turns out that $A$ is the $3$--twist spun of the trefoil knot, and 
it turns out that
$W_{+}$ is obtained from $W$ by the  Fintushel--Stern knot surgery operation \cite{fs}.

The reason why we have not  been able to decide  whether $H$
itself is the product cobordism is that we have not been able to put the handlebody of $H$ in a
suitable form to be able to apply our old reliable ``upside-down turning trick" (eg \cite{a2},
\cite{ak2}), which is used in our proofs. This might yet happen, but until then potentially
$H$ could be a fake $s$--cobordism. 

\medskip
{\bf Acknowledgements}\qua We would like to thank R Kirby for giving
us constant encouragement and being a friendly ear during development
of this paper, and also U Meierfrankenfeld for giving us generous help
with group theory which led us to prove the crucial fibration theorem.
We also want to thank IAS for providing a nice environment where the
bulk of this work was done.  The author was partially supported by NSF
grant DMS 9971440 and by IAS.

\section{Handlebody of $Q\times [0,1]$}

Let $I=[0,1]$. We draw $Q\times I$ by using the technique of \cite{a1}. $Q$ is
the obtained by surgering $ S^{3} $ along a link $L$ of two components linking each
other twice (the first picture of Figure 1), hence $ Q\times I $ is obtained by
attaching two $2$--handles  to $ (S^{3}-L) \times I $. 

Note that, since
$3$ and $4$--handles  of any $4$--manifold are attached in a canonical
way, we only need to visualize the $1$ and $2$--handles of $Q\times I$. Hence it
suffices to visualize $ (B^{3} -L_{0})\times I $,
where $ L_{0} $ is a pair of properly imbedded arcs linking each other twice,  plus
the two $2$--handles as shown in Figure 2 (the rest are $3$--handles).

\begin{figure}[ht!]
  \begin{center}
     \includegraphics{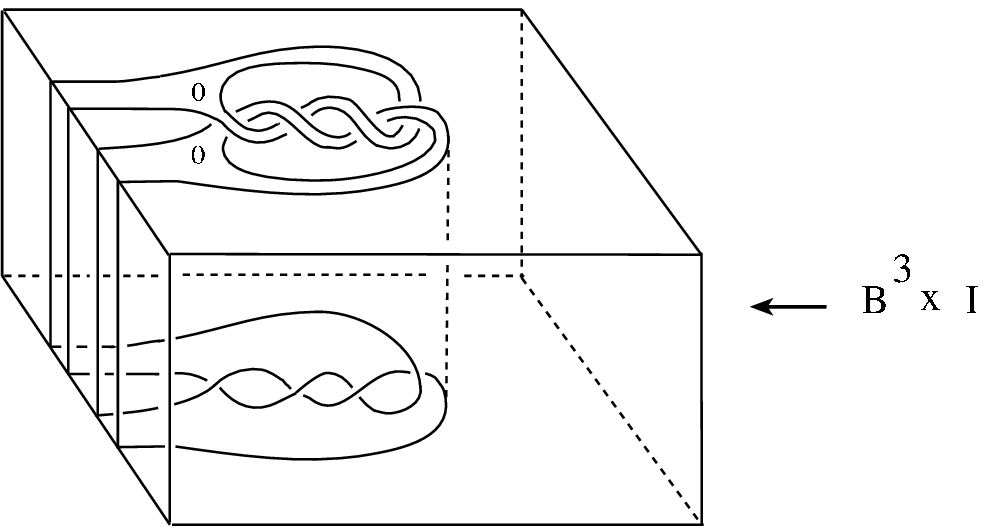}
   \caption{} \label{B2}
    \end{center}
  \end{figure}

Clearly Figure 2 is obtained first by removing the two obvious $2$--disks from
$ B^{4}=B^{3}\times I $ which $ L\#(-L) $ bounds, and then by attaching two
$2$--handles (here $-L$ denotes the the mirror image of $L$). This gives the 
first picture of Figure 3. In Figure 3 each circle with dot denotes a $1$--handle
(ie, the obvious disks it bounds is removed from $B^{4}$). The
second picture of Figure 3 is diffeomorphic to the first one, it is obtained by
sliding a
$2$--handle over a
$1$--handle as indicated in the figure. By an isotopy  of Figure 3 (pulling 1--handles
apart) we obtain the first picture of Figure 4, which is the same as
the second picture, where the $1$--handles are denoted by a different notation (as
pair of attaching balls). Hence Figure 4 gives $Q\times I$. 

\begin{figure}[ht!]
  \begin{center}
     \includegraphics{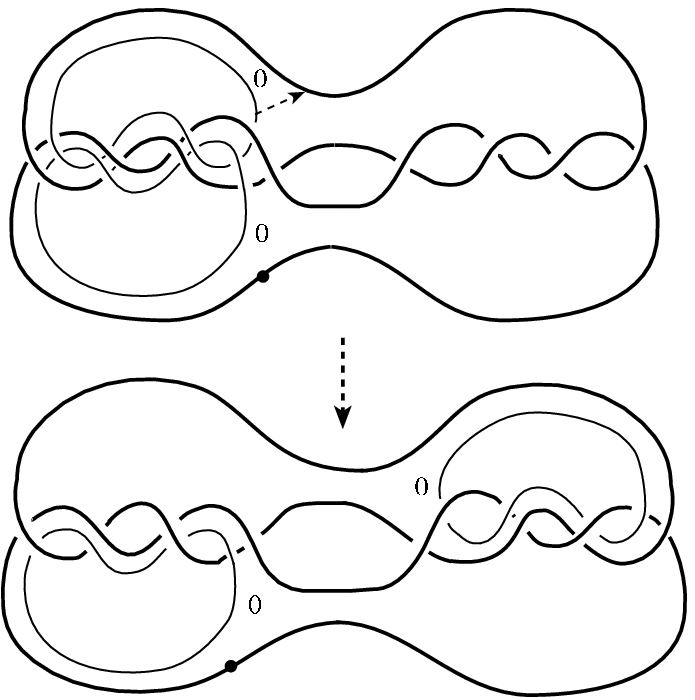}
   \caption{} \label{B3}
    \end{center}
  \end{figure}

\begin{figure}[ht!]
  \begin{center}
     \includegraphics{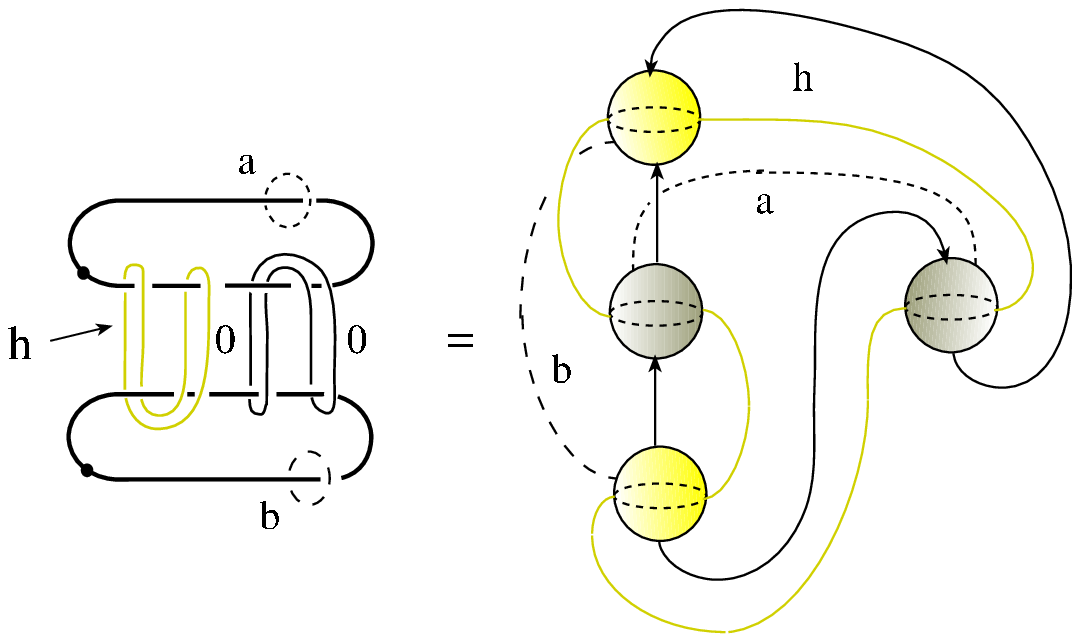}
   \caption{} \label{B4}
    \end{center}
  \end{figure}

For a future reference the
linking loops $a,b$ of the $1$--handles of the first picture of Figure 4 are indicated
in the second picture of Figure 4.  Notice that we can
easily see an imbedded copy of the Klein bottle ${\bf K}$ in $Q\times I $ as follows:
The first picture of Figure 5 denotes  ${\bf K}$ (a square with the opposite
sides identified as indicated). By thickening this to a four dimensional
handlebody we obtain the second picture of Figure 5 which is a
$D^2$--bundle $N$ over ${\bf K}$ (the orientation
reversing $1$--handle is indicated by putting ``tilde" in the corresponding balls). 
The third and the fourth pictures of Figure 5 are also $N$ , drawn in different
$1$--handle notations.

Clearly the handlebody of Figure 5 sits in Figure 4, demonstrating an imbedding of
the disk bundle $N$ over the Klein bottle into $Q\times I$. For the purpose of
future references, we indicated where the linking circle $b$ of the
$1$--handle lies in the various handlebody pictures of $N$ in Figure 5 .

Finally in Figure 6 we draw a very useful `vertical' picture of $Q\times I$  as
a product cobordism starting from the boundary of $W$ to itself like a collar. 
Though this
is a seemingly a trivial handlebody of $H$ it will be useful in a later construction.
Later, we will first construct a handlebody picture of the
 $s$--cobordism
$H$ from $Q$ to itself, and then view it like a collar sitting on the boundary of $W$,
ie, as a vertical picture of the cobordism starting from the boundary of a $W$ to $Q$.
Note that  $N$ is clearly visible in Figure 6, which is an alternative handlebody of $W$
(N is lying in the collar of its boundary). Also notice that the operation Figure 4
$\Rightarrow$  Figure 6,  ie, capping one end of $Q\times I$ by $W$, corresponds
to attaching a
$2$--handle to Figure 4 along the loop $b$. Similarly capping the other end of
$Q\times I$ by $-W$ corresponds to attaching a $2$--handle to it along the loop $a$.
Note that Figure 4 can also indicate the handlebody of $Q_{0}\times I$, where
$Q_{0}$ is the punctured $Q$ (in this case we simply ignore the $3$--handle).

\begin{figure}[ht!]
  \begin{center}
     \includegraphics[width=\hsize]{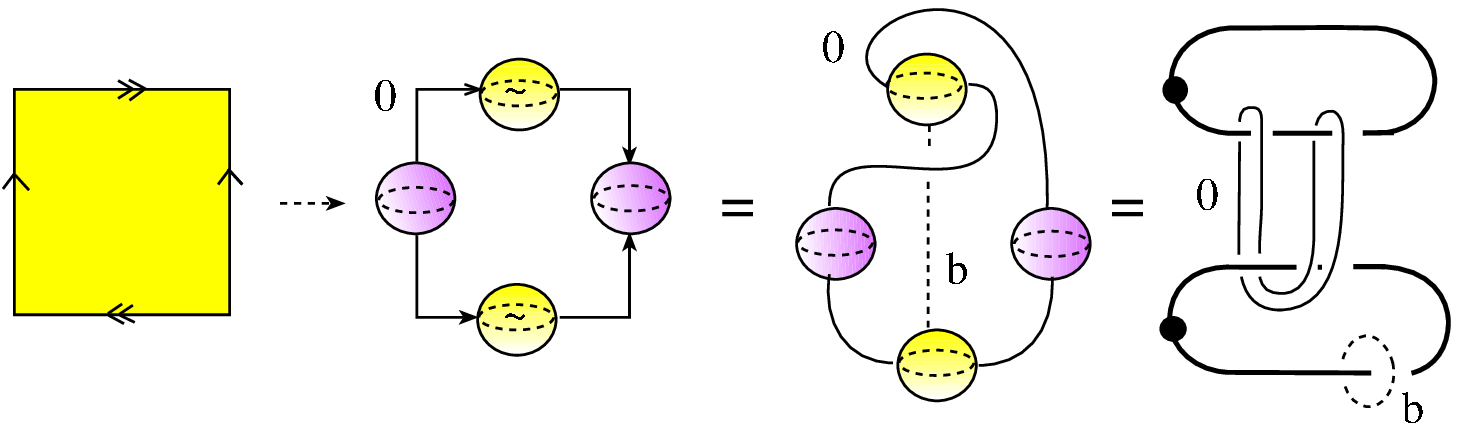}
   \caption{} \label{B5}
    \end{center}
  \end{figure}

\begin{figure}[ht!]
  \begin{center}
     \includegraphics[width=\hsize]{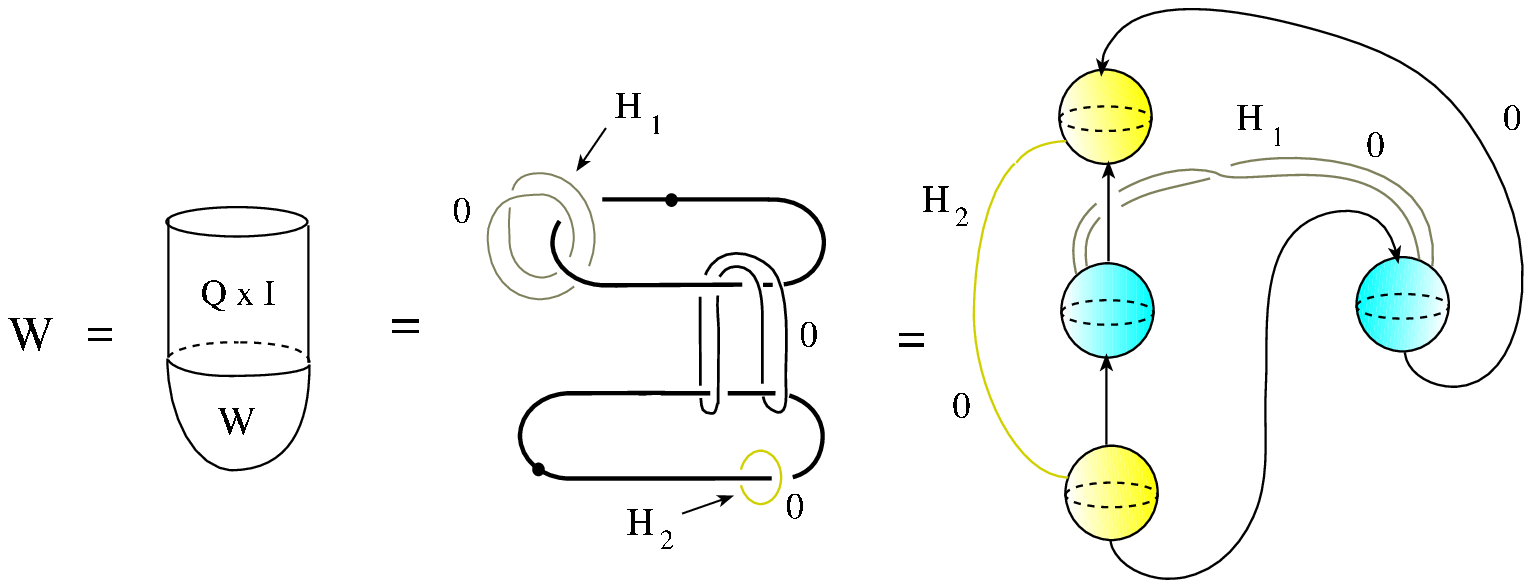}
   \caption{} \label{B6}
    \end{center}
  \end{figure}

\section{Construction of $H$ }

Here we will briefly recall the Cappell--Shaneson construction \cite{cs1}, 
and indicate why $H$ is an
$s$--cobordism: Let
$T_{0}$ denote the punctured
$2$--torus.
$M$ is constructed by gluing together two
$T_{0}$--bundles over Mobius bands given with the monodromies: 
$$ A=\left( \begin{array}{cc}0&1\\ 1&1\end{array}\right) \;\;\;\;  
B=\left( \begin{array}{cc}0&-1\\ -1&-1\end{array}\right)  $$
Since $ A^{2}=B^{2} $ these bundles agree over the boundaries of the Mobius bands,
hence they give a bundle $M$ over the the union of the two
Mobius bands along their boundaries (which is the Klein bottle). By using the handle description of
${\bf K}$ given by the second picture of Figure 7, we see that $M$ is the
$T_{0}$--bundle over ${\bf K}$, defined by the monodromies 
$$ A= \left( \begin{array}{cc}0&1\\ 1&1\end{array}\right) \;\;\;\;  C=B^{-1}A
=\left( \begin{array}{cc}-1&0\\ 0&-1\end{array}\right)  $$
\begin{figure}[ht!]
  \begin{center}
     \includegraphics{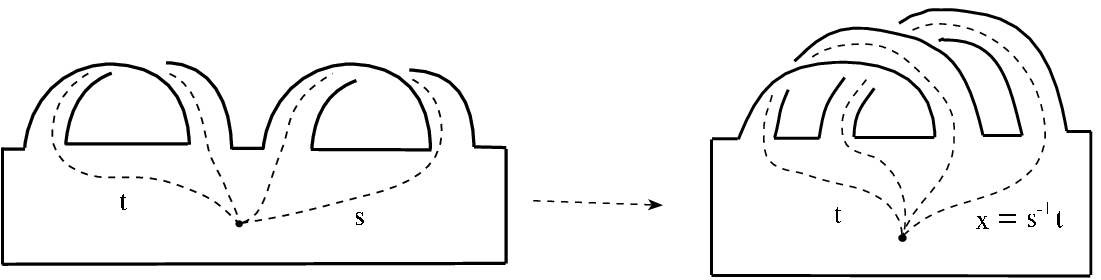}
   \caption{} \label{B7}
    \end{center}
  \end{figure}
\noindent  Let $t, x$ and
$\tau, \xi$ to be the standard generators of the fundamental groups of ${\bf K}$  and
$T_{0}$ respectively , then:
$$ \pi_{1} (M)= \left< t, x, \tau, \xi \; \;
| \; \; 
\begin{array}{c}  t \tau t^{-1}=\xi \;  \; \;
\;  t\xi t^{-1}=\tau \xi , \\ x\tau x^{-1}=
\tau^{-1}, \; x \xi x^{-1}=\xi ^{-1} \\ txt^{-1}=x^{-1}
\end{array}  \right> $$ 
So $\; x^{2}\tau x^{-2}=x\tau^{-1} x^{-1}=\tau $, 
$\;x^{2}\xi x^{-2}= x\xi^{-1}x^{-1}=\xi\; \Rightarrow$  
$x^{2}\tau =\tau x^{2} \;$  and $\; x^{2}\xi=\xi x^{2} $.

Recall that $\;H=  M \smile( Q\times [0,1] - \mbox{interior } N)\;$, which is the same
as  $M \smile h $, where $h$ is a $ 2 $--handle (Figure 4).
Let us briefly indicated why the boundary inclusion induces an isomorphism $
\pi_{1} (Q) \to
\pi_{1}(H) $: By Van-Kampen theorem attaching the $2$--handle $h$ introduces the
relation
$ xtx^{-1}=t^{-1} $ to $ \pi_{1}(M) $;  which together with
$txt^{-1}=x^{-1}$ gives $t^{-2}=x^{2}$. Therefore in
$\pi_{1}(H)$ the relations 
$\; t^{2}\tau t^{-2}=t\xi t^{-1}
$, and
$ t^{2}\xi t^{-2} =t\tau \xi t^{-1}$ become:  $\; \tau =t\xi t^{-1} $, and
$\;\xi  =t\tau \xi t^{-1}\;$. Then by substituting $\tau $ in $\xi$, and by
using the fact that $\xi$ commutes with $t^{-2}$,  we get $\xi=1$ and hence
$\tau=1$. Hence the boundary inclusions $\pi_{1} (Q)\to \pi_{1}(H) $ induce isomorphisms. In fact
$H$ is an
$s$--cobordism from
$Q$ to itself.  From now on let 
$ \hat{M} $ denote the corresponding $T^{2}$--bundle over ${\bf K}$ 
induced by $M$ by the obvious way,
clearly: $$M = \hat{M} - N $$

\section{Handlebody of $\hat{M}$ }

 From the last section we see that $M$ is obtained by first taking the
$T_{0}$--bundle $ \;T_{0}\times_{-I}S^{1} \to S^{1}$ with monodromy $-Id$ , and
crossing it by
$I$, and by identifying the ends of this $4$--manifold with the monodromy:
$$\varphi=\left(
\begin{array}{ccc}0&1&0\\ 1&1&0\\ 0&0&-1\end{array}\right) $$
This is indicated in Figure 8 (here we are viewing ${\bf K}$ as the handlebody of
the second picture of Figure 7). 

\begin{figure}[ht!]
  \begin{center}
     \includegraphics[width=4in]{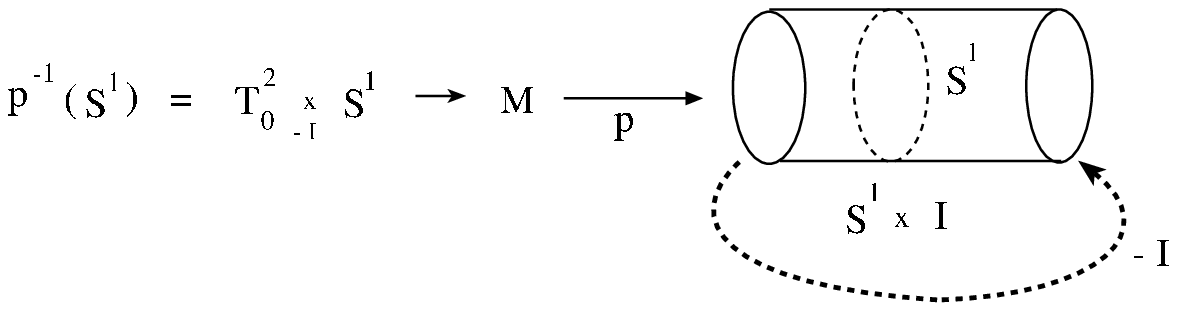}
   \caption{} \label{B8}
    \end{center}
  \end{figure}

Drawing the handlebody of $M= \;(T_{0}\times_{-I}S^{1})\times_{\varphi} S^{1} $
directly will present a difficulty in later steps, instead we will first draw the
corresponding larger $T^{2}$--bundle 
$$\hat{M}=\;(T^{2}\times_{-I}S^{1})\times_{\varphi} S^{1}  $$ 
then remove a  copy of $N$ from it. $ T^{2}\times_{-I}S^{1}$ is obtained by
identifying the two ends of $ T^{2}\times I $ with
$-Id$. Figure 9 describes a two equivalent pictures of the 
Heegaard handlebody of
$ T^{2}\times_{-I}S^{1} $. The pair of `tilde' disks 
describe a twisted $1$--handle
(due to $-Id $ identification  $ (x,y)\to (-x,-y) $). If need be, after
rotating the attaching map of the the twisted $1$--handle we can turn it to a
regular $1$--handle as indicated by the second picture of Figure 9.

\begin{figure}[ht!]
  \begin{center}
     \includegraphics[width=\hsize]{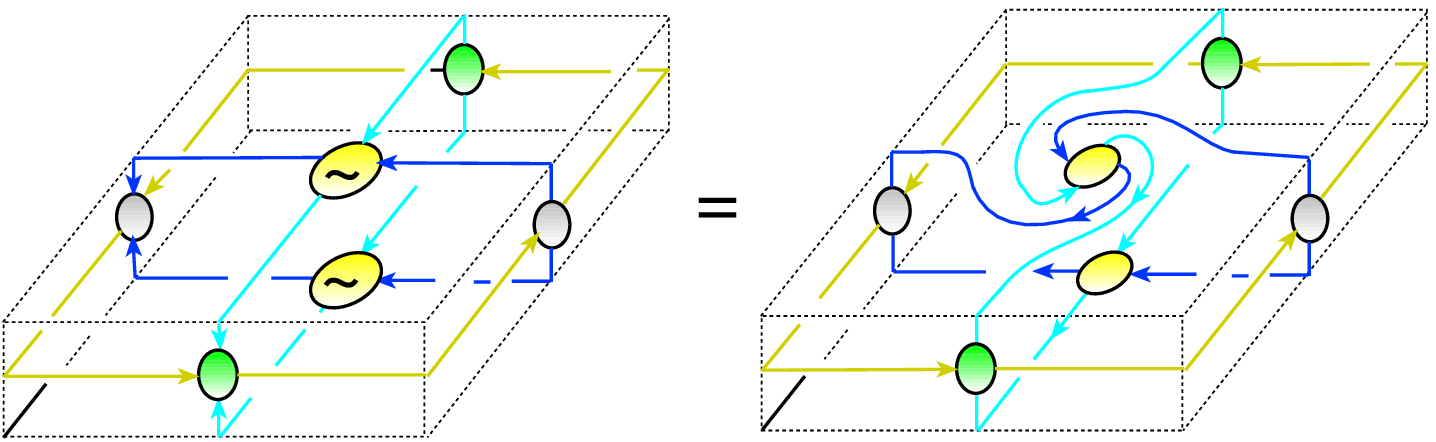}
   \caption{} \label{B9}
    \end{center}
  \end{figure}

Now we will draw $\hat{M} = (T^{2}\times_{-I}S^{1})\times_{\varphi}S^{1} $  
by using the technique introduced
in
\cite{ak}: We first thicken the handlebody 
$ T^{2}\times_{-I}S^{1} $ of Figure 9  to the $4$--manifold $ T^{2}\times_{-I}S^{1}\times I
$  (the first picture of Figure 10). Then
isotop   $\varphi \co  T^{2}\times_{-I}S^{1} \to T^{2}\times_{-I}S^{1}  $ so that it takes $1$--handles to
$1$--handles, with an isotopy, eg,
 $$ \varphi_{t}=\left( \begin{array}{ccc}-t&1-t&0\\ 1-t&1&0\\ 0&0&-1\end{array}\right) $$
Then  attach an extra $1$--handle and the $2$--handles as indicated in second 
picture of Figure 10 (one of the  
attaching balls of the new  $1$--handle  is not visible in the picture
since it  is placed at the point of infinity). The extra $2$--handles are induced from the
identification of the $1$--handles of  the two boundary components of $ \;
T^{2}\times_{-I}S^{1}\times I
\;$  via $\varphi$. So, the second picture of Figure 10 gives the handlebody of
$\hat{M}$.

\begin{figure}[ht!]
  \begin{center}
     \includegraphics[width=\hsize]{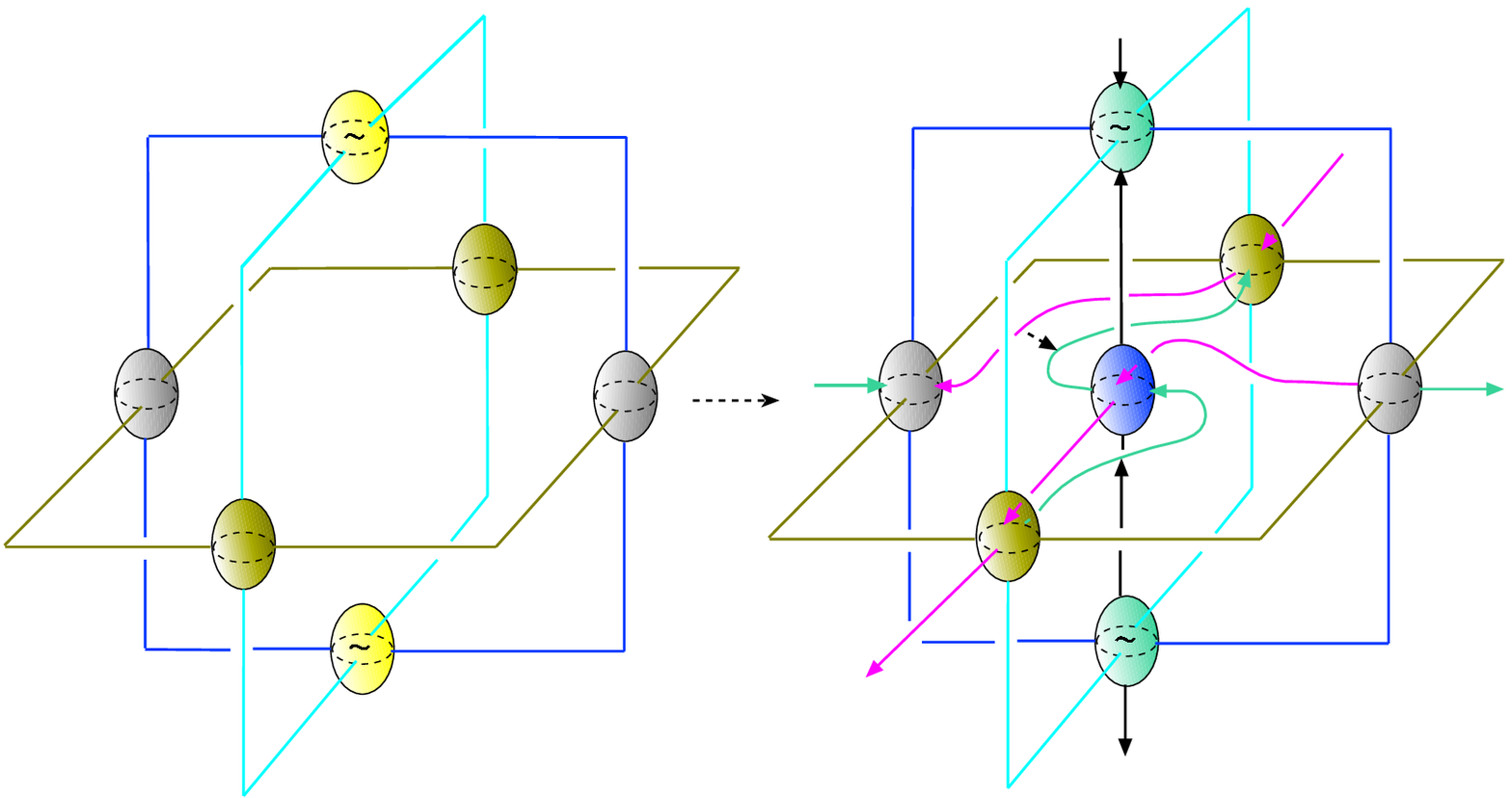}
   \caption{} \label{B10}
    \end{center}
  \end{figure}
\noindent

We want to emphasize that
the new $1$--handle identifies  the $3$--ball at the center of Figure 10 with the
$3$--ball at the infinity by the following diffeomorphism as indicated in Figure 11: 
$$ (x,y,z) \to (x,-y,-z)$$ 
\begin{figure}[ht!]
  \begin{center}
     \includegraphics{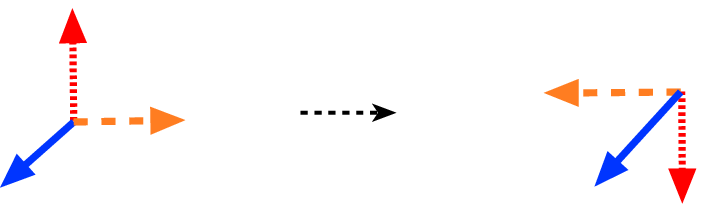}
   \caption{} \label{B11}
    \end{center}
  \end{figure}%
Figure 12 describes how part of this isotopy $\varphi_{t}$ acts on $T^{2}$ 
(where $T^{2}$ is represented by a disk with opposite sides identified).
This is exactly the reason why we started with $T^{2}\times_{-I}S^{1} $ instead of
$T^{2}_{0}\times_{-I}S^{1} $ (this isotopy 
takes place in $T^{2}$ not in $T^{2}_{0} $ !).

\begin{figure}[ht!]
  \begin{center}
     \includegraphics[width=\hsize]{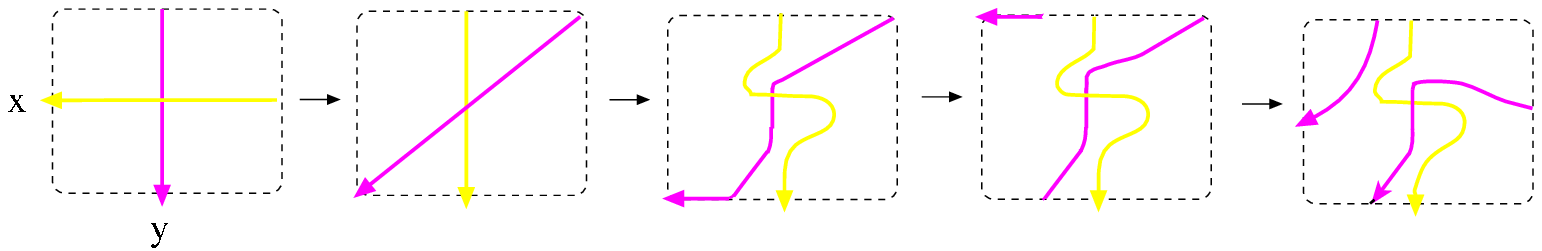}
   \caption{} \label{B12}
    \end{center}
  \end{figure}

\section{Simplifying the handlebody of $\hat{M}$ }

We now want to simplify the handlebody of Figure 10 by cancelling some
$1$-- and 2--handle pairs and by isotopies: We first perform the
$2$--handle slide as indicated (by the short arrow) in Figure 10 and
obtain Figure 13. By doing the further $2$--handle slides as indicated
in Figures 13--15 we obtain Figure 16.  Note that while going from
Figure 15 to 16 we cancelled a $2$--handle with a $3$--handle (ie, we
erased a zero framed unknotted circle from the picture).

Figure 17 is the same as Figure 16 except that the
twisted $1$--handle (two balls with `tilde' on it) is drawn in the standard way. By an isotopy
we go from Figure 17 to  Figure 18. Figure 19 is the same as 
Figure 18 except that we drew one of the $1$--handles in a different $1$--handle notation
(circle with a dot notation). 

Note that in our figures, if the framing of a framed knot is
the obvious ``black-board framing'' we don't bother to indicate it, 
but if the framing deviates from the obvious black-board framing we indicate 
the deviation from to the black-board framing by putting a number 
in a circle on the knot ($-1$'s in the case of Figure 19).

Figure 20 is obtained from Figure 19 by simply
leaving out one of the $2$ handles. This is because the framed knot corresponding 
to this $2$--handle is the unknot with $0$--framing !, hence it is cancelled by a
$3$--handle (this knot is in fact  the `horizontal' framed knot of Figure 10). 

Figure 21 is the desired handlebody of $ \hat{M} $, it is the same as the
Figure 20, except that one of the attaching balls of a $1$--handle which had been
placed at the point of infinity is isotoped into  ${\bf R}^{3}$.

\section{Checking that the boundary of $\hat{M}$ is correct} 

Now we need to check that the boundary of the closed manifold $ \hat{M} $ 
 (minus the
three and four handles) is correct. That is, the boundary of Figure 21 is the  
connected sum of copies of
$ S^{1}\times S^{2} $; so that after cancelling them with
$3$--handles we get
$S^{3}$, which is then capped by a $4$-- handle. This process is done
by changing the interior of $\hat {M} $ so that boundary becomes visible: By
changing a $1$--handle to a
$2$--handle in Figure 21 (ie, turning a `dotted circle' to a zero framed
circle) we obtain Figure 22. Then by doing the indicated handle slides and 
isotopies we arrive to Figures 23, 24 and 25. Then by operation of turning a
$2$--handle to a $1$--handle by a surgery (ie, turning a zero framed circle to a
`dotted circle') and cancelling the resulting $1$-- and $2$--handle pair we get
Figure 26. By isotopies  we obtain the
Figure 28, which after surgering the obvious $2$--handle becomes
 $ S^{1}\times B^{3} \;\# \;  S^{1}\times B^{3} $ with the desired boundary.

\eject

\hbox{}\par\includegraphics[width=\hsize]{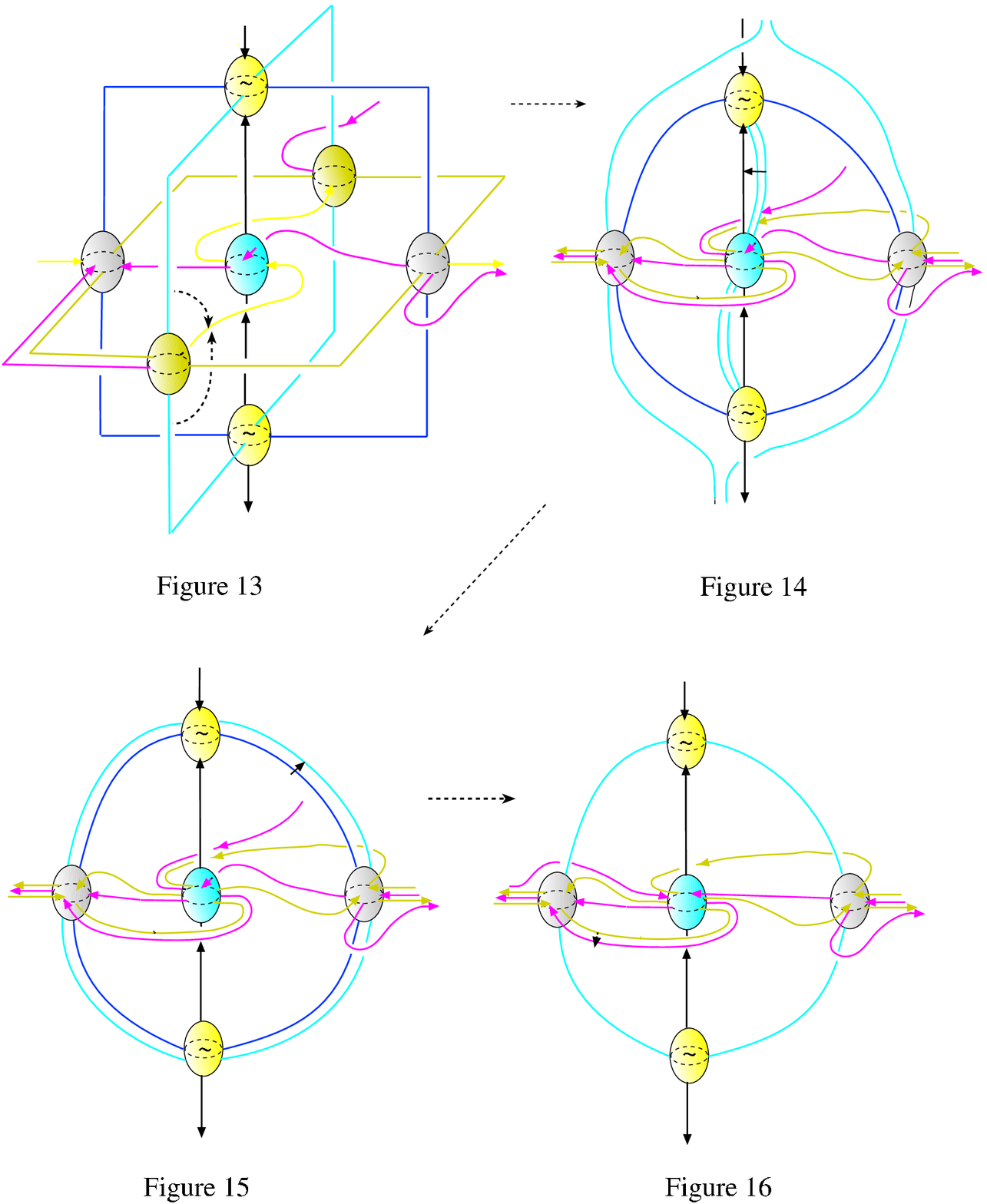}

\newpage
\hbox{}\par\includegraphics[width=\hsize]{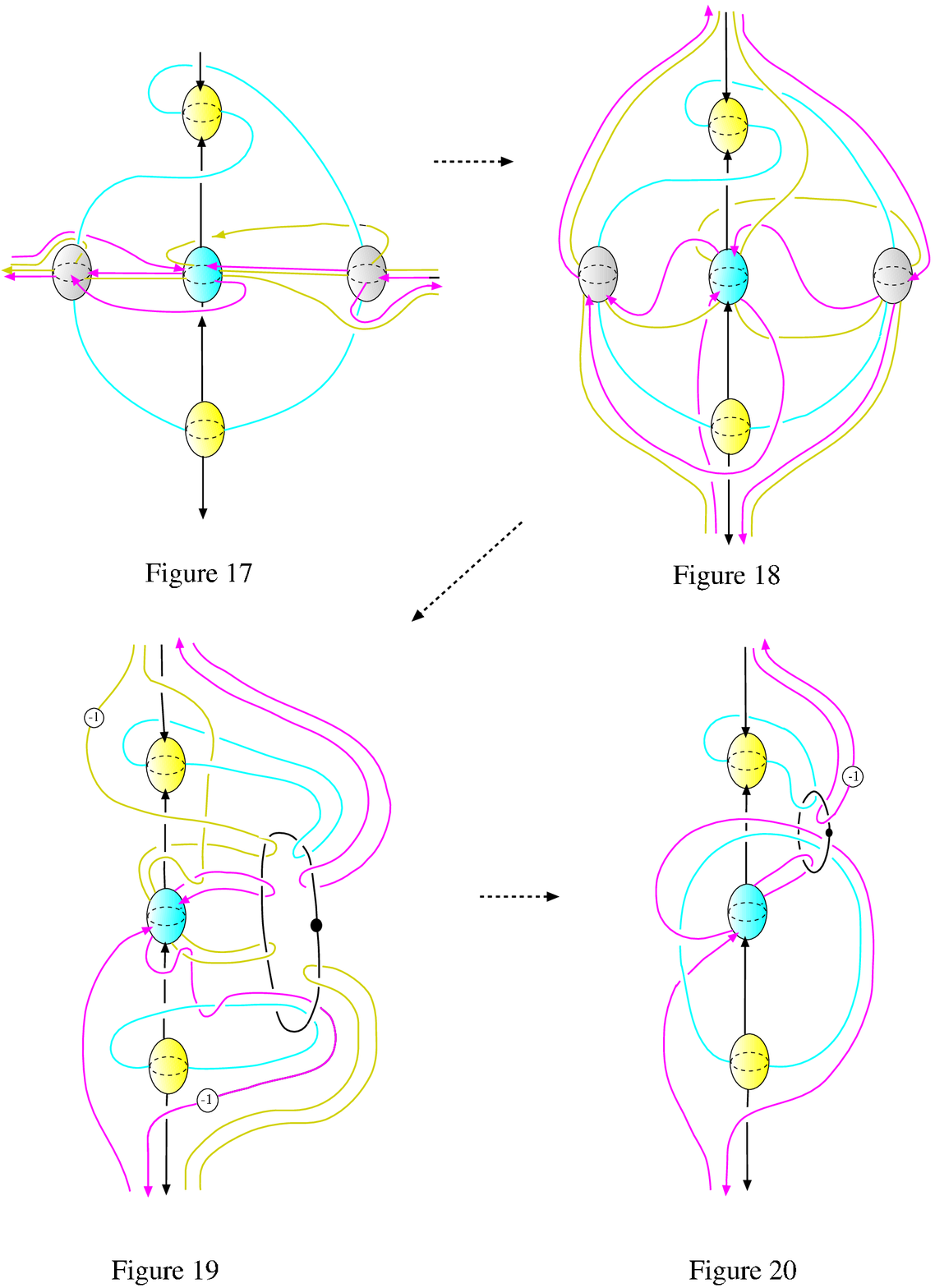}

\newpage
\hbox{}\par\includegraphics[width=\hsize]{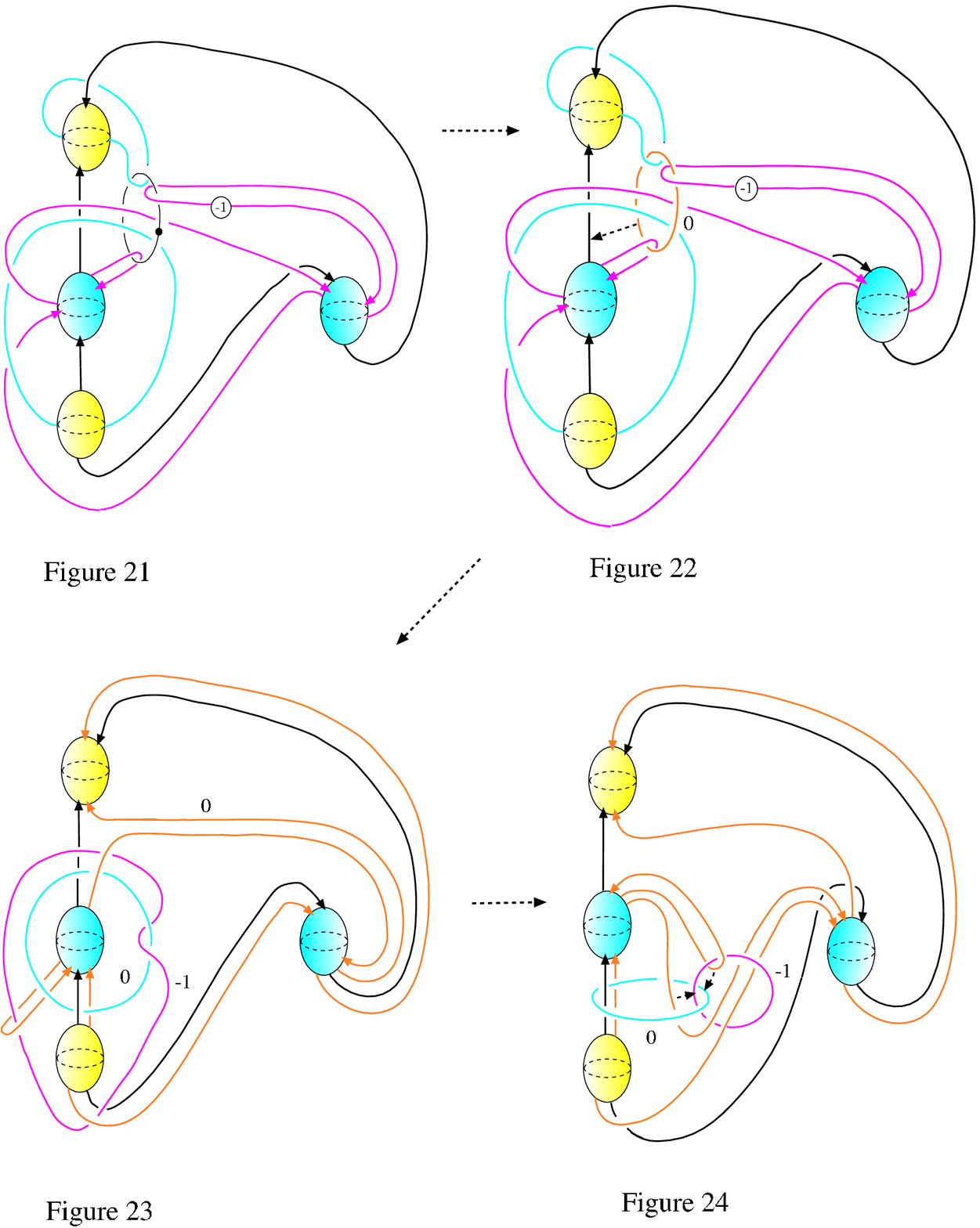}

\newpage
\hbox{}\par\includegraphics[width=\hsize]{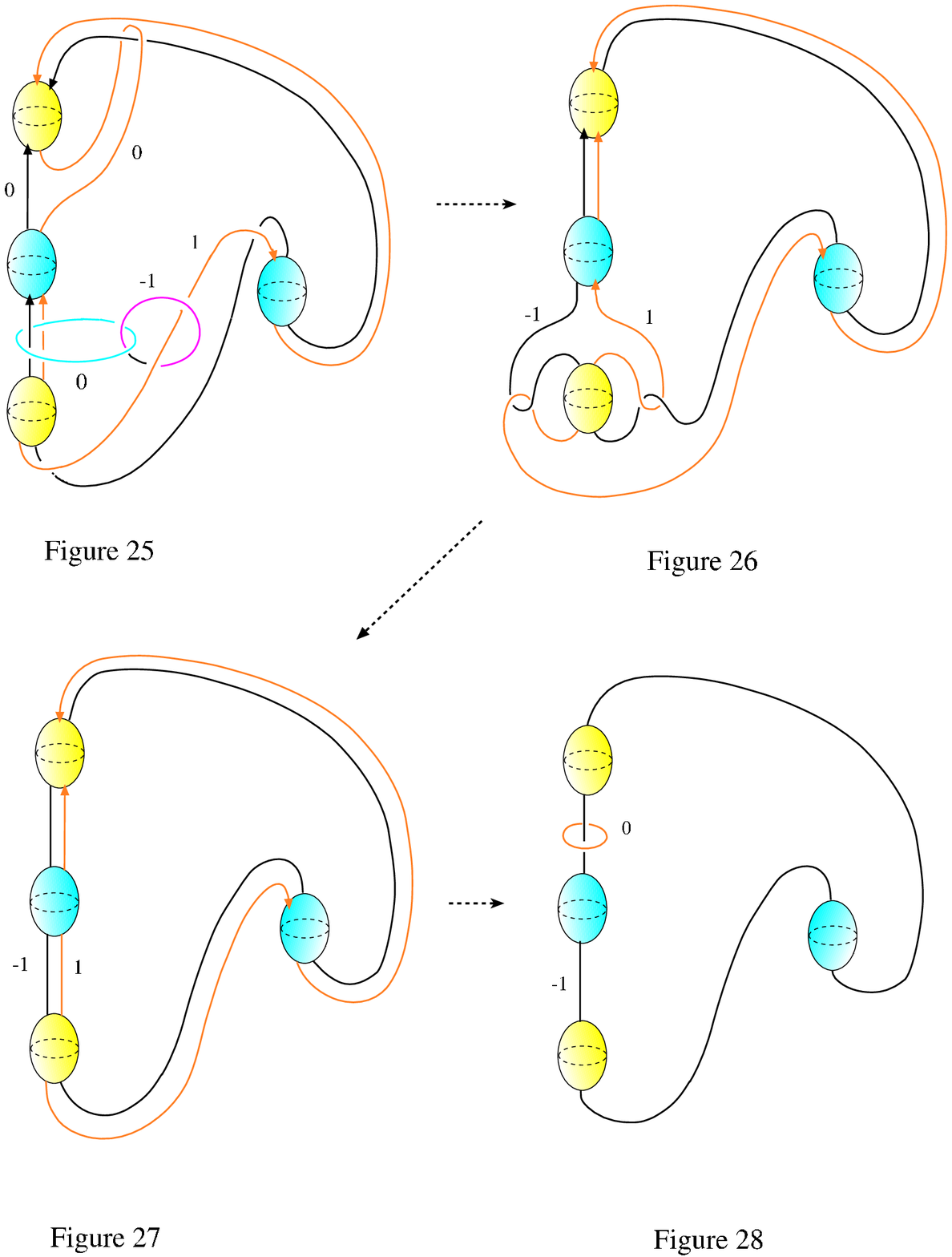}

\newpage

\section{Turning $\hat{M}$ upside down and constructing a handlebody of $H$ }

 Evidently it is not so easy to obtain the
handlebody of
$M$ from the closed manifold $\hat{M} $ even though $ M= \hat{M} - N $ and  $N$ 
is clearly visible
as a subset in the  handlebody of $\hat{M} $. For this, we will turn the
handlebody of
$ \hat{M} $ upside-down and take all the handles up to $N$ (excluding $N$), and
then attach a
$2$--handle $h$ as indicated in Figures 30 and 31.  
Recall that this last $2$--handle $h$ is attached
along the loop $h$ of Figure 4, ie, $h$  is the loop on $\partial N $ along which 
attaching a $2$--handle to $N$ gives  $ Q\times [0,1] $ (here we are using the same
notation
$h$ for the $2$ handle and for its attaching circle).

 To turn $ \hat{M} $  upside down we simply take the dual loops (attaching loops of
the $2$--handles of Figure 21, indicated by the small
$0$--framed circles in Figure 31) and the loop h, and then 
trace them via the diffeomorphism from the
boundary of Figure 21 to the boundary of Figure 28, which is 
$(S^{1}\times S^{2})\# (S^{1}\times S^{2})
$ (ie, the steps Figure 21
$\leadsto $ Figure 28); and then attach  $2$--handles to $(S^{1}\times B^{3})\#
(S^{1}\times B^{3})$ along the image of these loops. Note that, during this process
we are free to isotope these dual loops over the other handles. The reader unfamiliar
with this process can consult \cite{a2}.

Figure 31 $\leadsto $ Figure 38 is the same as the isotopy Figure 21 
$\leadsto$  Figure 28, except that we
carry the dual loops along and isotope them over the handles 
as indicated by the short arrows in these
figures. 

More explanation: By isotoping the dual loops as indicated in Figure 32 
we arrive to
Figure 33.  The move Figure 33
$\leadsto $ Figure 34 is the same as Figure Figure 25  $\leadsto$ Figure 26. 
The move Figure 34 $\leadsto $
Figure 35 is just an isotopy (rotating the lover ball by $360^{0}$ around the $y$--axis). 
By performing the handle slides as indicated by the arrows in the figures we obtain  Figure 35 
$\leadsto$ Figure 38. By changing the $1$--handle notation we obtain Figure 39, 
by rotating
the upper attaching ball of the $1$--handle  by $360^{0}$ we obtain Figure 40. 
Then by a handle
slide (indicated by the arrow) we obtain Figure 41. Changing the notation of the  
remaining $1$--handles to `circles-with-dots' we get Figure 42. Then by the indicated
handle slide we get Figure 43,  which is the same as Figure 44 (after an isotopy). The
indicated handle slides gives the steps Figure 44 $\leadsto$ Figure 47. Figure 47 is our desired
handlebody picture of $H$. The reader is suggested to compare this picture with the
picture of $Q_{0}\times I$ in Figure 4. Here we also traced the position of the loop
$a$ lying on the boundary of $Q_{0}\times I$. 

\eject

\hbox{}\par\includegraphics[width=\hsize]{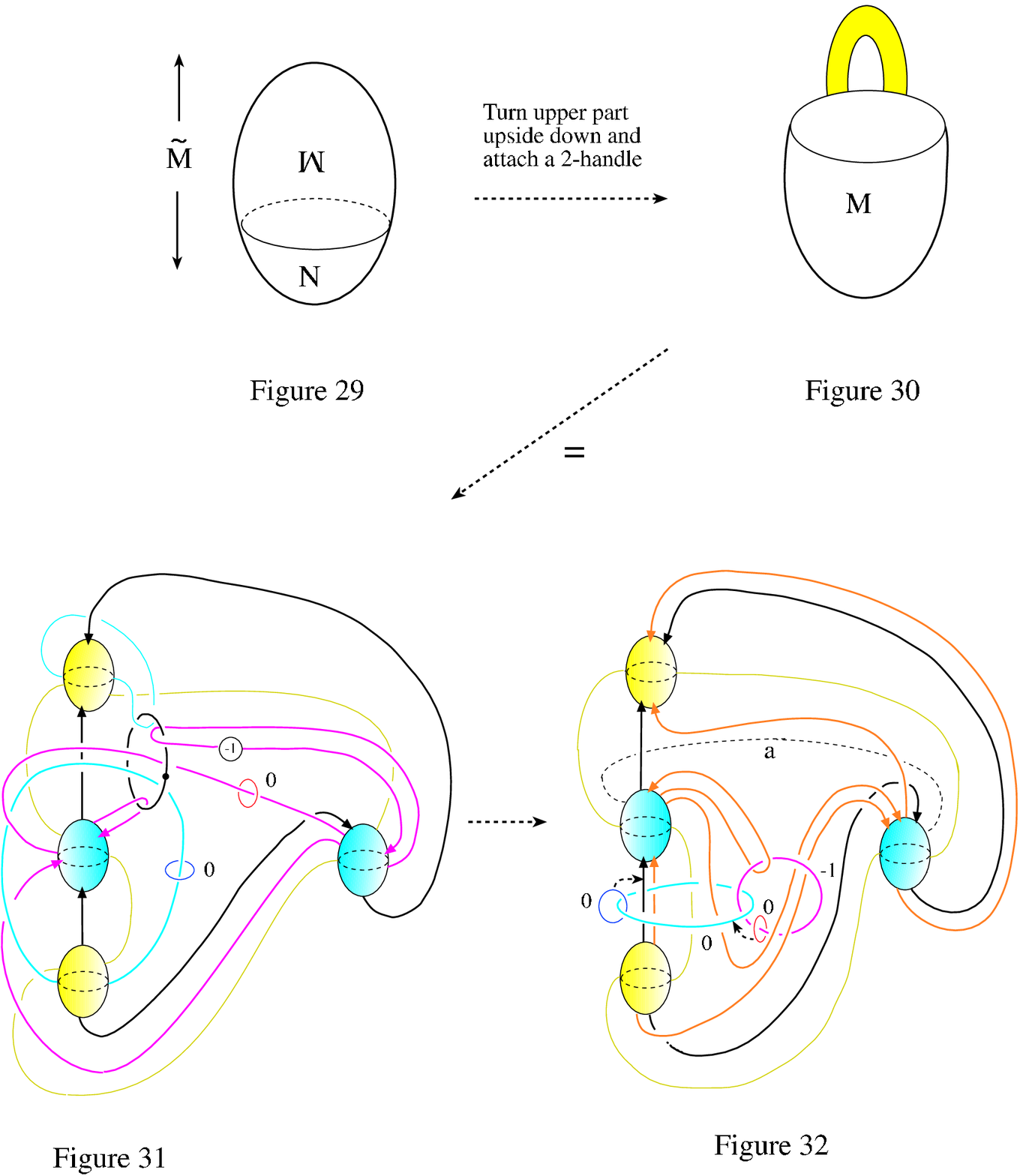}

\newpage
\hbox{}\par\includegraphics[width=\hsize]{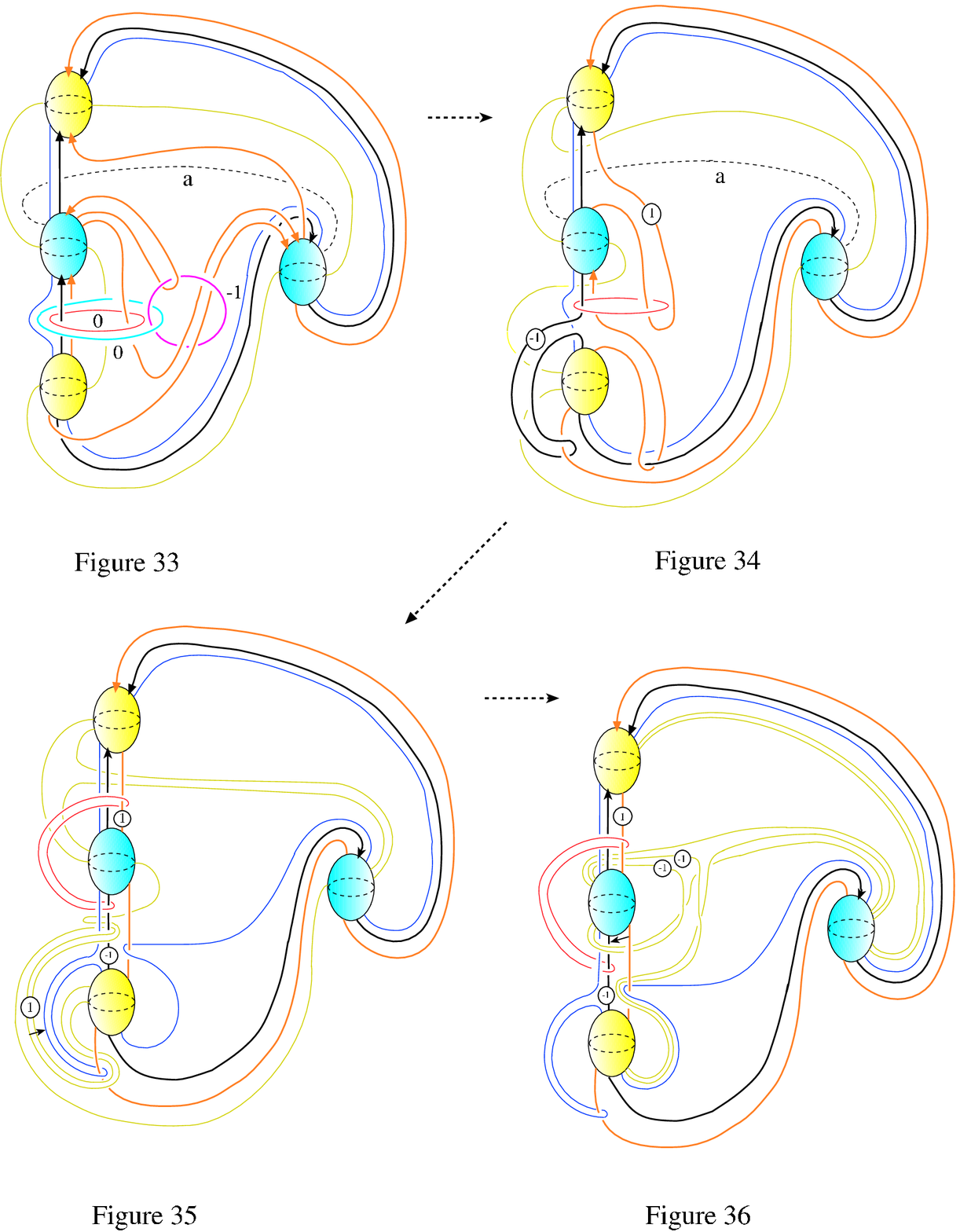}

\newpage
\hbox{}\par\includegraphics[width=\hsize]{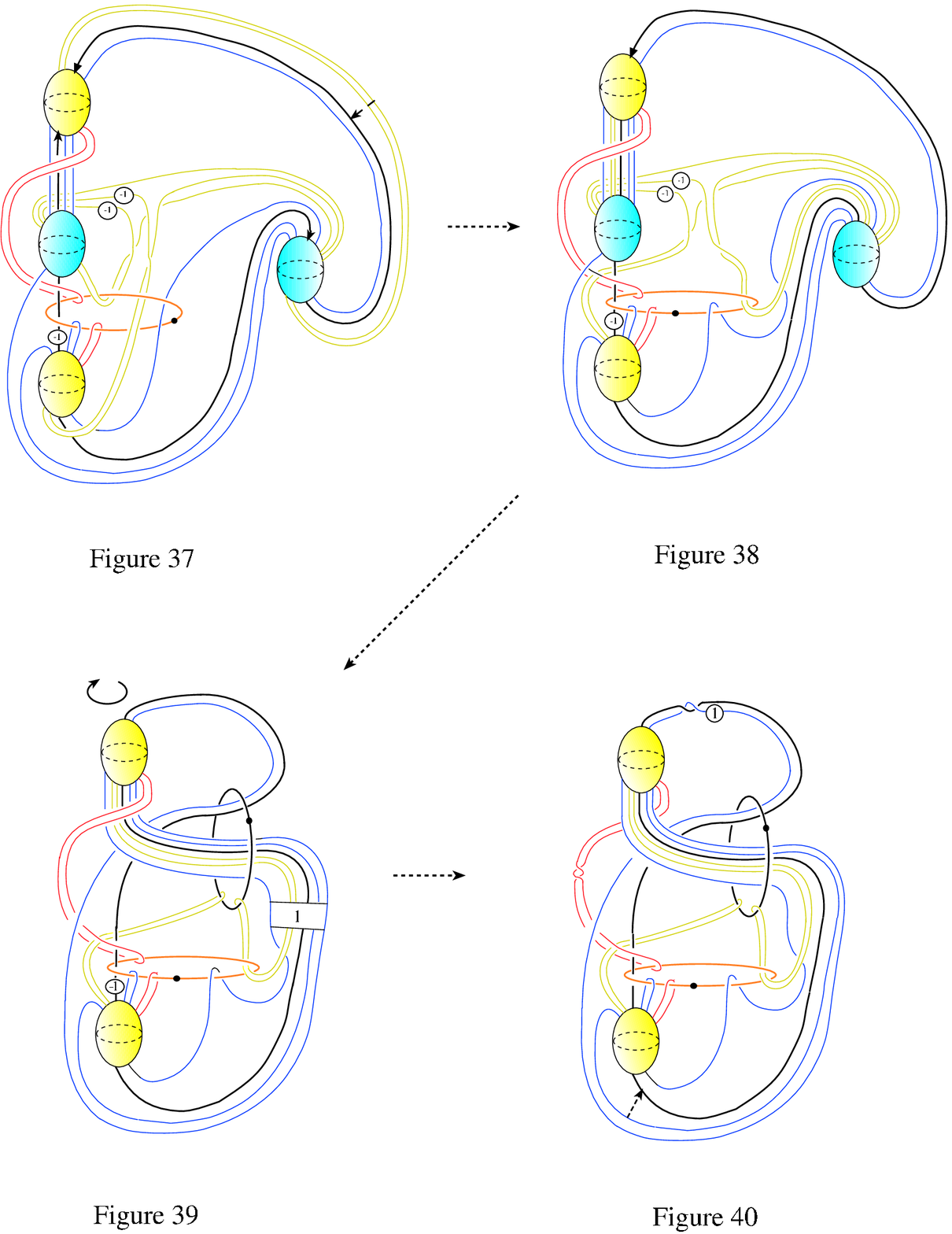}

\newpage
\cl{\includegraphics[width=4in]{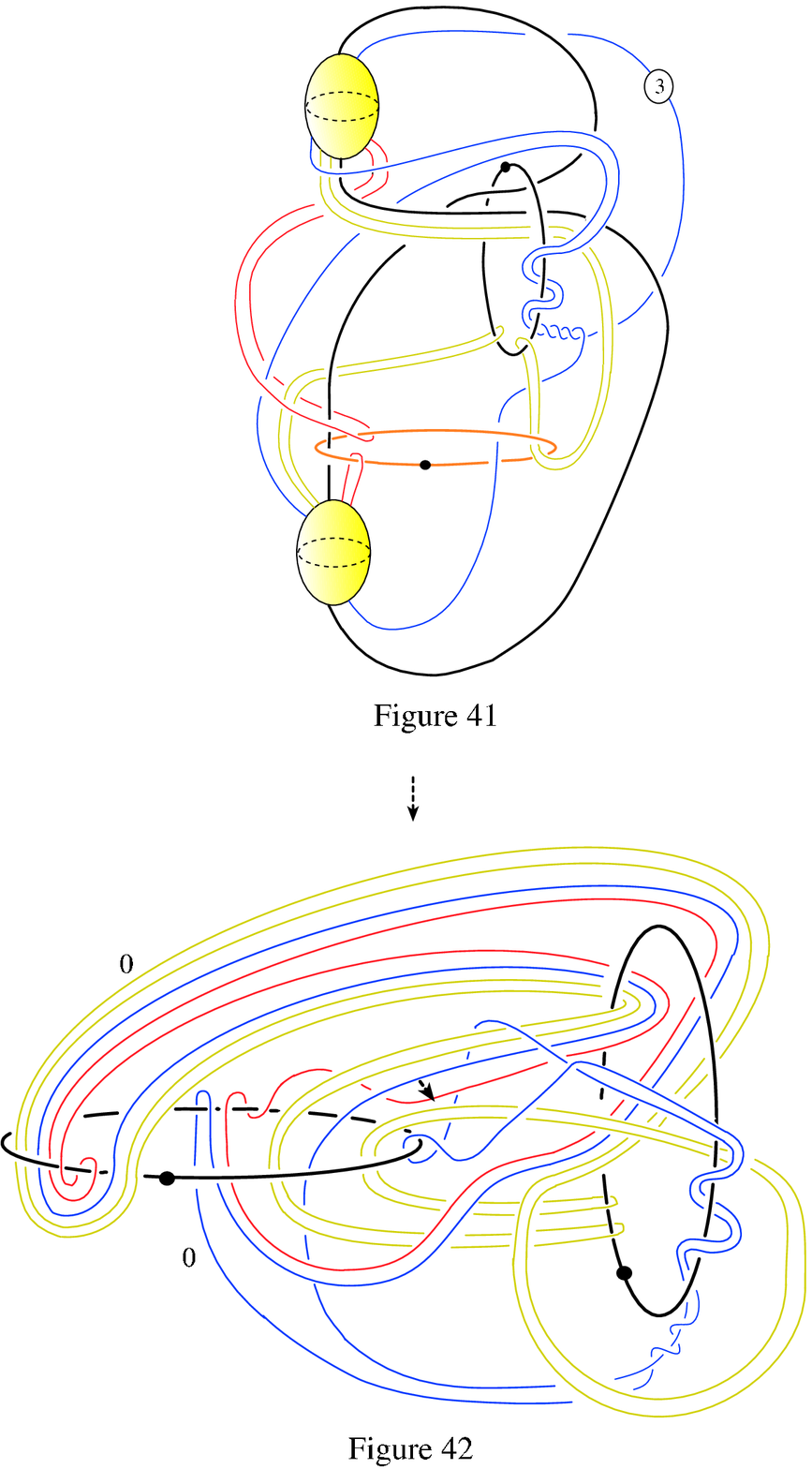}}

\newpage
\cl{\includegraphics[width=4.5in]{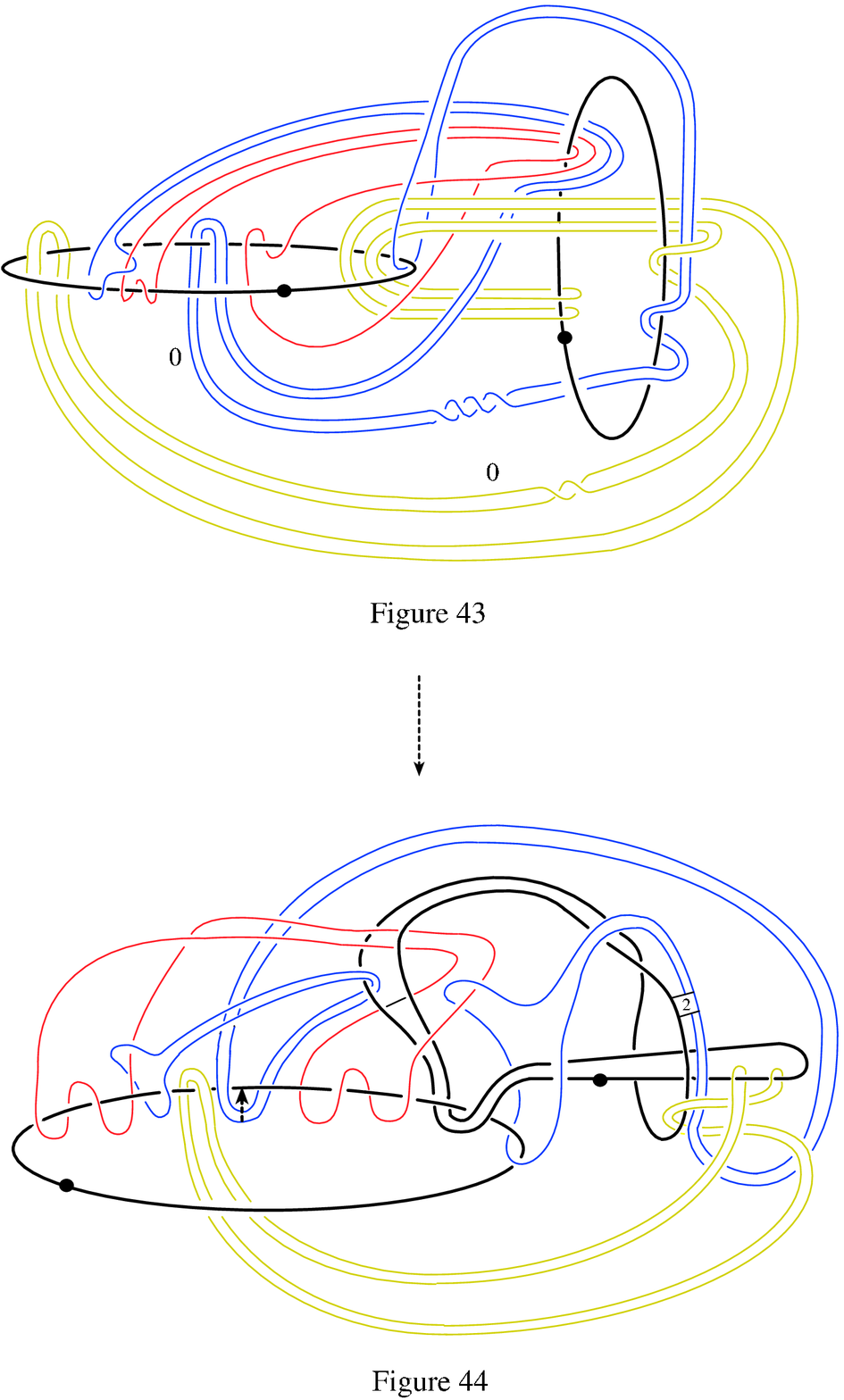}}

\newpage
\cl{\includegraphics[width=4.8in]{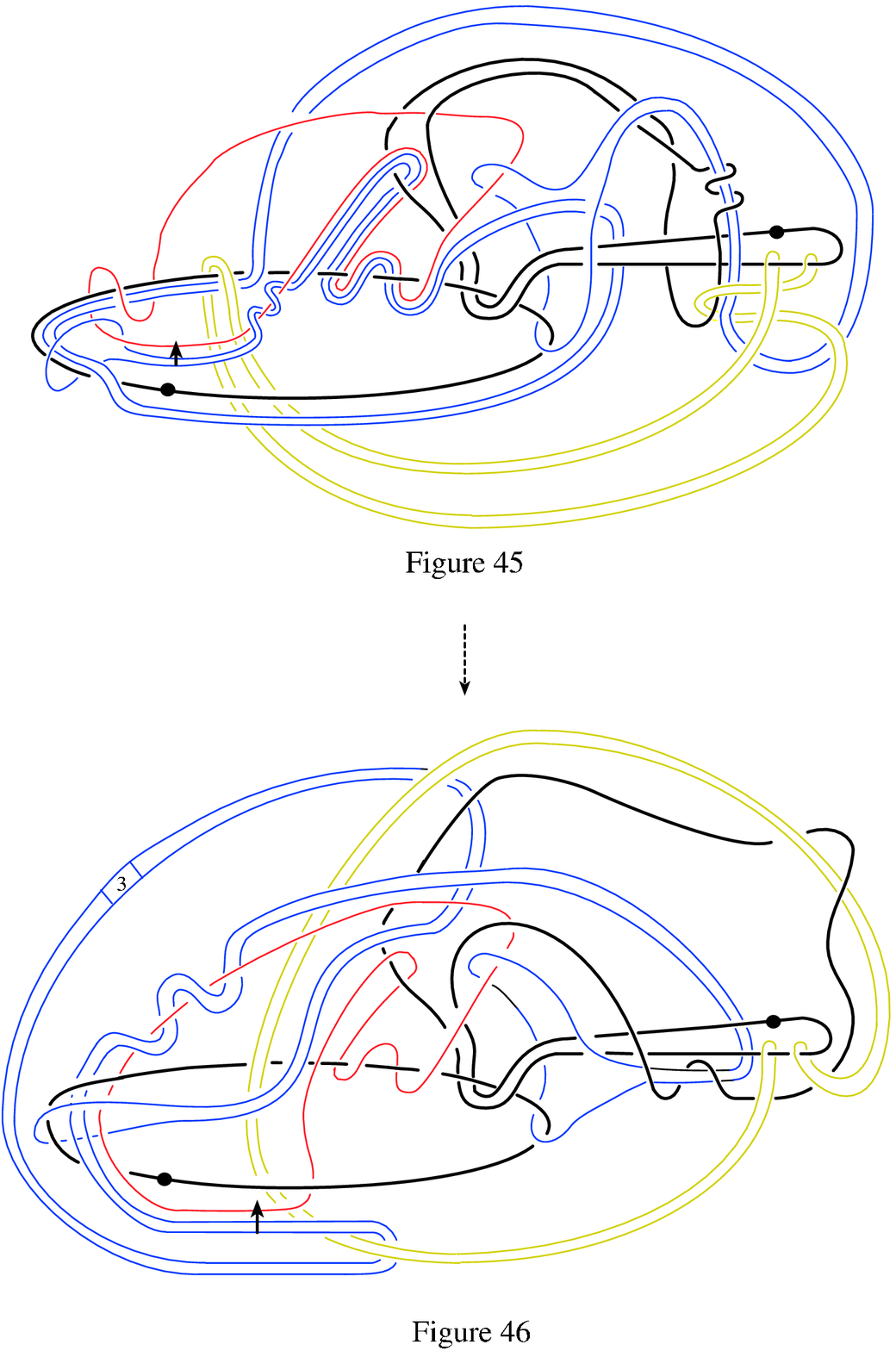}}

\newpage
\hbox{}\par\includegraphics[width=\hsize]{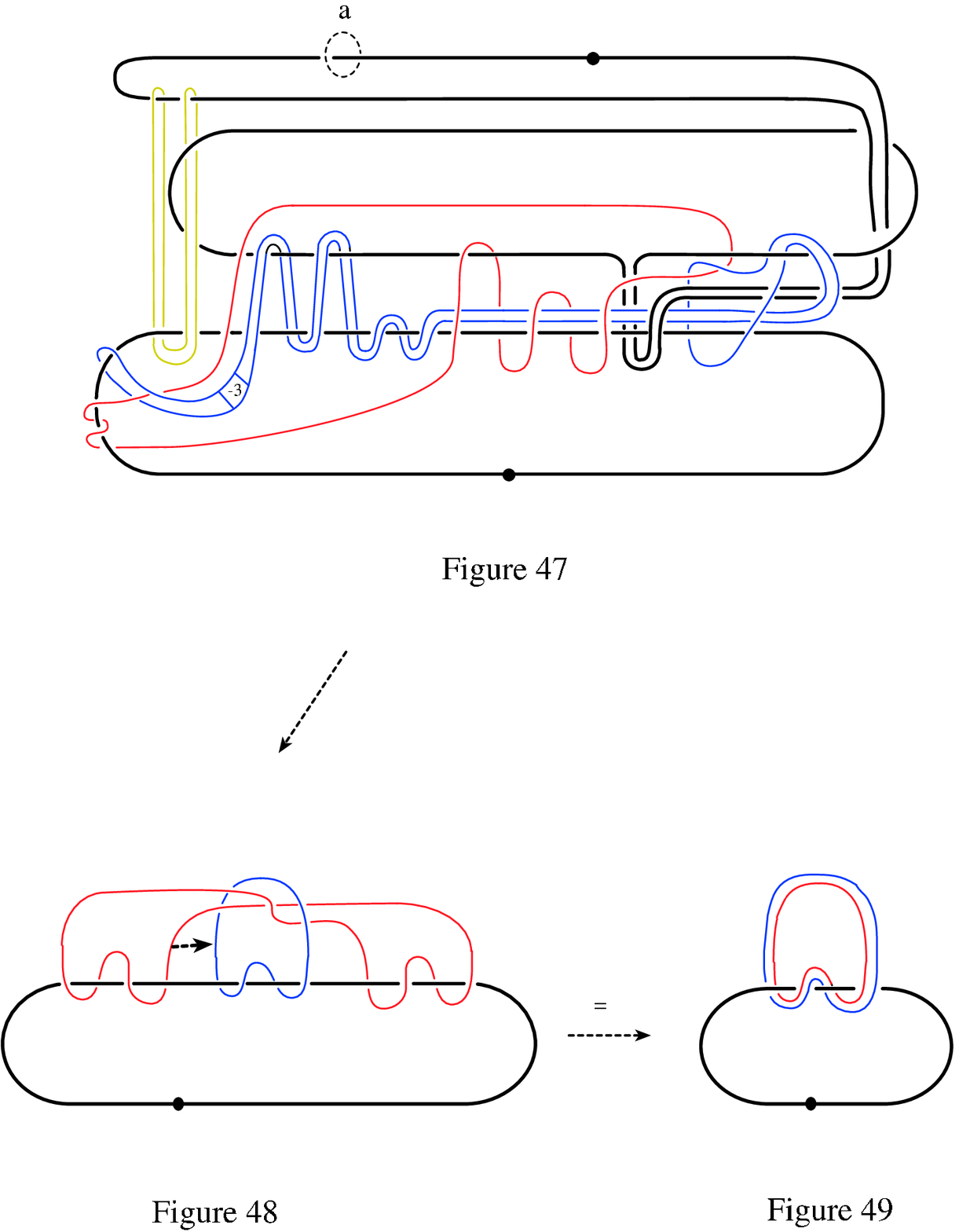}

\newpage

\section{Capping the boundaries of $H$ with $W$}

$H$ has two boundary components homeomorphic
to $Q$, $\partial H = Q_{-}\cup Q_{+}$.
Recall that by capping the either ends of $H$ with $W$ we obtained
$W_{\pm} = H\smile_{Q_{\pm}}W $, and the handlebodies of $W_{-}$, 
$W_{+}$ are obtained by
attaching 
$2$--handles to $H$ along the loops $a$, $b$  in $\partial H =\partial
(Q\times I)$, respectively (Figure 4).

{\Prop  $W_{-} =W$}

\begin{proof} The diffeomorphism $ \partial (Q\times I) \approx \partial H $, takes the loop
in $a$ of Figure 4 to the loop $a$ of Figure 47. By attaching a $2$--handle to
Figure 47, along the $0$ framed loop $a$ (and cancelling the resulting unknotted $0$--framed circles by
$3$--handles) we get Figure 48. By the further indicated handle slide we obtain 
Figure 49. One of the
$2$--handles of Figure 49 slides-off over the other and becomes free, and hence gets
cancelled by a
$3$--handle. So we end up with $W$.
\end{proof}

The story with $W_{+}$ evolves differently, in coming sections we
will see that  $W_{+}$ has a more amusing nontrivial  structure. In the next
section we will use $W_{+}$ to examine $H$ more closely. 

\section{Checking that the boundary of $H$ is correct}

A skeptical reader might wonder how she can verify that the boundary 
of Figure 47 is
the same as the boundary of $Q\times I$ ?  We will check directly
from Figure 47 that it has the same boundary as Figure 4.  This will also be useful
for locating the position of the loops $a$, $b$ in Figure 47. By turning the
$1$--handles to 
$2$--handles (ie, by replacing dotted-circles by 
$0$--framed circles), and by blowing up -- then handle sliding -- then a blowing down operation (done twice)
we obtain Figure 47 $\leadsto$ Figure 50.  By isotopies and 
the indicated handle sliding operations  we
obtain Figure 50 $\leadsto$ Figure 58. By the indicated handle sliding operation, and by 
surgering the $2$--handles of Figure 58 we
obtain Figure 59 which is $ Q\times I$. 

By tracing back
the boundary diffeomorphism Figure 59 $\leadsto $ Figure 47 gives the positions of
the curves $a$ and $b$ on the boundary of Figure 47, which is indicated in Figure 60.
Recall that attaching a $2$--handle to Figure 60 along $b$ (with $0$--framing) gives
$W_{+}$. Attaching a $2$--handle to $b$ and handle slidings and cancelling a $1$-- and $2$-- handle pair gives
Figure 61, a picture of $W_{+}$ which is seemingly different than $W$. Figure 61 
should be considered a
``vertical picture" of $H$ built over the handlebody of $W$ (notice $W$ is visible in
this picture). In the next section we will construct a surprisingly simpler
``vertical" handlebody picture of $H$. 

\section{Vertical handlebody of $H$}

Here we will construct a simpler handlebody picture of $H$ 
as a vertical cobordism starting from the boundary of $W$ to $Q$. This is 
done by stating with Figure 6, which
is $ W=W\smile_{\partial} (Q\times I )$, and then by replacing 
in the interior an imbedded copy of $N\subset Q\times I$ by
$M$. This process gives us $W_{+}$, with an imbedding $ W\subset W_{+} $,
such that $H= W_{+} -W$. This way we will not only
simplifying the handlebody of $H$ but also demonstrate the crucial difference
between $W{+}$ and $W_{-}$   

 We proceed as in Figure 31, except that when we turn
$M$  upside down we add pair of
$2$--handles to the boundary along the loops $H_{1}$ and 
$H_{2}$ of Figure 6 (instead of the loop $h$
of Figure 4). This gives Figure 62. We then apply the boundary diffeomorphism 
Figure 21 $\leadsto $ Figure 28, by carrying the
2--handles $H_{1}$ and $H_{2}$ along the way (we are free to slide $H_{1}$ and $H_{2}$
over the other handles). For example, Figure 68 corresponds to Figure 36. 

By performing the indicated handle slides (indicated by the short arrows) we obtain Figure
68 $\leadsto$ Figure 69, which corresponds to Figure 37. Then by  performing
the indicated handle slides  we obtain Figure 69 $\leadsto$ Figure 71. We then change
the 1--handle notation from pair of balls to the dotted-circles to
obtain Figure 72; and by the indicted handle slides Figure 72 $\leadsto $ Figure 77 we arrive to Figure 77.
Figure 77 is the desired picture of $W_{+}$.

Next we check that the
boundary of the manifold of Figure 77 is correct. This can easily be done by
turning one of the dotted circles to a $0$--framed circle, and turning a
0--framed circle to a dotted circle as in Figure 78 and then by cancelling the
dotted circle with the $-1$ framed circle which links it geometrically
once  (ie, we cancel a 1-- and 2--handle pair). It easily checked that
this operation results $W$ and a disjoint $0$--framed unknotted circle, which is then
cancelled by a $3$--handle. So we end up with $W$, hence
$W_{+}$ has the correct boundary.

\eject

\cl{\includegraphics[width=4in]{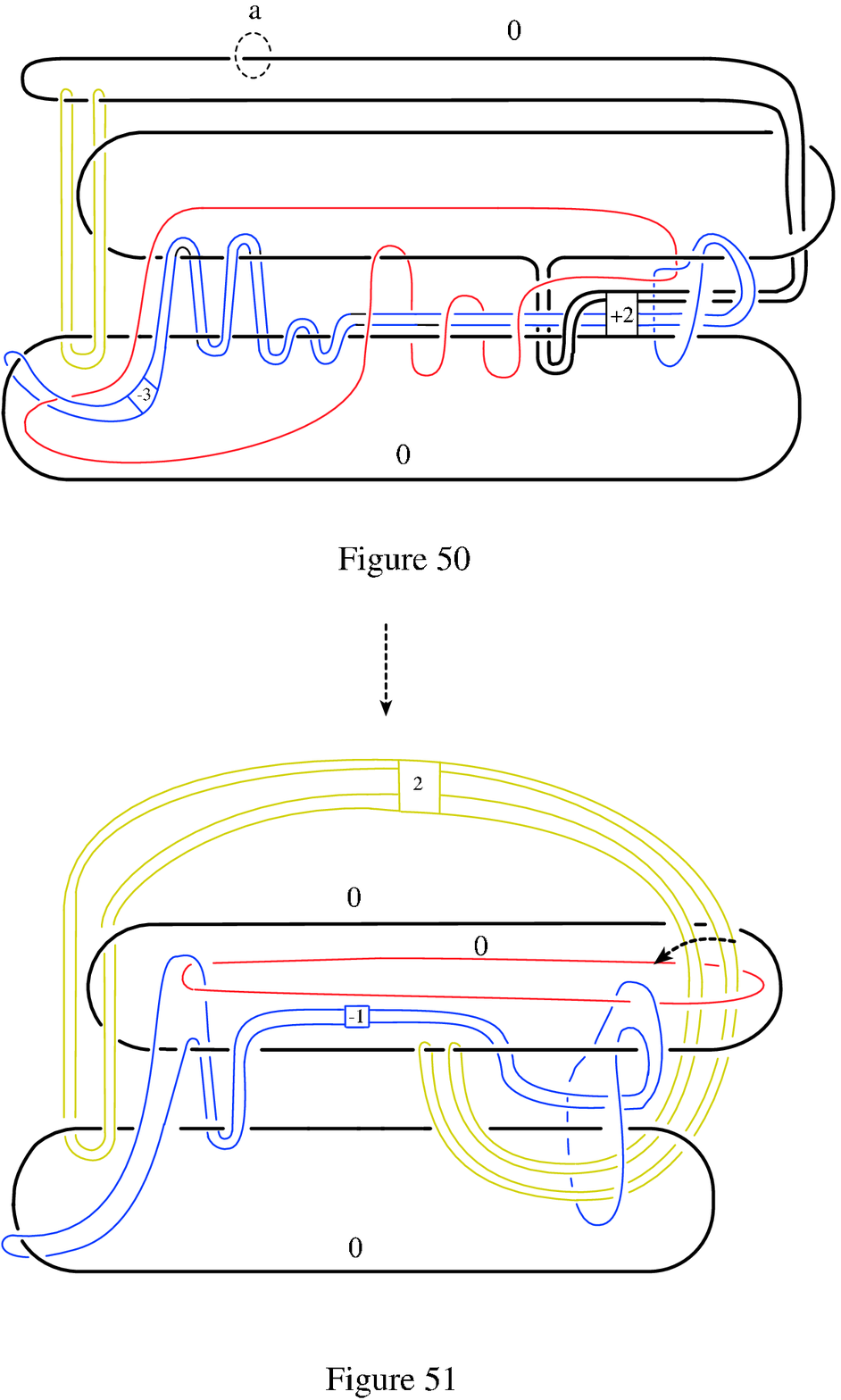}}

\newpage
\cl{\includegraphics[width=4in]{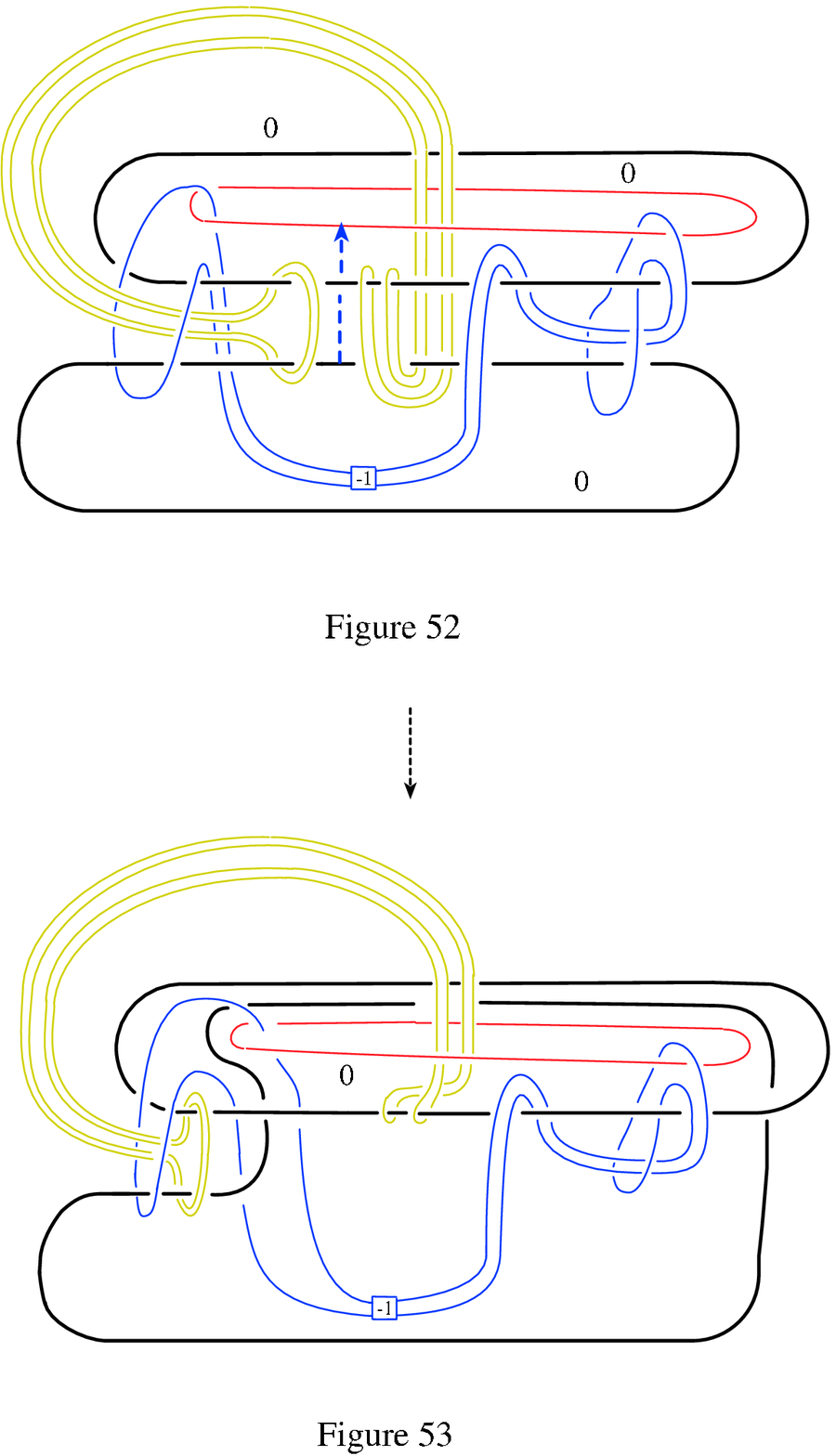}}

\newpage
\cl{\includegraphics[width=4in]{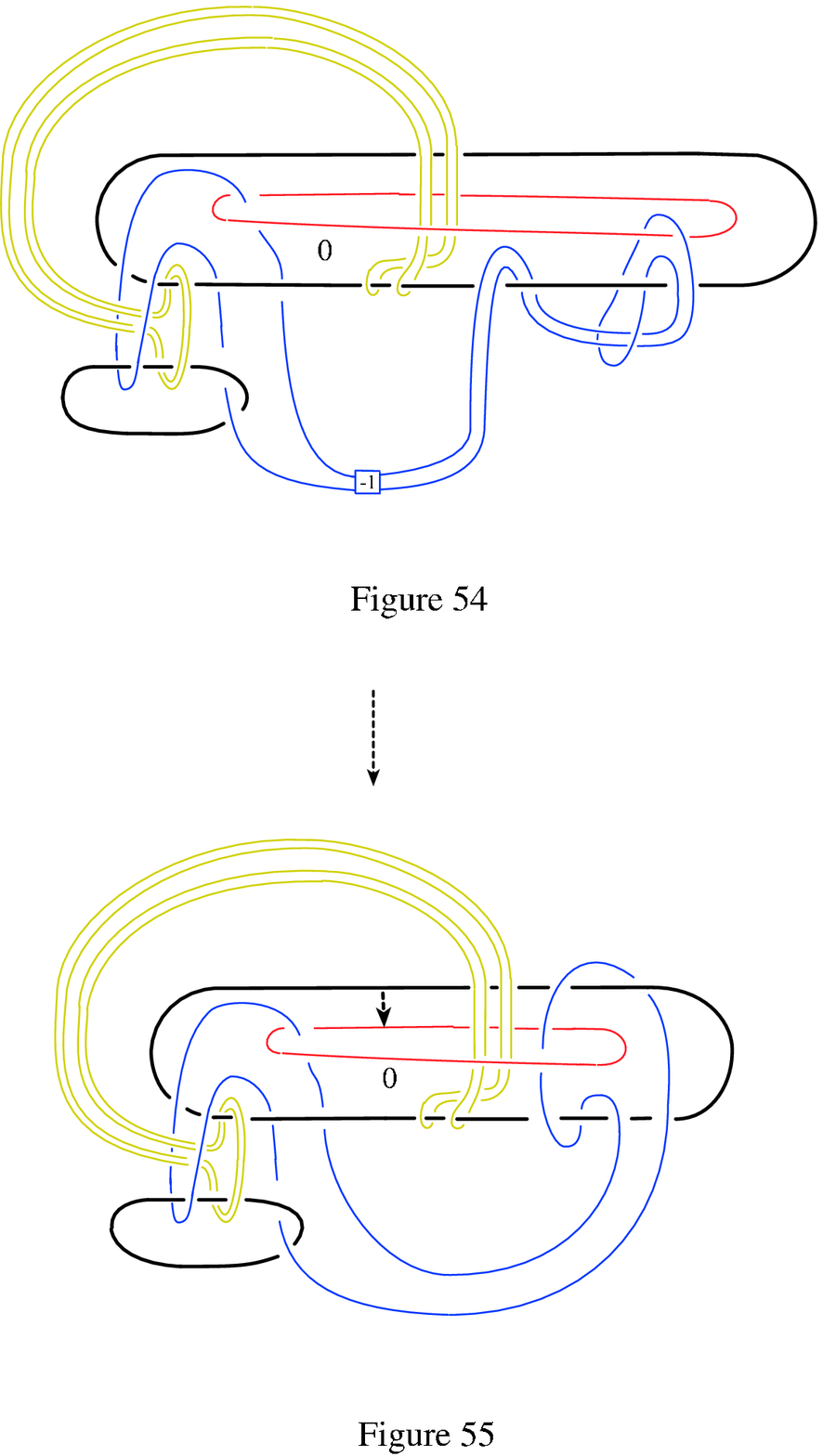}}

\newpage
\cl{\includegraphics[width=3.3in]{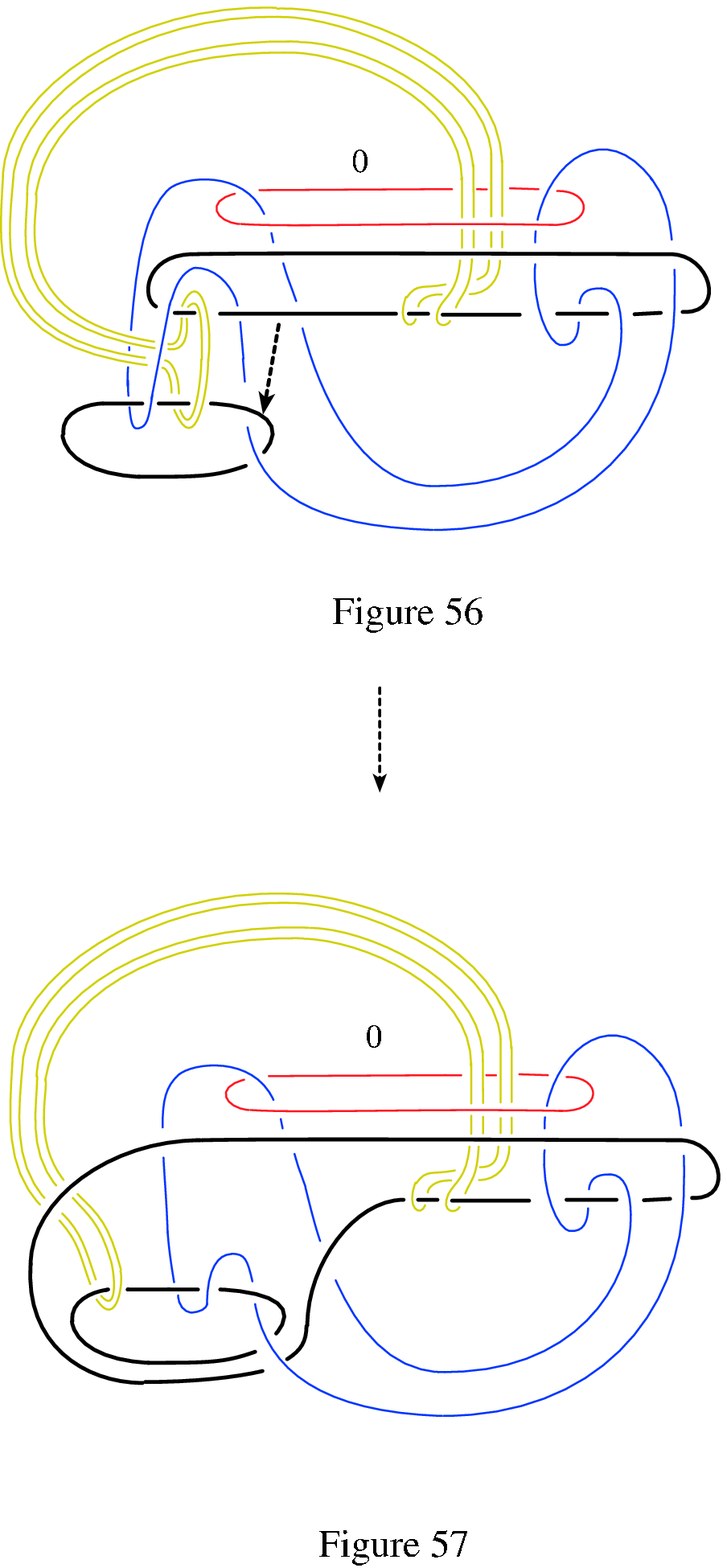}}

\newpage
\cl{\includegraphics[width=3.5in]{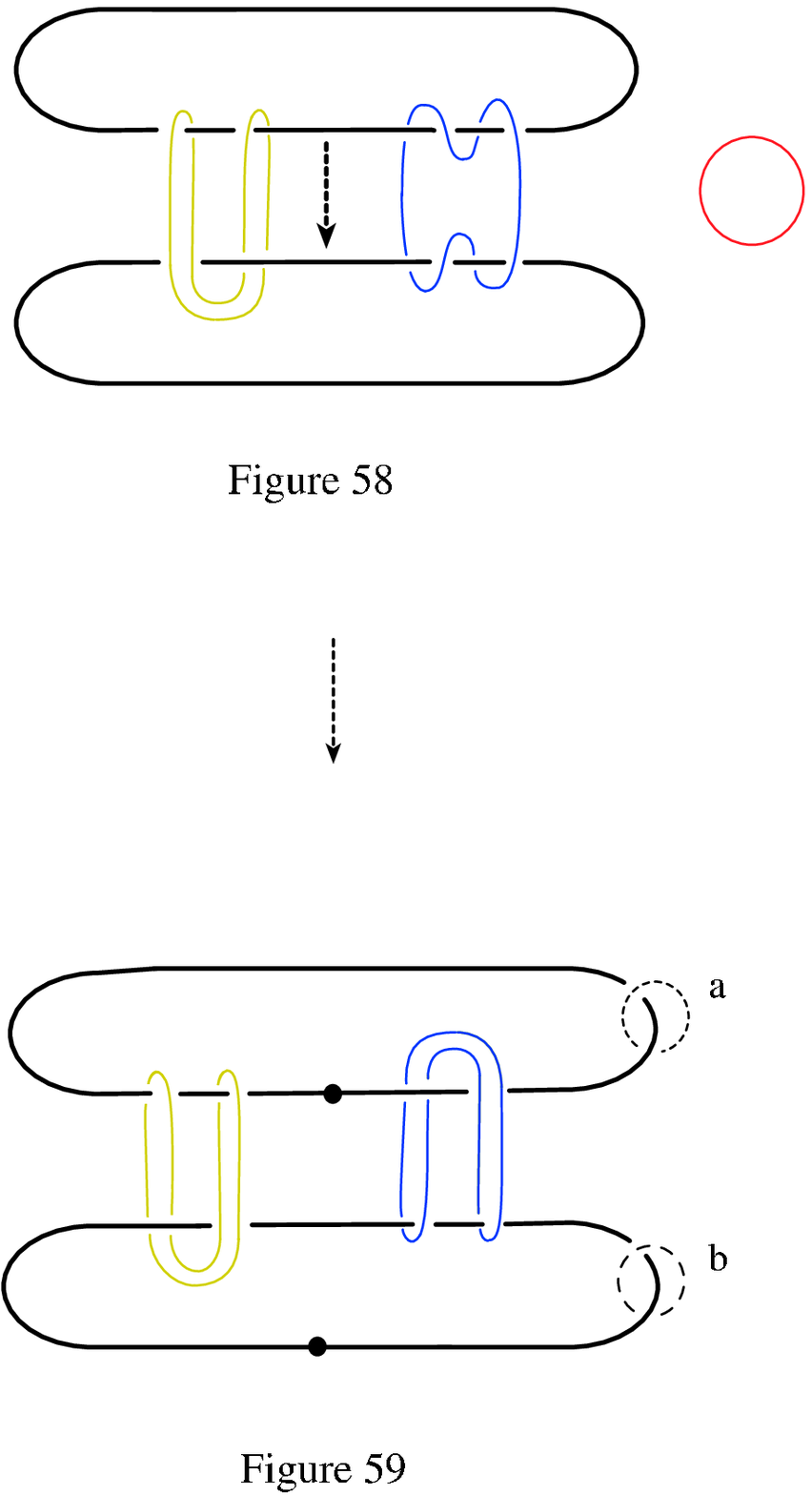}}

\newpage
\cl{\includegraphics[width=4in]{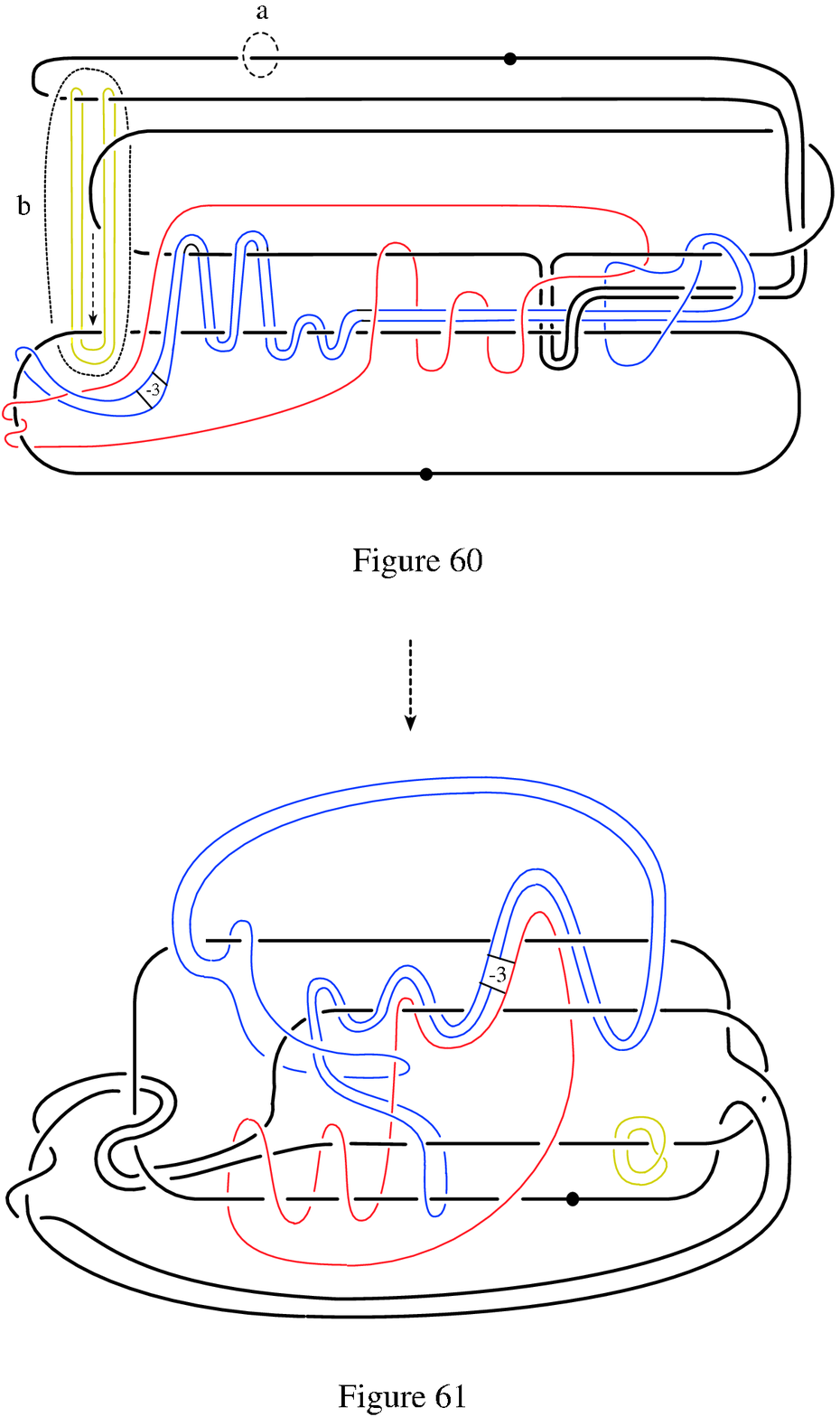}}

\newpage
\hbox{}\par\includegraphics[width=\hsize]{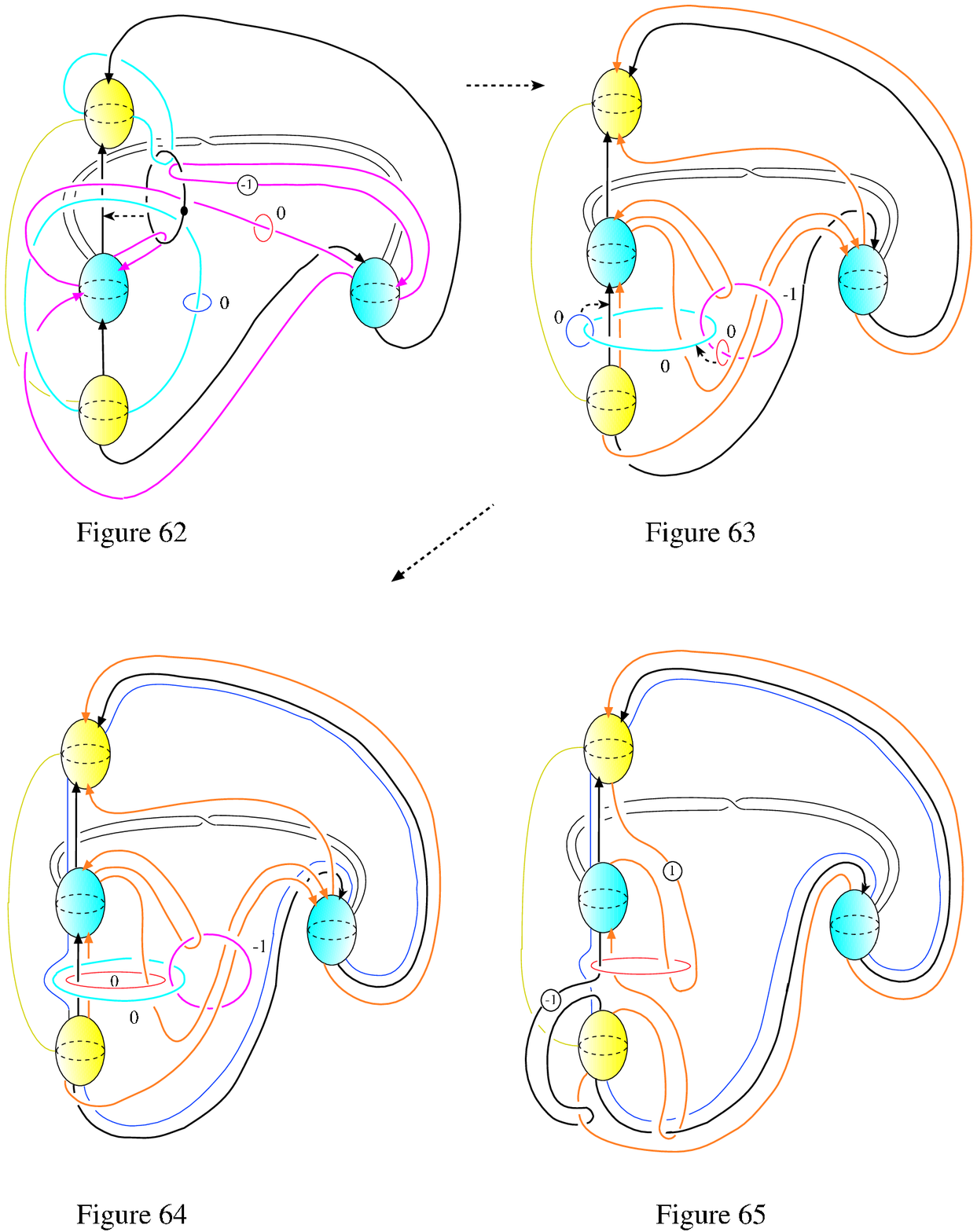}

\newpage
\hbox{}\par\includegraphics[width=\hsize]{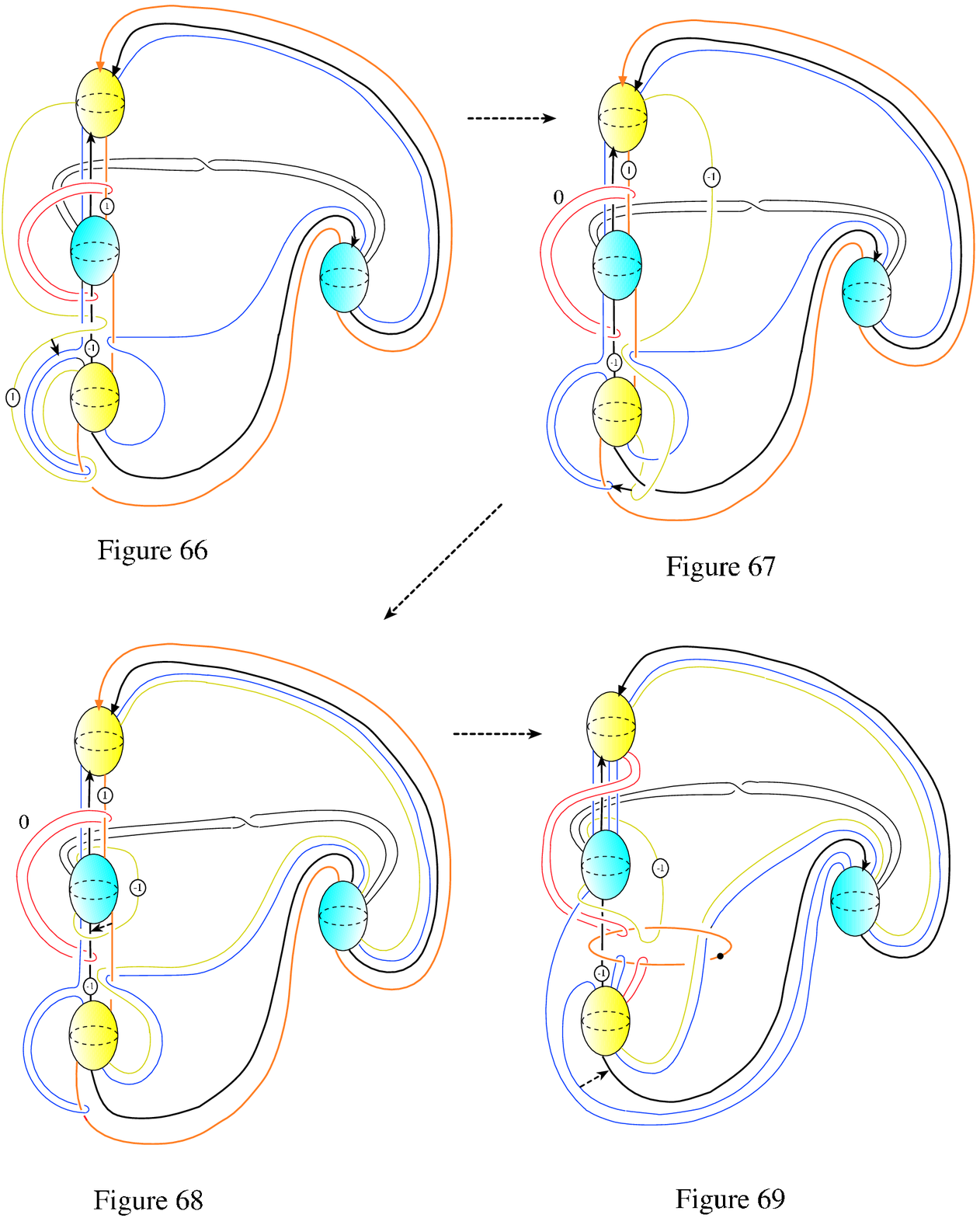}

\newpage
\hbox{}\par\includegraphics[width=\hsize]{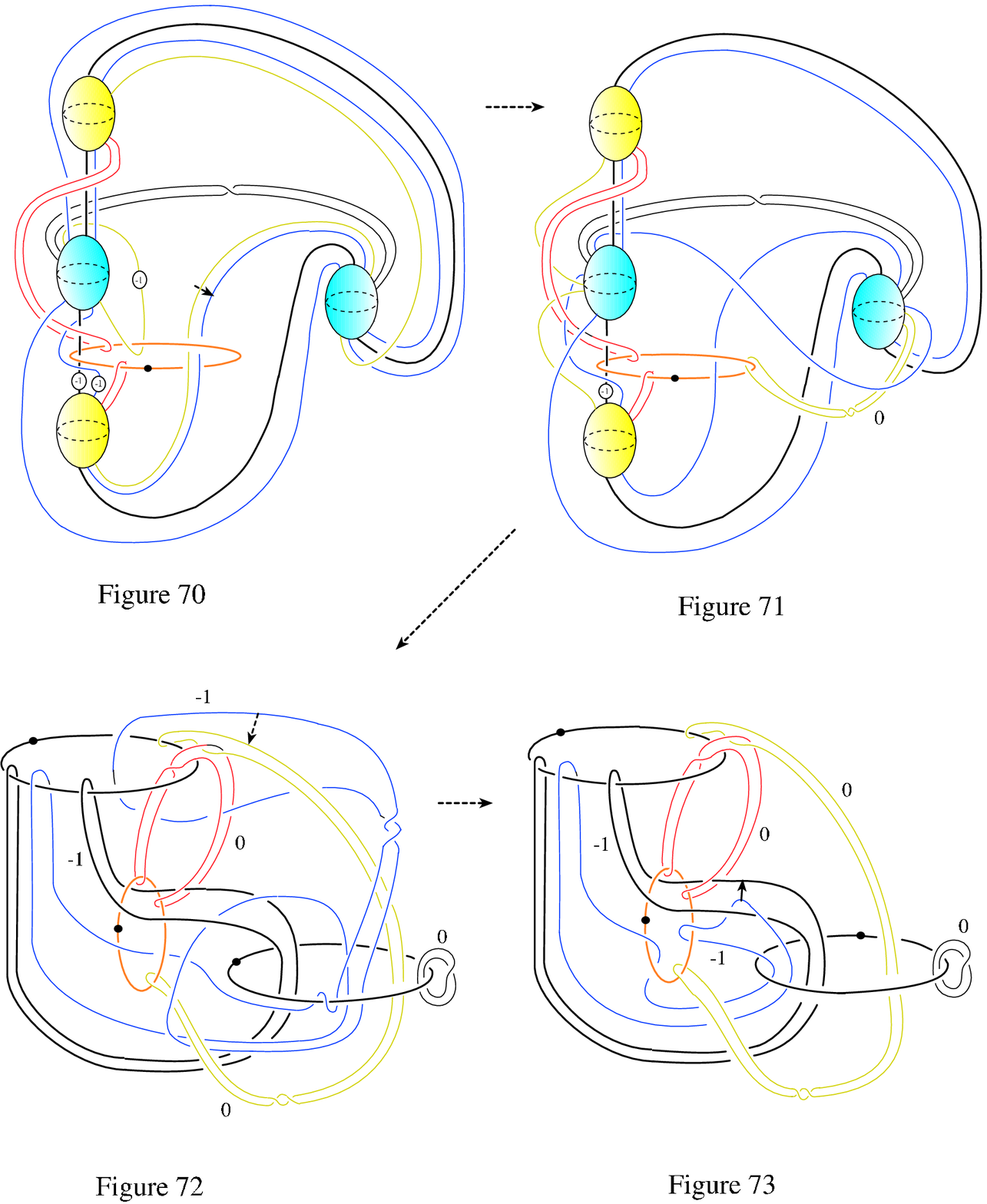}

\newpage

\hbox{}
\vglue 0.5in
\includegraphics[width=\hsize]{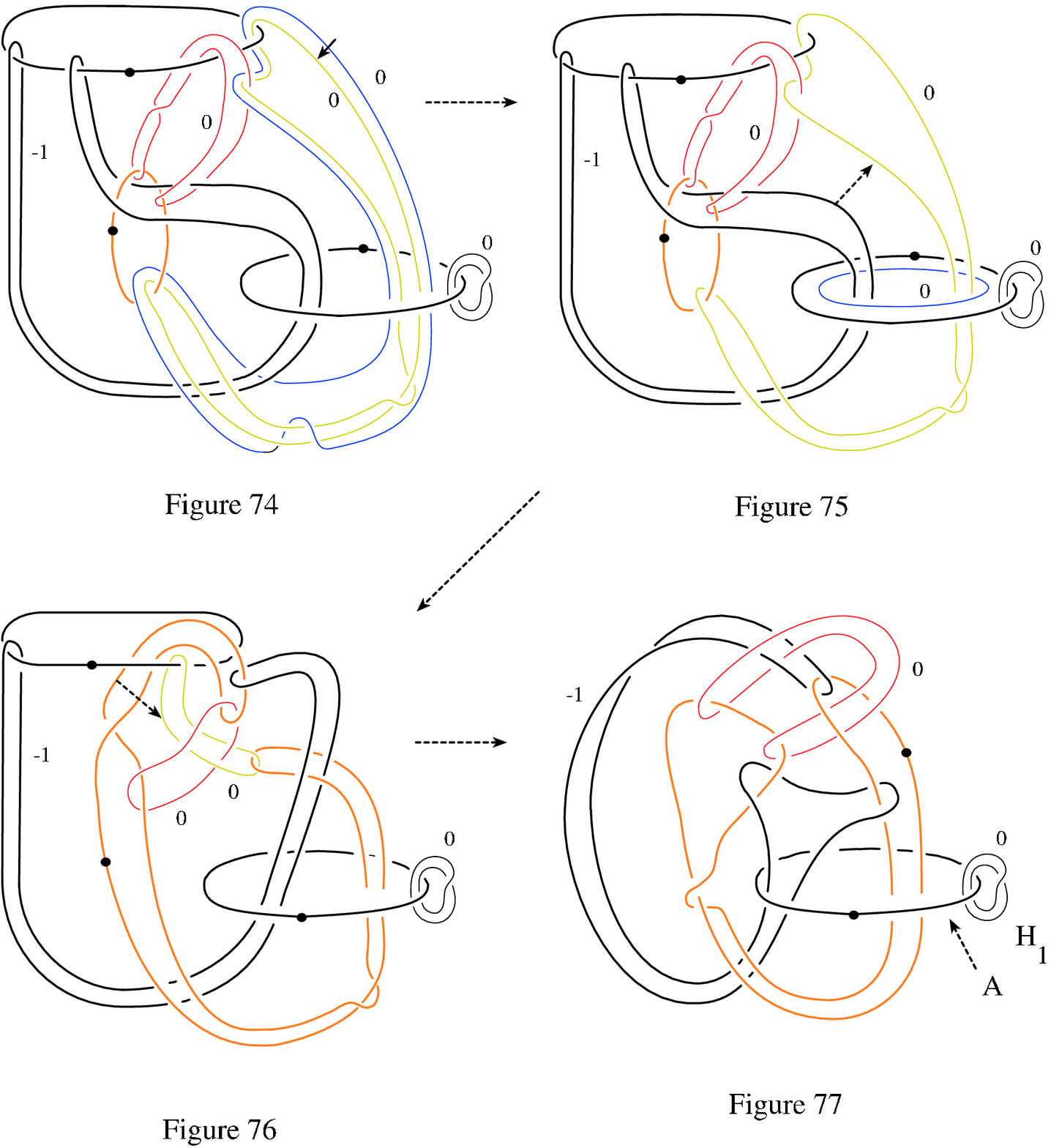}

\newpage

\section{A knot is born}

By twisting the strands going through the middle $1$--handle, and then
by sliding $H_{1}$ over this $1$--handle we see that Figure 77 is
diffeomorphic to the Figure 79 (similar to the move at the bottom of
the page 504 of \cite{a2}).  Figure 79 demonstrates the complement of
the imbedding $f\co B^{2} \hookrightarrow B^{4}$ with the standard
boundary, given by the ``dotted circle" A (the $1$--handle). This
follows from the discussion on the last paragraph of the last
section. Because there by changing the interior of Figure 77 we
checked that as a loop A is the unknot on the boundary of the
handlebody $X$ consisting of all the handles of Figure 77 except the
$1$--handle corresponding to $A$ and the last $2$--handle $H_{1}$.
The same argument works for Figure 79. In addition in this case, the
handlebody consisting of all the handles of the Figure 79 except the
1-handle corresponding to A and the the 2-handle $H_{1}$, is $B^4$.  So
$W_{+}$ is obtained from $B^4$ by carving out of the imbedded disk bounded
by the unknot A (i.e creating a 1--handle A) and then by attaching the
$2$--handle $H_{1}$. Hence by capping with a standard pair $(B^4,
B^2)$ we can think of $f$ as a part of an imbedding $ F\co S^{2}
\hookrightarrow S^4 $. Let us call $A=F(S^2)$.

We can draw a more concrete picture of the knot $A$ as follows: During the next few
steps, in order not to clog up the picture, we won't draw the last $2$--handle $H_{1}$.
By an isotopy and cancelling $1$-- and $2$--handle pair we get a diffeomorphism from
Figure 79, to Figures 80, 81 and finally to Figure 82. In Figure
82 the `dotted' ribbon knot is really the unknot in the presence of a cancelling $2$
and $3$--handle pair (ie, the unknotted circle with $0$--framing plus the $3$--handle
which is not seen in the figure). So this ribbon disk with the boundary the unknot in
$S^{3}$, demonstrates a good visual picture of the imbedding 
$f\co  B^{2} \hookrightarrow B^{4}$.

\subsection{A useful fundamental group calculation}

We will compute the fundamental group of the $2$--knot complement $S^4-A$, ie, we will
compute the group $G:=\pi_{1} (Y)$, where $Y$ is the handlebody consisting of all the
handles of Figure 79 except the
$2$--handle
$H_{1}$. Though this calculation is not necessary for the rest of the paper, it is
useful to demonstrate why $W_{+} -W$ gives an $s$--cobordism. By using the generators
drawn in Figure 83 we get the following relations for $G$:
$$\begin{array}{lcll}
(1)&x^{-1}yt^{-1}x^{-1}t&=&1 \\
(2)&x^{-1}yxy&=&1\\
(3)&txt^{-1}y^{-1}&=&1
\end{array}$$
From $(1)$ and $(2) $ we get $ t^{-1}x^{-1}t=xy $ , then by using $(3)$
$\Longrightarrow t^{3}=(tx)^{3}
$. Call
$ a=tx $ so
$t^{3}=a^{3} $. By solving $y$ in $(3)$ and plugging into $(2)$ and substituting
$x=t^{-1}a$ we get $ata=tat^{-2}a^{3}=tat^{-2}t^{3}=tat$. Hence we get the presentation: 
$$ G=\< t,a | t^3=a^3, ata=tat\>$$ 
Notice that attaching the $2$--handle $H_{1}$ to $Y$ (ie, forming $W_{+}$) 
introduces the extra relation
$t^{4}=1$ to 
$G$, which then implies $t=x$, demonstrating an $s$--cobordism from the boundary of W to the boundary to
$W_{+}$.  The following important observation of Meierfrankenfeld has motivated 
us to prove the crucial fibration structure for $A$ in the next section.

{\Lem {\rm\cite{m}}\qua $G$ contains normal subgroups ${\bf Q_{8}}$ and ${\bf Z}$ 
giving the exact sequences}:
$$  1\to {\bf Q_{8}}\to G\to {\bf Z}\to 1$$
$$ 1\to {\bf Z}\to G\to SL(2,{\bf Z}_{3})\to 1 $$
\begin{proof} Call $u:=ta^{-1}=y^{-1}$ and $v:=a^{-1}t = x^{-1}$. First notice that the group
$\<u,v\>$ generated by $u$ and $v$ is a normal subgroup. For example, since
$u=ta^{-1}=a^{-1}t^{-1}at=a^{-1}v^{-1}t=a^{-1}v^{-1}av $ $\Longrightarrow$ 
 $a^{-1}va=uv^{-1}\in \; \<u,v\>$.
Also since $v=a^{-1}t=a^{-1}ua$ $\Longrightarrow$  $a^{-1}ua=v\in \;\<u,v\>$. 
Now we claim that $ \<u,v\>={\bf Q_{8}} $. This follows from $a^{3}=t^{3}\in \mbox{Center}(G)$
$\Longrightarrow u=a^{-3}ua^{3}= a^{-1} (vu^{-1})a=vu^{-1}v^{-1}$ $\Longrightarrow uvu=v$. So
$vu^{-1}=a^{-1}va= a^{-1}(uvu)a= (a^{-1}ua)(a^{-1}va)(a^{-1}ua)= v(vu^{-1})v$, implying $vuv=u$. So
$\<u,v\>=\<u,v \;| \; uvu=v, vuv=u\>$, which is a presentation of ${\bf Q_{8}}$.

For the second exact sequence take ${\bf Z}=\<t^{3}\>$ and then observe that $ G/\<t^3\>=SL(2,{\bf Z}_{3}) $ 
(for example, by using the symbolic manipulation program GAP one can check that $G/\<t^{3}\>$ has order
$24$, then use the group theory fact that $SL(2,{\bf Z}_{3})$ is the only group of order 24 generated by
two elements of order 3)
\end{proof}

\subsection{Fiber structure of the knot A}

Consider the order three self diffeomorphism $ \phi \co Q\to Q $ of Figure 84.  As
described  by the the pictures of Figure 85, this diffeomorphism is obtained by the
compositions of blowing up, a handle slide, blowing down, and
another handle slide operations. $\phi$
permutes the circles $ P,Q,R $ as indicated in Figure 84, while twisting the tubular
neighborhood of $R$ by $-1$ times. Note that $Q$ can be obtained by doing $-1$ surgeries to
three right-handed Hopf circles, then $\phi$ is the map induced from the map
$S^{3}\to S^{3}$ which permutes the three Hopf circles.
Let $Q_{0}$ denote the punctured
$Q$, then:

{\Prop The knot $A\subset S^{4}$ is a fibered knot with fiber $Q_{0}$ and monodromy $\phi$.}

\begin{proof} We start with Figure 86 which is the knot complement $Y$. By introducing a zero framed
unknotted circle (ie, by introducing a cancelling pair of $2$--and $3$--handles) we arrive to the
Figure 87. Now something amazing happens!: This new zero framed unknotted circle isotopes to the
complicated looking circle of  Figure 88, as indicated in the figure. The curious reader can check
this by applying the boundary diffeomorphism Figure 77 $\Rightarrow$ Figure 78 to Figure 88
(replacing dotted circle with a zero framed circle) and tracing this new loop along the way back to
the trivial loop!. By isotopies and the indicated handle slides, from Figure 88 we arrive to Figure 92. 

Now in Figure
92 we can clearly see an imbedded copy of $ Q\times [0,1]$ (recall Figures 3 and 4). We claim
that, in fact the other handles of this figure has the role of identifying the two ends of 
$ Q\times [0,1]$ by the monodromy $\phi$. To see this,  recall from \cite{ak} how to draw the
handlebody of picture of :
$$Q_{0}\times [0,1]/(x,0)\sim (\phi(x),1)$$
For this, we attach a $1$--handle to $Q_{0}\times [0,1]$ connecting the top to the bottom, and
attach $2$--handles along the loops $\gamma \;\# \;\phi(\gamma)$ where
$\gamma$ are the core circles of the $1$--handles of $Q\times{0}$ and $\phi(\gamma)$ are their 
images in $Q\times{1}$ under the map
$\phi$ (the connected sum is taken along the $1$--handle).  By
inspecting where the $2$--handles are attached on the boundary of $Q_{0}\times [0,1]$ (Figure 93), we
see that in fact the two ends are identified exactly by the diffeomorphism $\phi $. Note
that, by changing the monodromy of Figure 94 by $\phi^{-1}$ we obtain Figure 95, which
is the identity monodromy identification $Q\times S^{1}$. 
\end {proof}

\section{The Gluck Construction}

Recall that performing the Gluck construction to $S^{4}$ along an imbedded $2$--sphere
$S^{2} \hookrightarrow S^{4}$ means that we first thicken the imbedding 
$S^{2}\times B^{2} \hookrightarrow S^{4}$ and then form: 
$$\Sigma=(S^{4} -  S^{2}\times B^{2} )\smile_{\psi} S^{2}\times B^{2}$$
where $\psi\co  S^{2}\times S^{1} \to S^{2}\times S^{1} $ is the diffeomorphism given by
$\psi(x,y)=(\rho(y)x,y) $, and
$\rho\co  S^{1}\to SO(3)$ is the generator of $\pi_{1}(SO(3))={\bf Z}_{2}$. 

\eject

\hbox{}\par\includegraphics[width=\hsize]{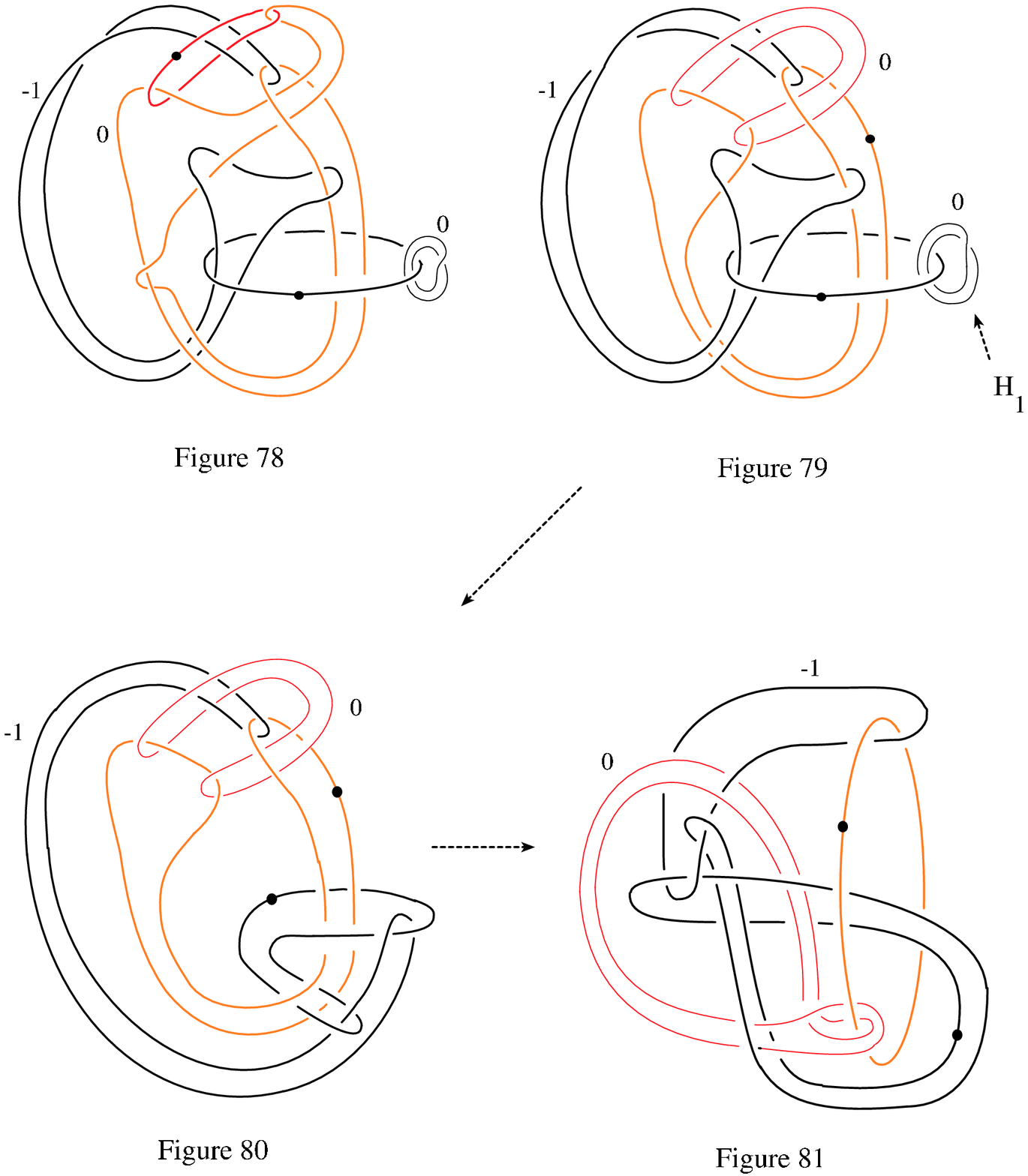}

\newpage
\cl{\includegraphics[width=3.5in]{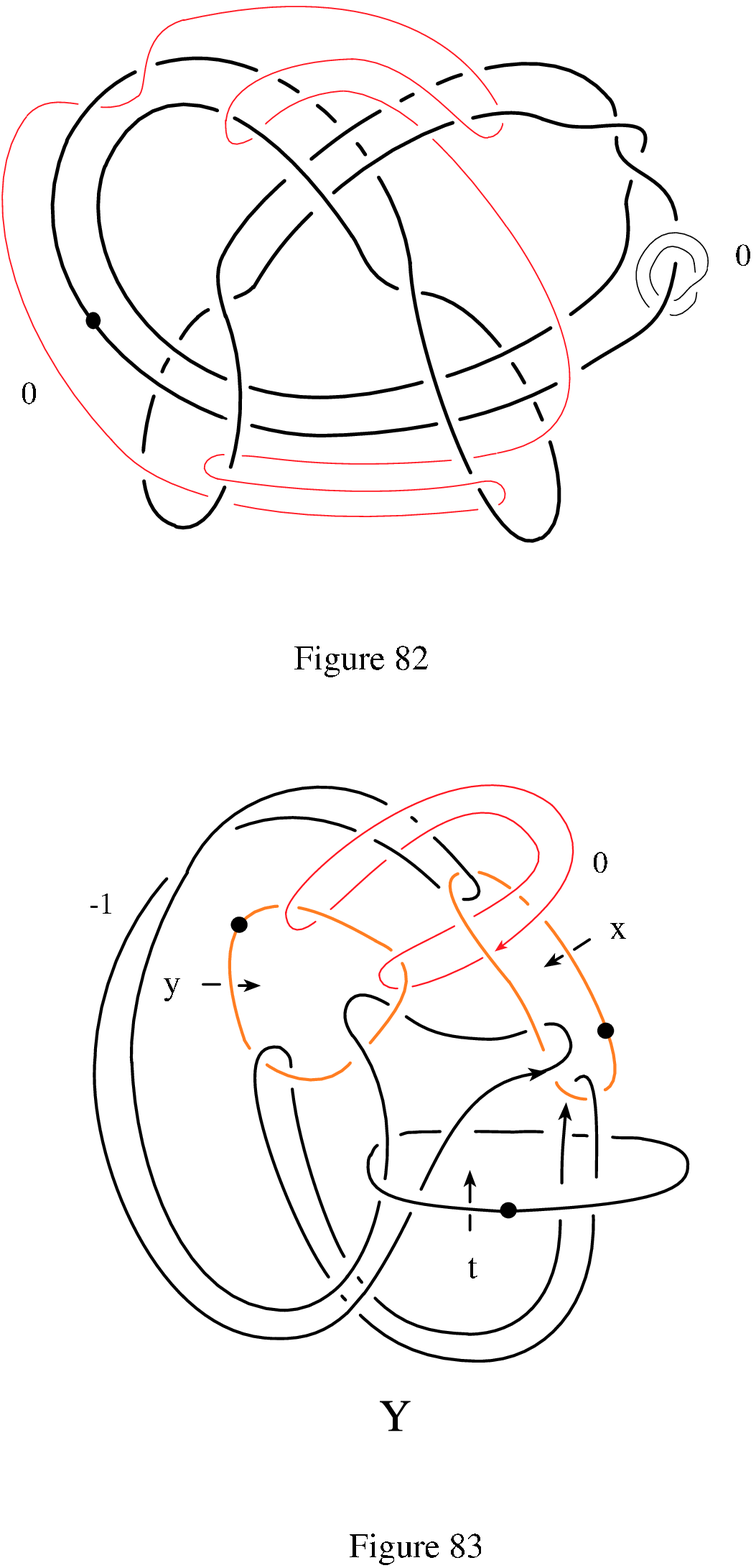}}

\newpage
\includegraphics[width=\hsize]{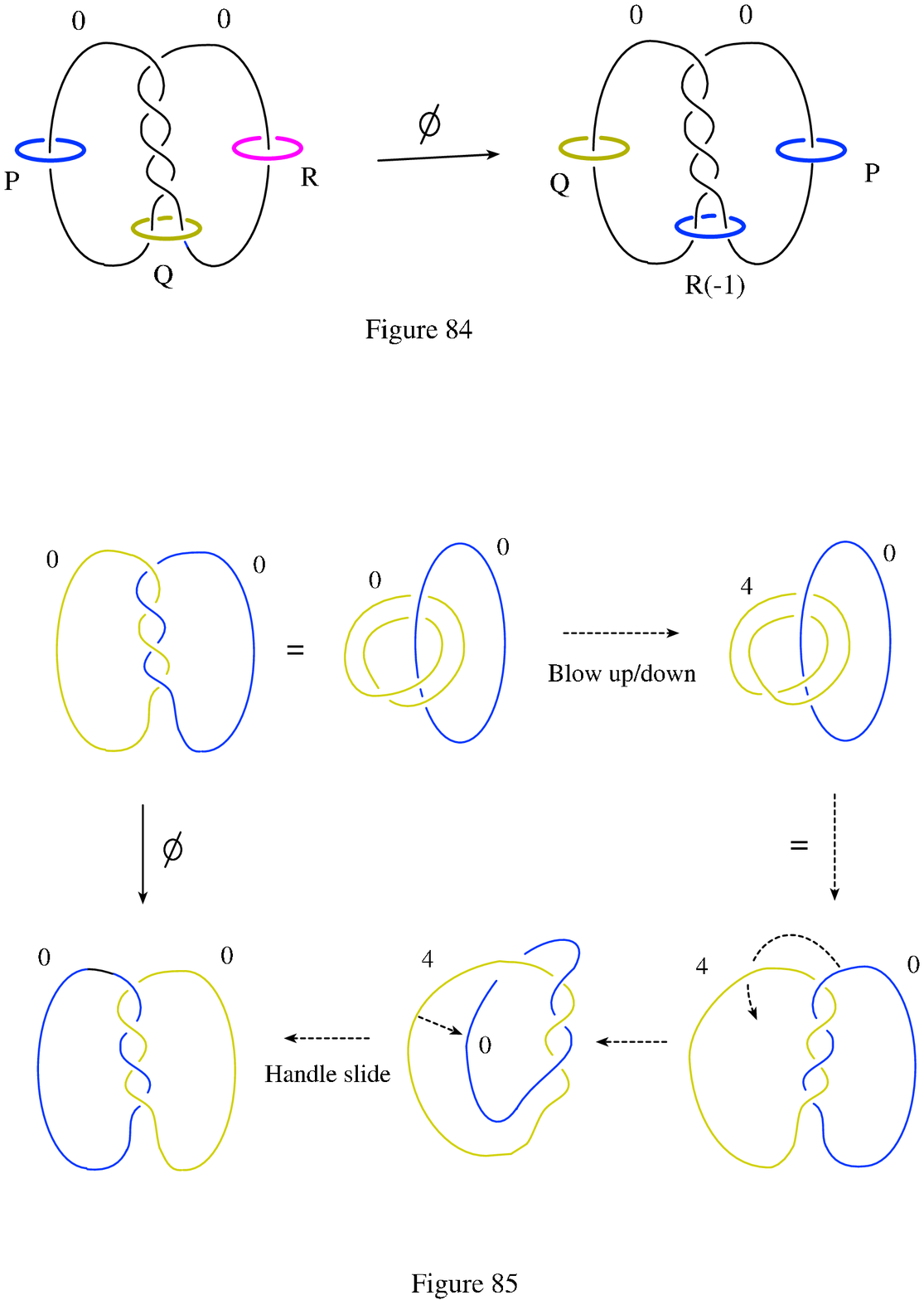}

\newpage
\hbox{}\par\includegraphics[width=\hsize]{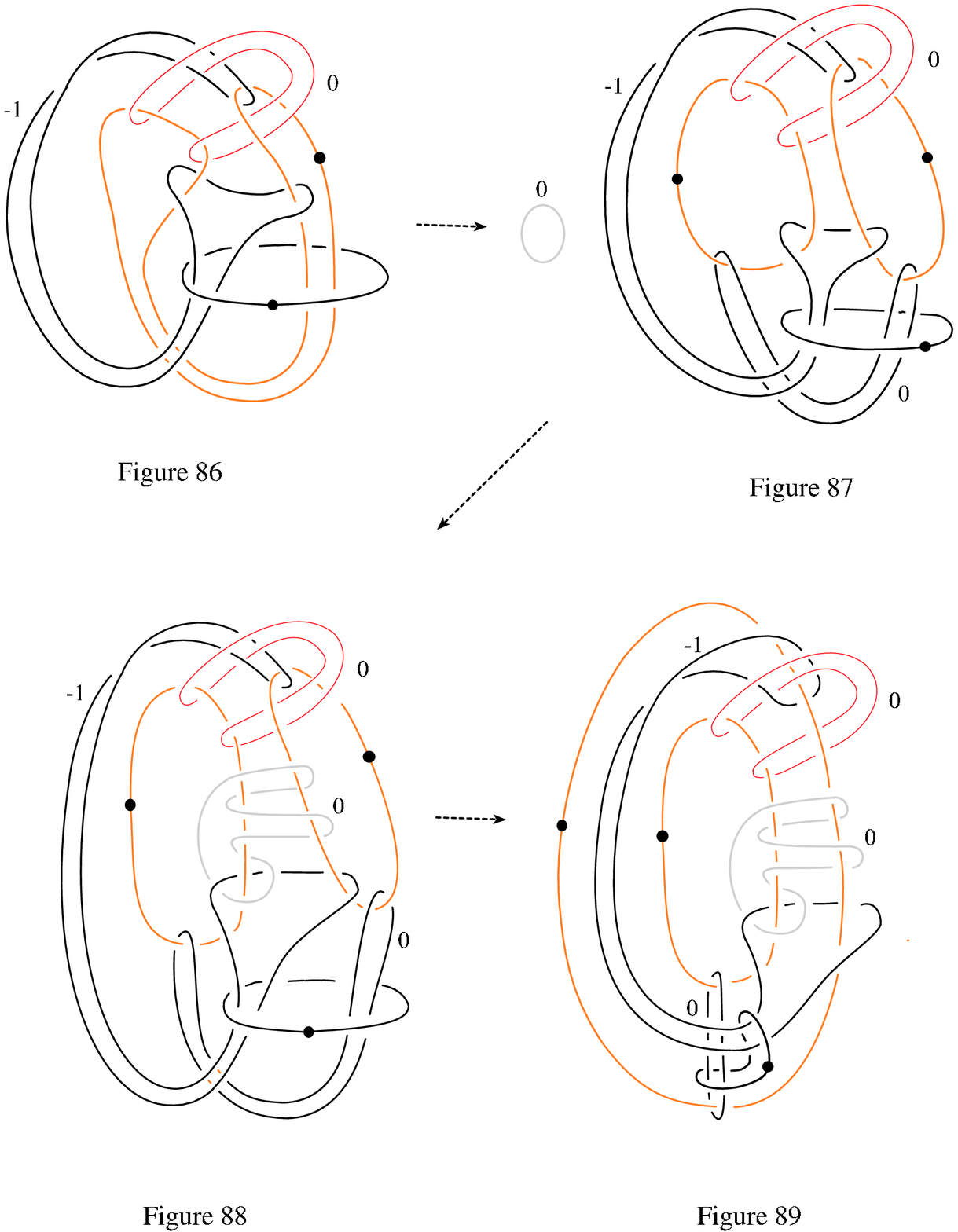}

\newpage
\hbox{}\par\includegraphics[width=\hsize]{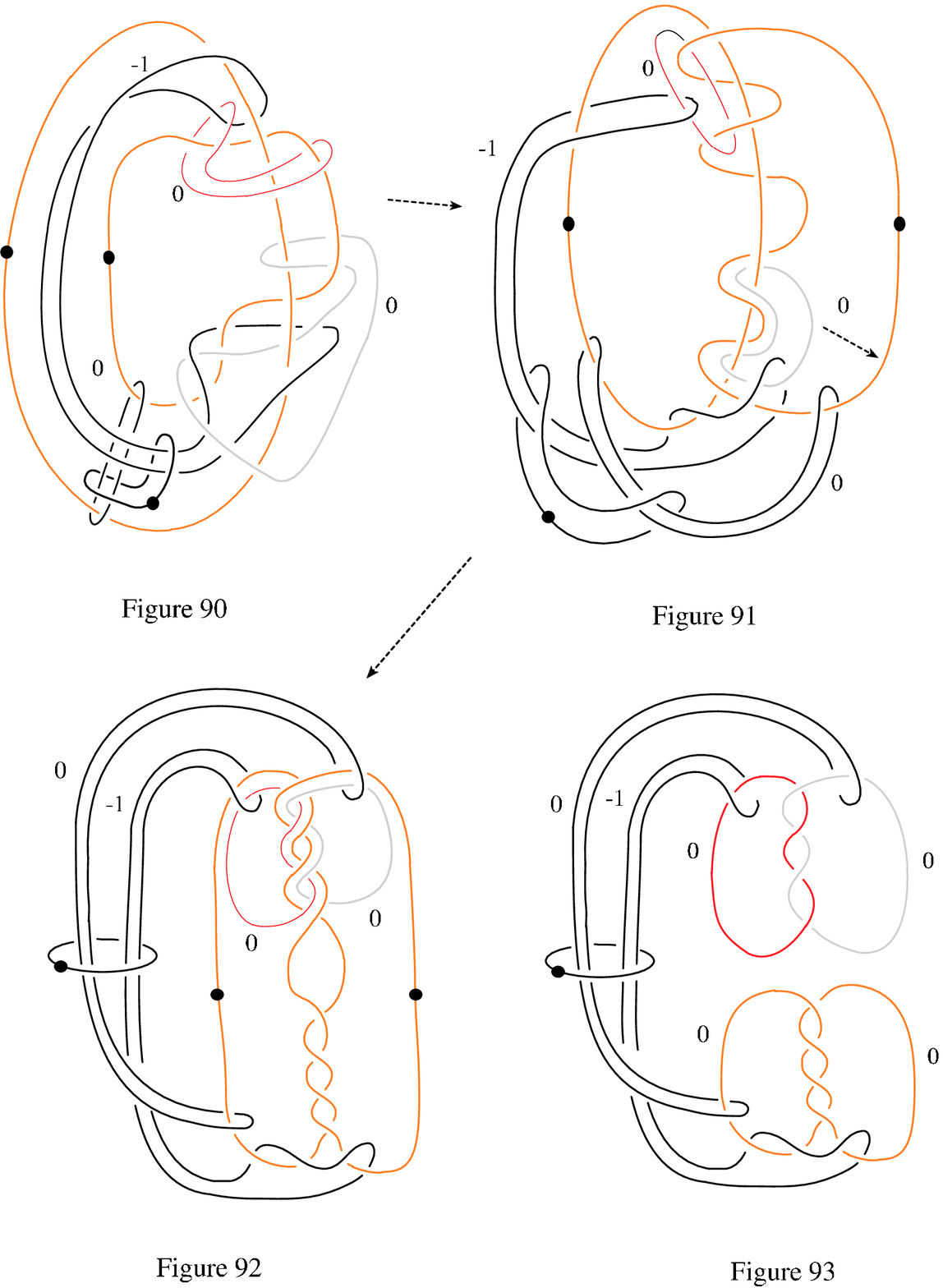}

\newpage
\hbox{}\par\includegraphics[width=\hsize]{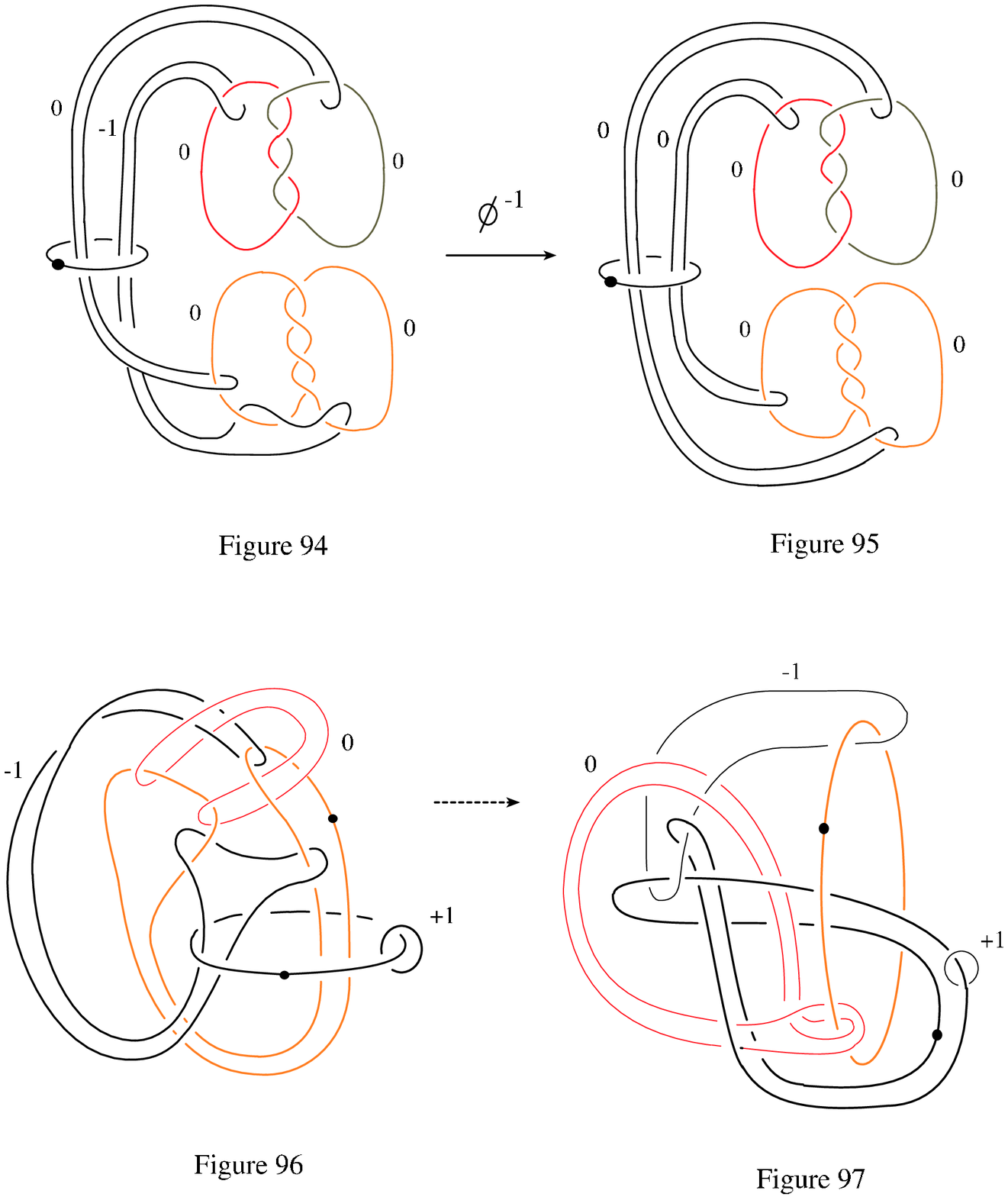}

\newpage

{\Prop The Gluck construction to $S^4$ along the $2$--knot $A$ gives
back $S^4$.}

\begin{proof} The handlebody of Figure 96 describes the Gluck
construction  of $S^4$ along $A$, which can be equivalently described 
by Figure 97 (recall  Figure 79 $\leadsto $ Figure 82 identification). 
By introducing a cancelling  $1$ and $2$--handle pair we get Figure 98, 
and by cancelling $1$ and $2$--handle pair from Figure 98 we obtain Figure 99. 
 
Now in Figure 100 we introduce a cancelling pair of $2$ and $3$--handles. In Figure 100
only the new $2$--handle ``m'' is visible as the zero framed unknotted circle. It is
important to check that this new $2$--handle is attached along the unknot! (this can be
checked by tracing m along the boundary diffeomorphisms Figure 100 $\leadsto $ Figure
101, and  Figure 98 $\leadsto $ Figure 97, and Figure 81 $\leadsto $ Figure 78). 

In Figure 100 by sliding the $+1$ framed $2$--handle over the $0$ framed $2$--handle
m, we obtain Figure 102. Now comes an important point!: Notice that the $1$--handles
of Figure 102 are cancelled by $2$--handles (ie, through each circle-with-dot there is
a framed knot going through geometrically once). So in fact after cancelling $1$ and
$2$--handle pairs, Figure 102 becomes a handlebody consisting of only two $2$--handles and
two $3$--handles. Now, rather than performing these $1$ and $2$--handle cancellations and
drawing the resulting handlebody of $2$ and $3$--handles, we will turn the handlebody of
Figure 102 upside down. This process is performed  by taking the dual loops of the
$2$--handles as in Figure 103 (ie, the small $0$--framed circles), and by tracing them
under the boundary  diffeomorphism from the  boundary of the handlebody of Figure 103
to  $\partial (S^{1}\times B^{3}\;\# \;S^{1}\times B^{3})$, and then by attaching
$2$--handles to $ S^{1}\times B^{3}\;\# \;S^{1}\times B^{3} $ 
along the images of these dual loops. It is important to note that along this process we
are allowed to slide the dual $2$--handles over each other and over the other
handles. 

By a blowing up and down operation, and by isotopies and the
indicated handle slides Figure 103 $\leadsto$ Figure 112 we arrive to the handlebody of 
Figure 112, and by sliding one dual $2$--handle over the other one we obtain Figure
113. Now by applying Figure 113 to the boundary diffeomorphism Figure 77 $\leadsto$
Figure 78  we obtain Figure 114 (note that the handlebody  of Figure 114 is
just $ S^{1}\times B^{3}\;\# \;S^{1}\times B^{3} $, it happens to look complicated
because of the presence of the dual $2$--handles).  By
sliding the dual $2$--handles over each other (as indicated in the figures), and by a blowing up and down
operation and isotopies  we arrive to Figure 121 which is $B^4$.
\end{proof}

\Rm  Note that there is an interesting similarities between this proof and the steps 
Figure 19 $\leadsto $ Figure 29 of
\cite{ak2}, which was crucial in showing that the $2$--fold covering space of
the Cappell--Shaneson's  fake ${\bf RP}^4$ is $S^4$, \cite{g2}.

{\Cor $\widetilde{\widetilde{H}}= S{^3} \times [0,1]$}

\begin{proof} Showing  $\widetilde{\widetilde{H}}= S{^3} \times [0,1]$ is equivalent to showing
that the $4$--manifold obtained by capping the boundaries of $\widetilde{\widetilde{H}}$ with $4$--balls is
diffeomorphic to $S^4$. Observe that under the $8$--fold 
covering map $\pi \co S^{3}\to Q $ the loop $C$ of Figure 122 lifts to a pair of linked Hopf
circles in $S^3$, each of it covering $C$ four times (this is explained in Figure 123).
By replacing $C$ in $\partial (W)$ by the whole stands of $1$-- and
$2$--handles going through that middle $1$--handle as in Figure 79 (all the handles other than
$H_{1}$), and by lifting  those $1$ and $2$--handles to $S^{3}$ we obtain the $8$--fold covering of
$H$, with ends capped off by $4$--balls. Since the monodromy
$\phi$ has order $3$, each strand has the monodromy $\phi^{4}=
\phi $. So we need to perform the Gluck construction as in Figure 124, which  after handle
slides becomes $\Sigma \# \Sigma =S^4$ (because we have previously shown that
$\Sigma =S^4$). Note that the bottom two 
 handlebodies of Figure 124 are nothing but $S^4$ Glucked along $A$, along with a cancelling pair of
$2$ and $3$--handles (as usual in this pair the $2$--handle is attached to the unknot on the boundary, ie, the
horizontal zero-framed circle, and the $3$--handle is not drawn). 
\end{proof}

{\Cor $\widetilde{W}_{+} =\widetilde{W}$}

\begin{proof} By inspecting the $2$--fold covering map in Figure 123, and by observing that
$\phi^{2}=\phi^{-1}$ we get the handlebody of  $\widetilde{W}_{+}$ in Figure 125. As before, since the
$-4$ framed handle is attached along the trivial loop on the boundary we get $\widetilde{W}_{+} =
\widetilde{W}\#\Sigma$, where
$\Sigma $ is the $S^4$ Glucked along $A$ (recall the previous Corollary), hence we have
$\widetilde{W}_{+}=\widetilde{W} = \mbox{ Euler class} -4 \mbox{ disk bundle over} \; S^{2}$ \end{proof}

\Rm An amusing fact: It is not hard to check that the $2$--knot complement $Y$ is
obtained by the $0$-{\it logarithmic transformation operation} performed along an imbedded
Klein bottle
${\bf K}$ in $\hat{M} - S^{1}\times B^{3}$ (which is $\hat{M}$ minus a 
$3$--handle) ie, in Figure 21. This
is done by first changing the
$1$--handle notation of Figure 21 (by using the arcs in Figure 126) to circle-with-dot notation, then by simply
exchanging a dot with the zero framing as indicated by the first picture of Figure 127. The
result is the second picture of Figure 127 which is $Y$. This operation is nothing other than
removing the tubular neighborhood
$N$ of  ${\bf K}$ from $\hat{M} - S^{1}\times B^{3}$ and putting it back by a
diffeomporphism which is the obvious involution on the boundary. 
It is also easy to check that by performing yet another $0$--logarithmic transformation
operation to $Y$ along an imbedded ${\bf K}$ gives $S^{1}\times B^{3}$ (this is 
Figure 77 $\Rightarrow $ Figure 78). So the operations 
$$S^{1}\times B^{3} \Rightarrow Y \Rightarrow \hat{M} - S^{1}\times B^{3}$$ are nothing but $0$--logarithmic
transforms along ${\bf K}$. Note that all of the $4$--manifolds $S^{1}\times B^{3}$, $Y$ and $M$
are bundles over $S^{1}$,with fibers $B^{3}$, $Q_{0}$, and
$T_{0}\times_{-I}S^{1}$  respectively.

\Rm Recall \cite{go1}  that a knot $\Sigma^{n-2} \subset S^{n}$ is said 
to admit a {\it strong ${\bf Z}_{m}$--action}
if there is a diffeomorphism $h\co  S^{n} \to S^{n}$ with
\begin{enumerate}
\item[(i)]  $ h^{m} =1 $ \smallskip
\item[(ii)] $ h(x)=x $ for every $\;x\in \Sigma^{n-2}$ \smallskip
\item[(iii)] $ x, h(x), h^{2}(x),..,h^{m-1}(x) $ are all distinct for every $ \;x\in
S^{n}-\Sigma^{n-2} \;$
\end{enumerate} 

By the proof of the Smith conjecture when $n=3$ the only knot that admits a strong ${\bf Z}_{m}$ action is the
unknot. For $n=4$ in \cite{gi} Giffen found knots that admit strong ${\bf Z}_{m}$ actions when $m$ is odd. Our
knot $A\subset S^4$ provides an example of knot which admits a strong ${\bf Z}_{m}$ action for $m \neq 0$ mod
$3$. This follows from Proposition 7, and from the fact that  $A$ is a fibered knot with an order $3$
monodromy.

\Rm Recall the vertical picture of H in Figure 77, appearing as $ W_{+}-W $. We can place $H$
vertically on top of $Q\times I$ (Figure 4) by identifying  $Q_{+}$ with $Q\times 1$
$$ Z:=H\smile_{Q_{+}} (Q\times I)$$ 
Resulting handlebody of $Z$ is Figure 128. As a smooth manifold $Z$ is nothing other than a copy of $H$.
So Figure 128 provides an alternative handlebody picture of $H$ (the other one being Figure 47).

\Rm Let $X^{4}$ be a smooth $4$--manifold, and $C\subset X^{4}$ be any loop with the property that $[C]\in
\pi_{1}(X)$ is a torsion element of order $\pm 1$ mod $3$, and $U\approx S^{1}\times B^{3}$ be the open tubular
neghborhood of $C$. We can form:
$$\hat{X}:= (X-U)\smile_{\partial} Y $$
Recall that $ \pi_{1}(Y)=\< t,a | t^3=a^3, ata=tat\>$,  so by Van-Kampen theorem we
get $\pi_{1}(\hat{X})=\pi_{1}(X)$, in fact
$\hat{X}$ is homotopy equivalent to $X$. In particular by applying this process to $X= M^{3}\times I$, where
$M^{3}$ is a $3$--manifold whose fundamental group contains a torsion element of order $\pm 1$ mod $3$, we
can construct many examples of potentially nontrivial $s$--cobordisms $\hat{X}$ from $M$ to itself.  

\eject

\hbox{}\par\includegraphics[width=\hsize]{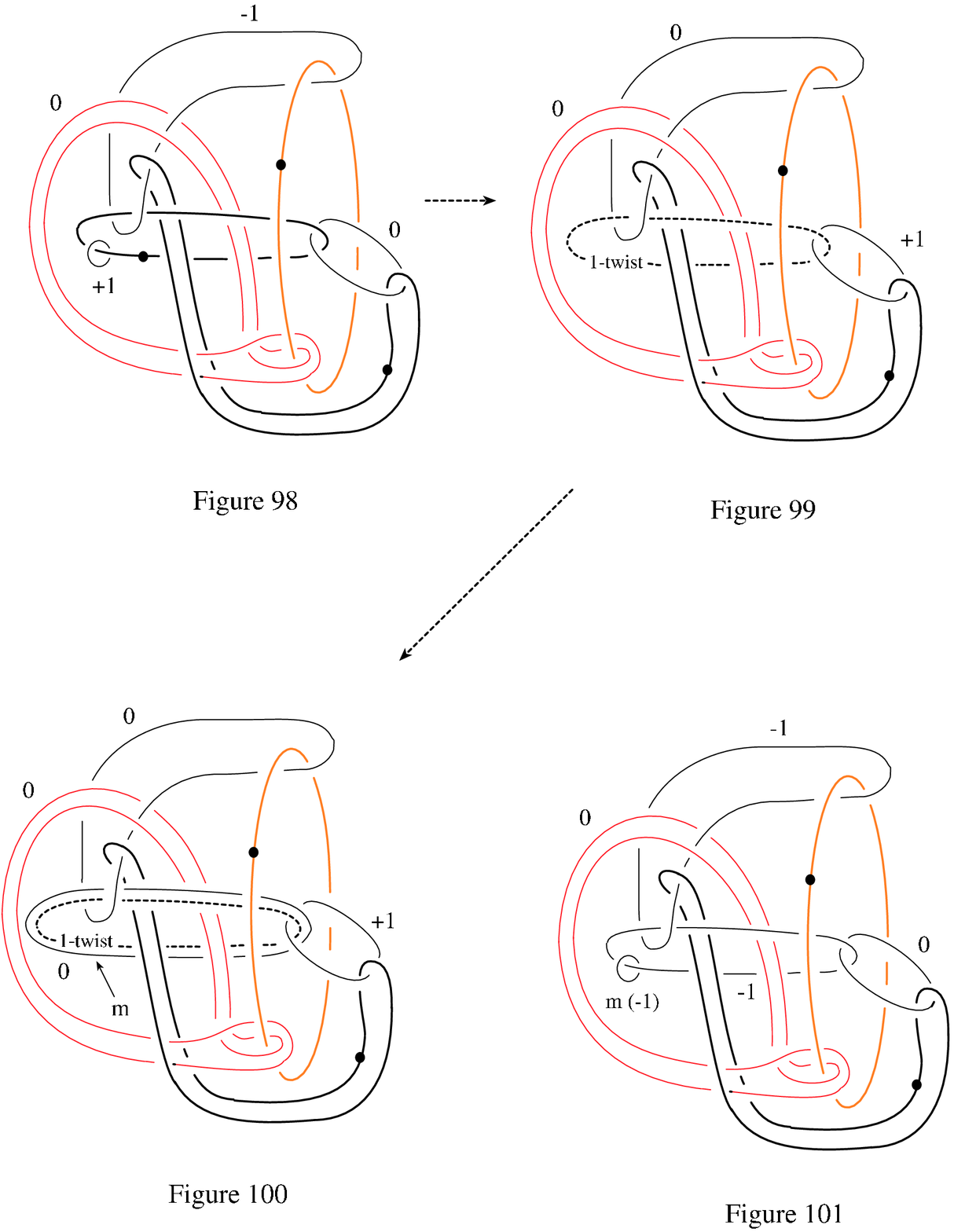}

\newpage
\hbox{}\par\includegraphics[width=\hsize]{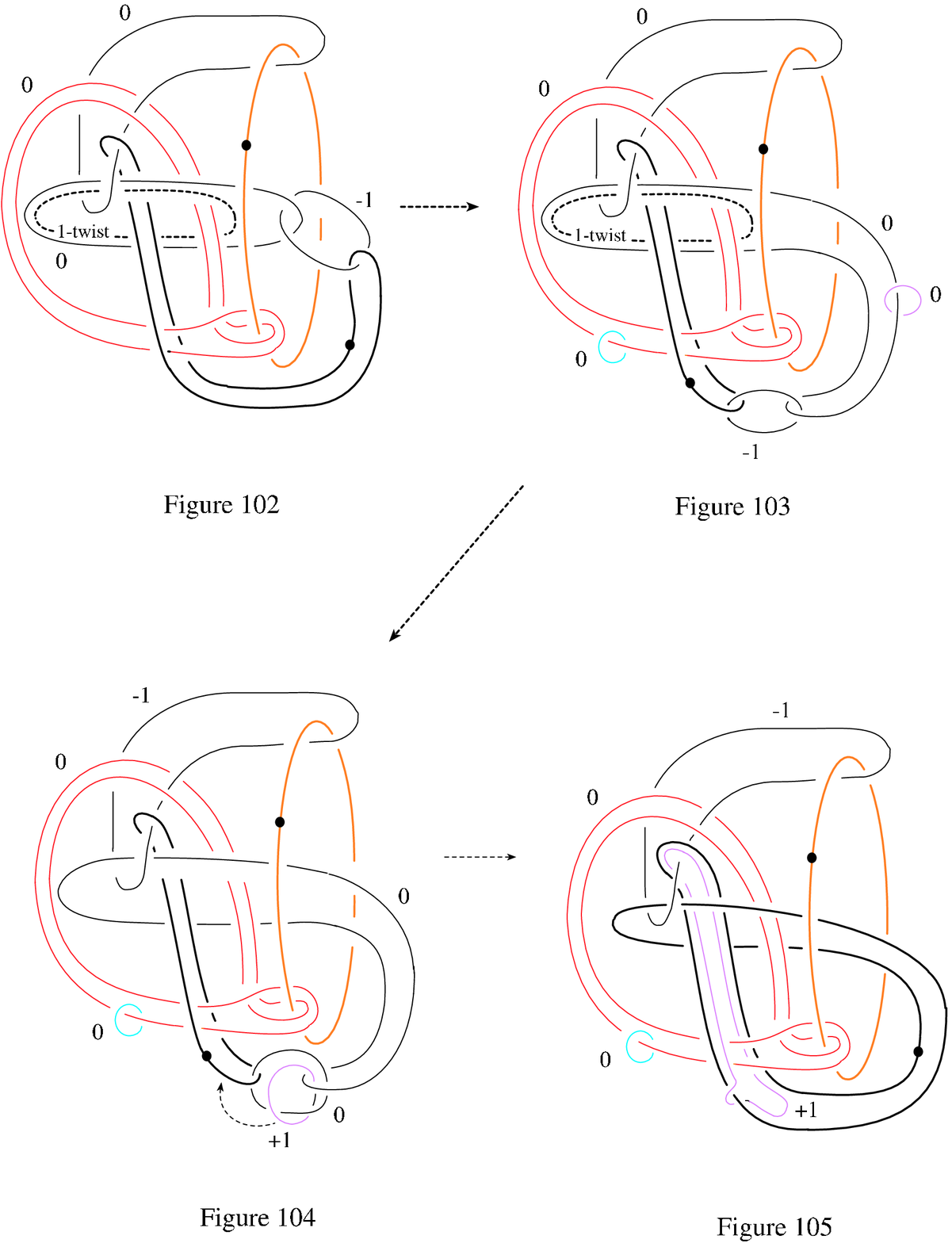}

\newpage
\hbox{}\par\includegraphics[width=\hsize]{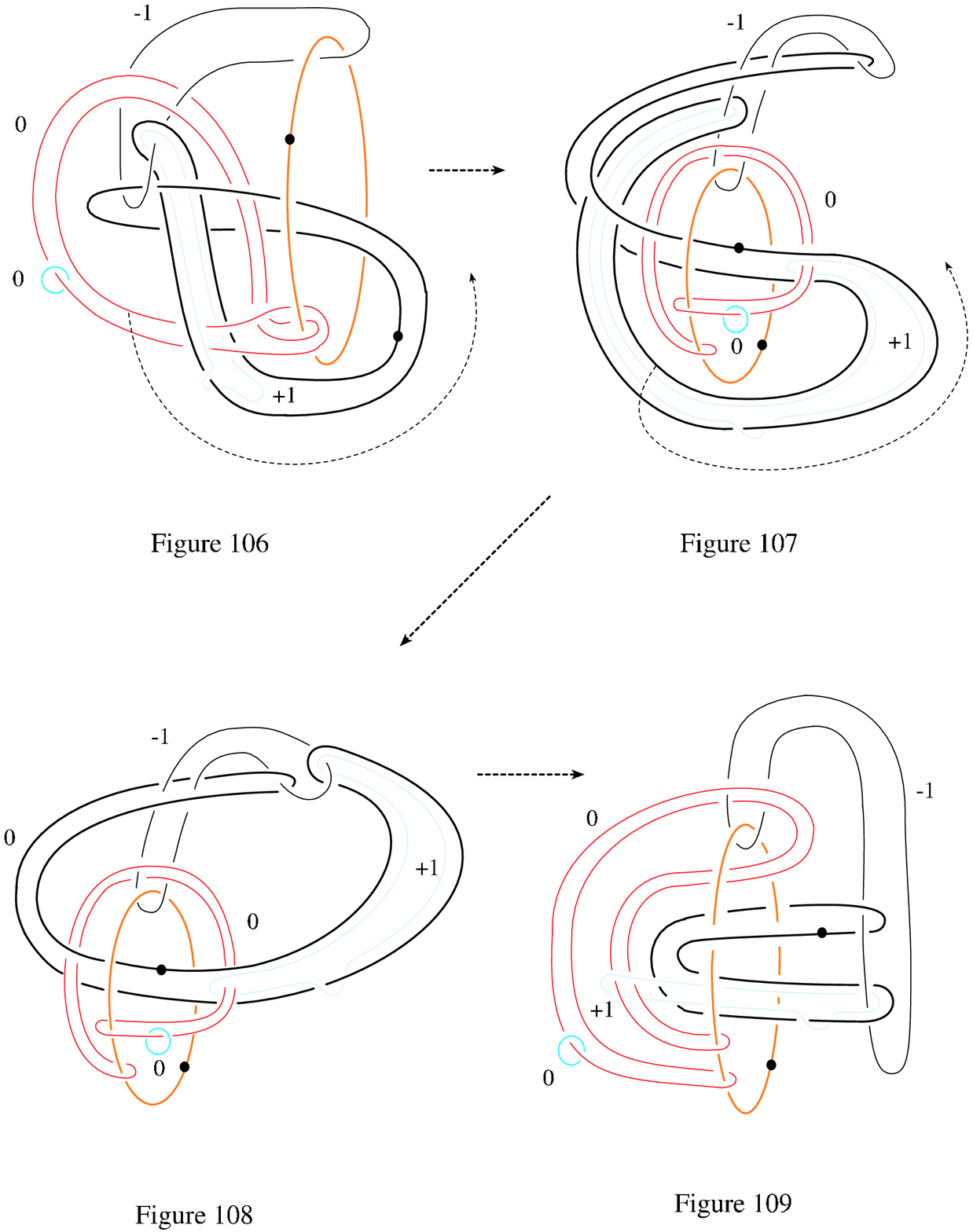}

\newpage
\hbox{}\par\includegraphics[width=\hsize]{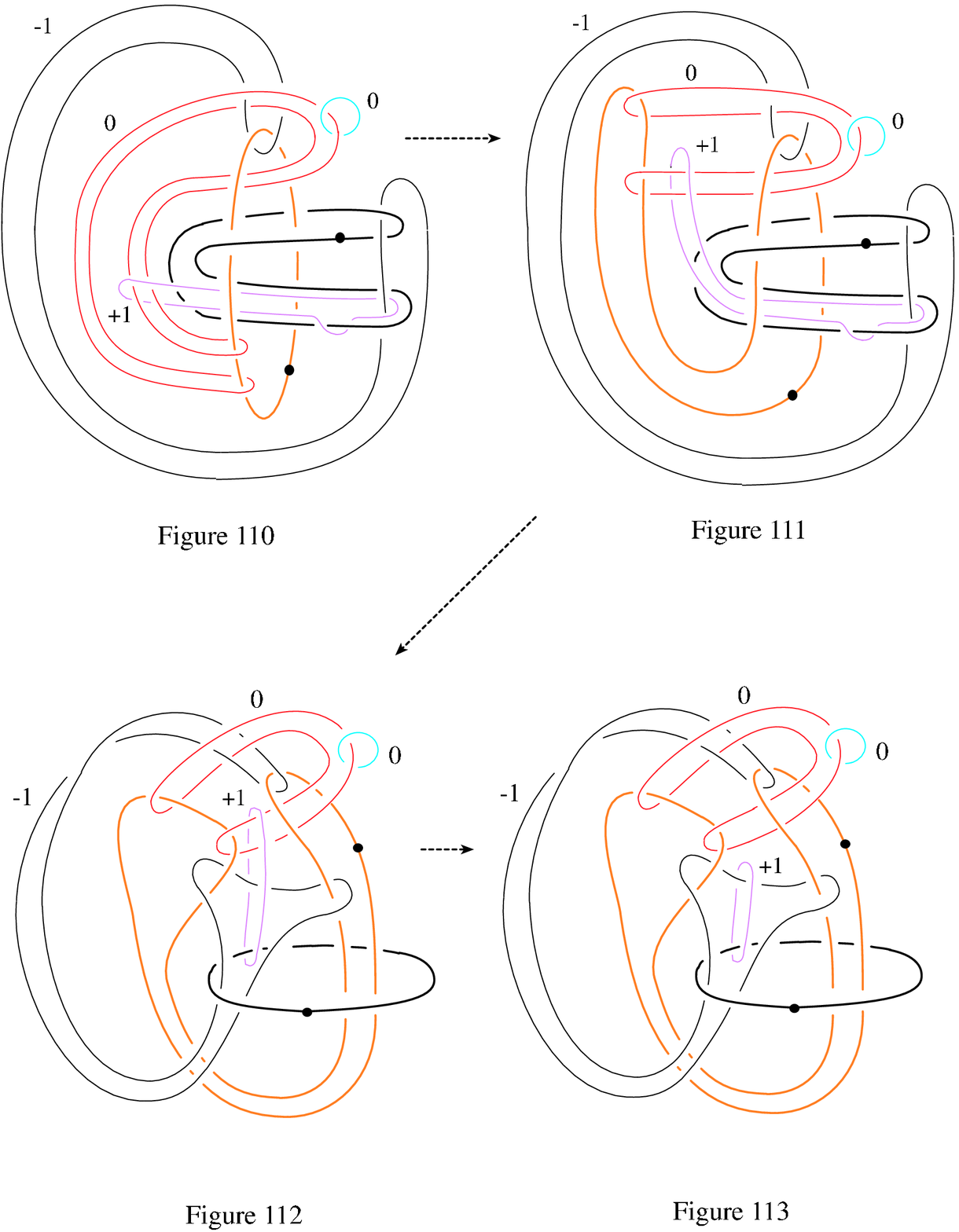}

\newpage
\hbox{}\par\includegraphics[width=\hsize]{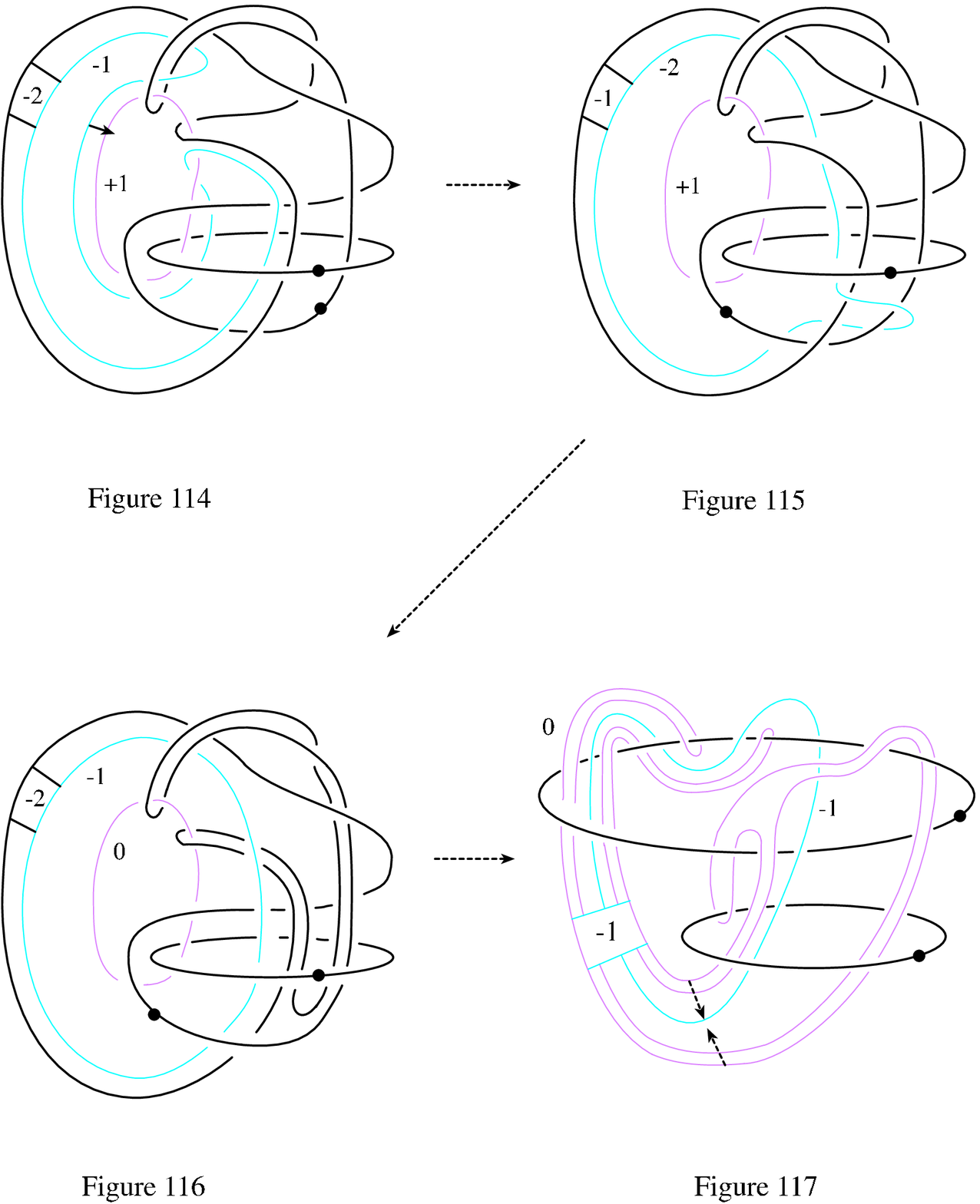}

\newpage
\hbox{}\par\includegraphics[width=\hsize]{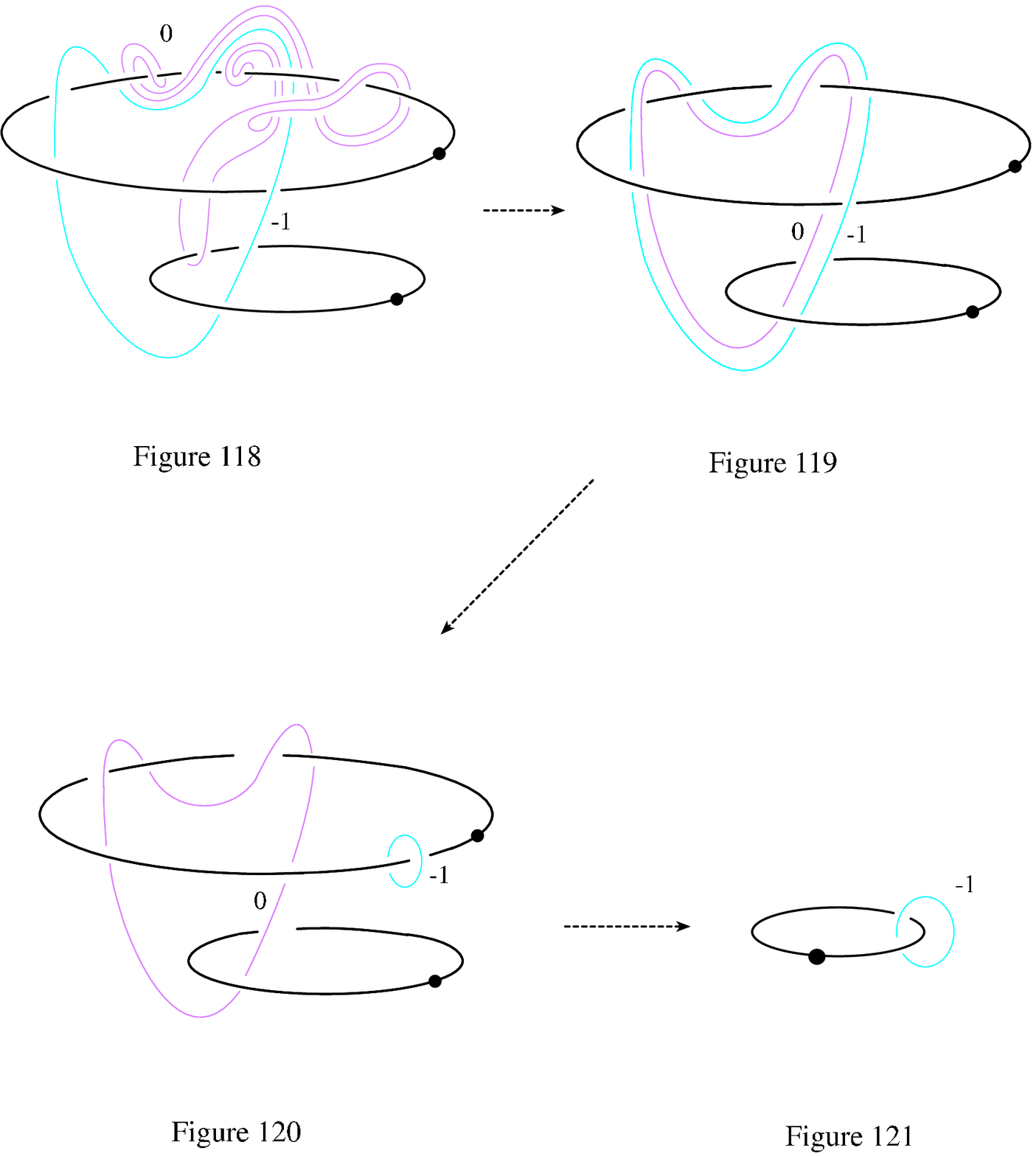}

\newpage
\cl{\includegraphics[width=4in]{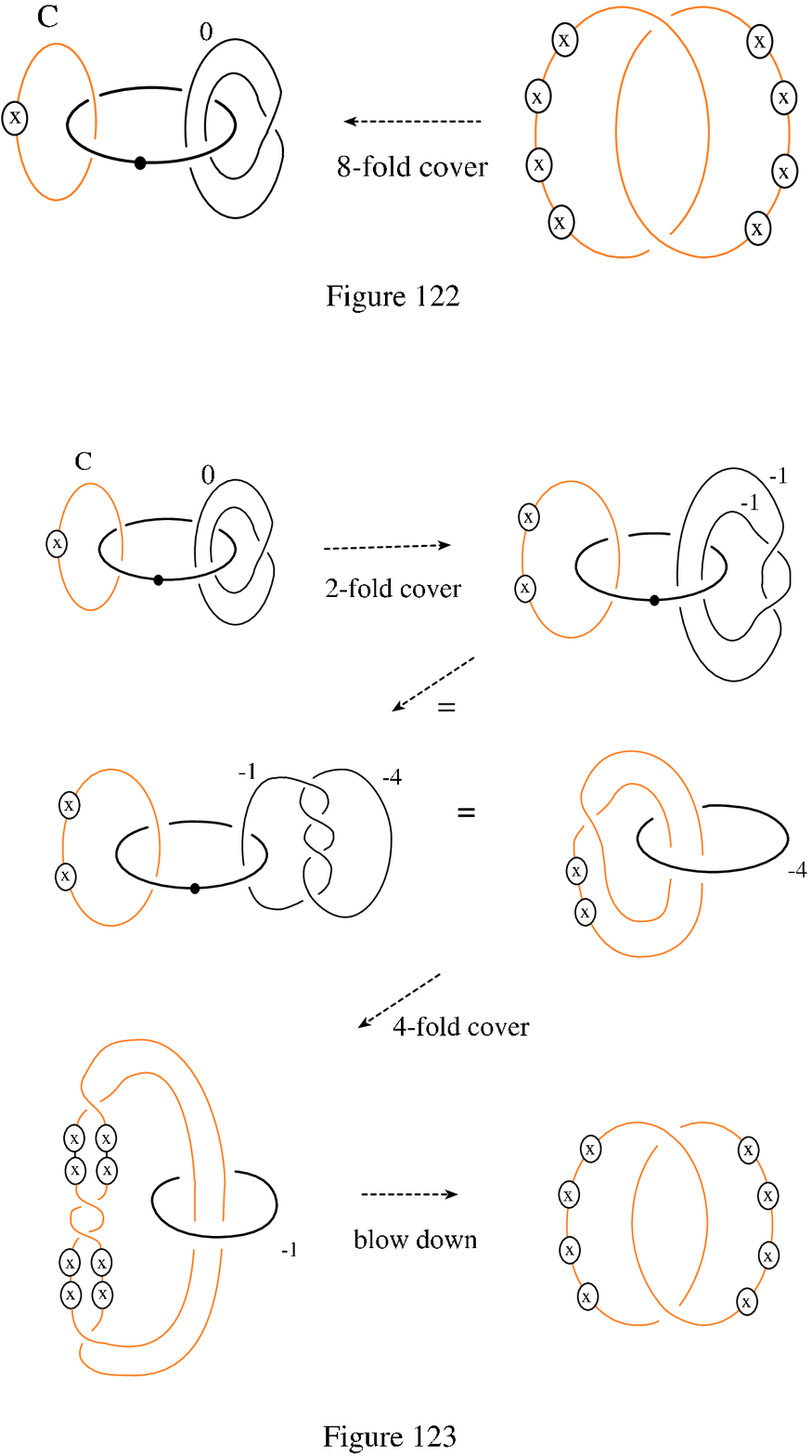}}

\newpage
\hbox{}\par\includegraphics[width=\hsize]{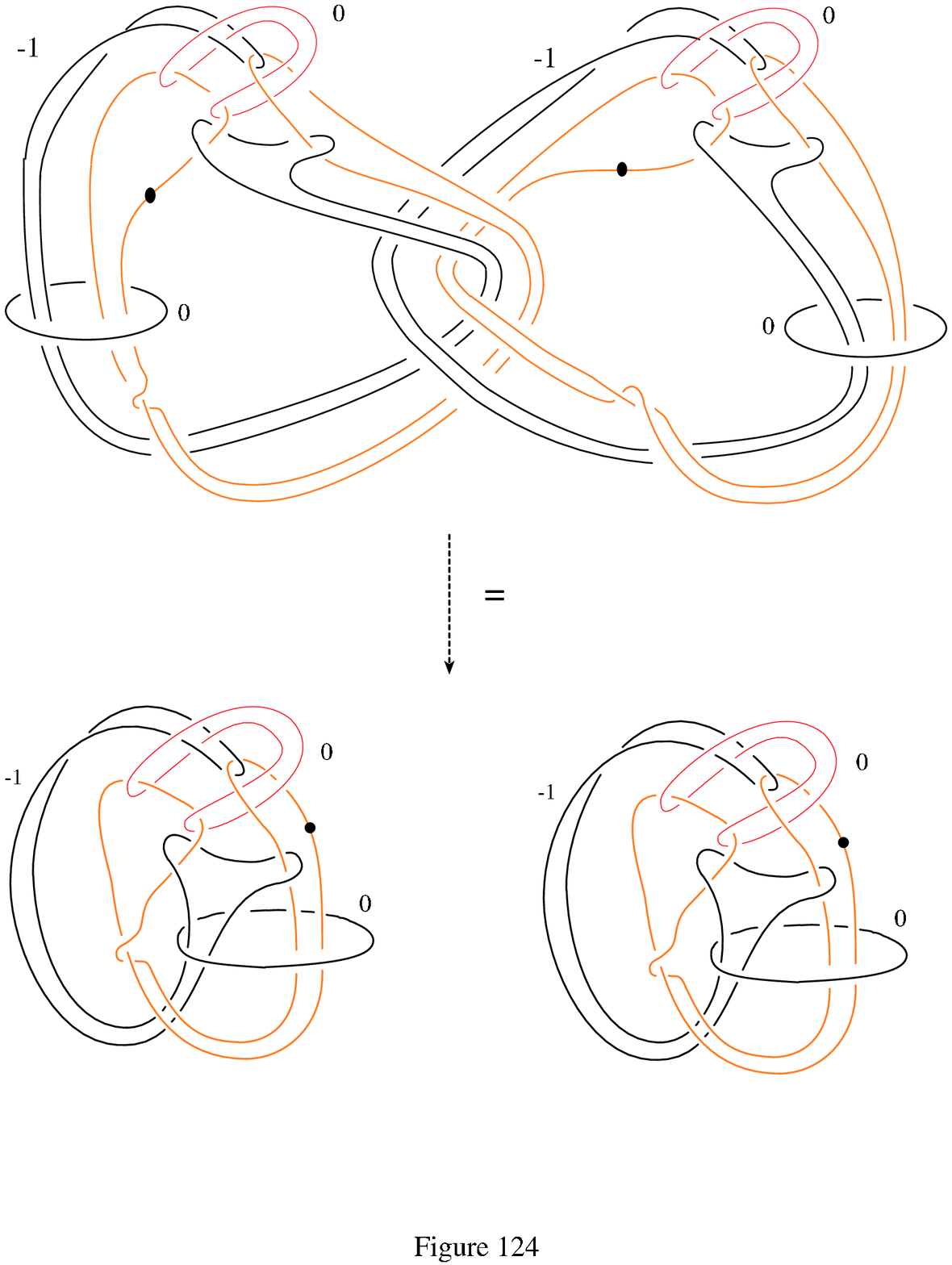}

\newpage
\hbox{}\par\includegraphics[width=\hsize]{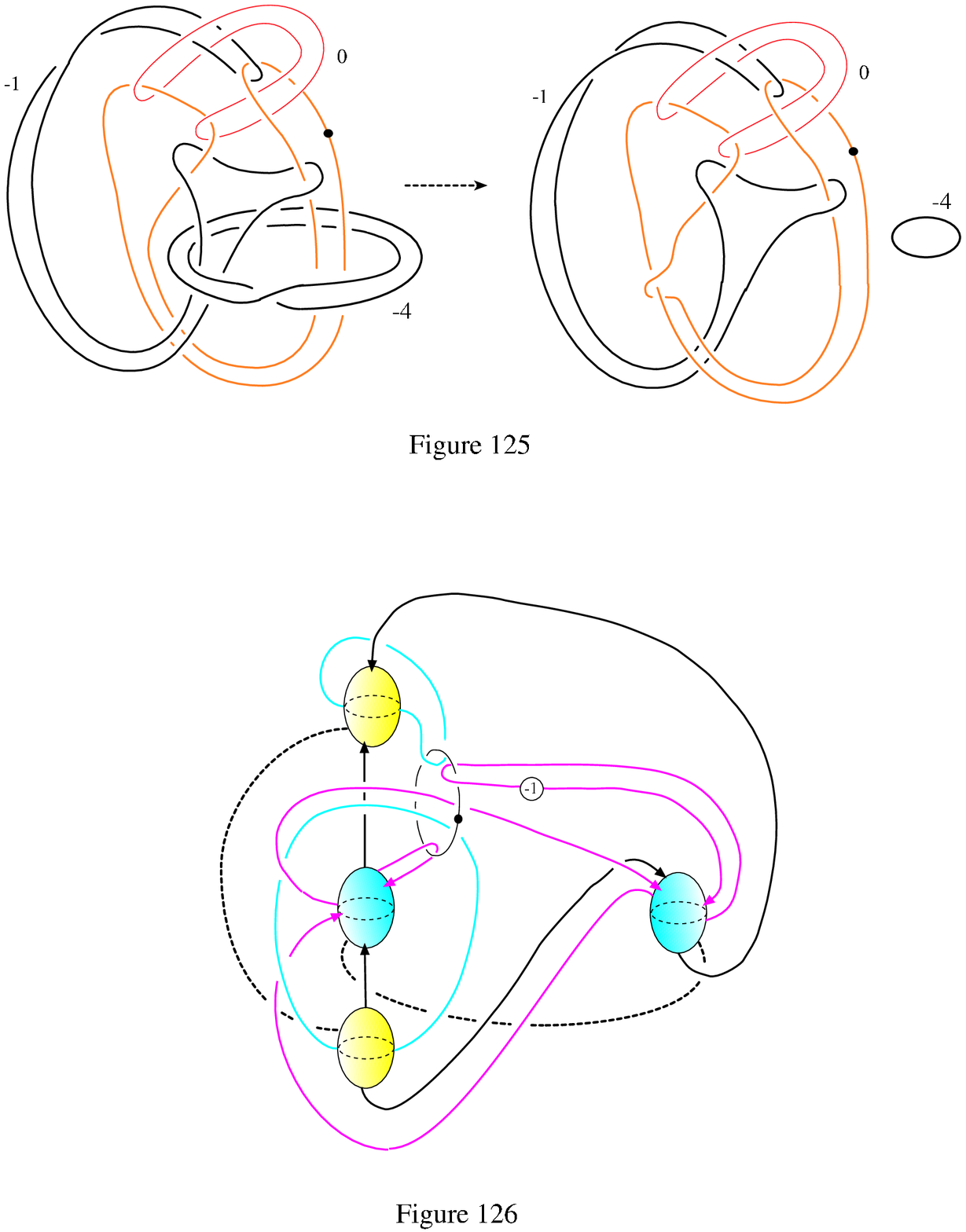}

\newpage
\hbox{}\par\includegraphics[width=\hsize]{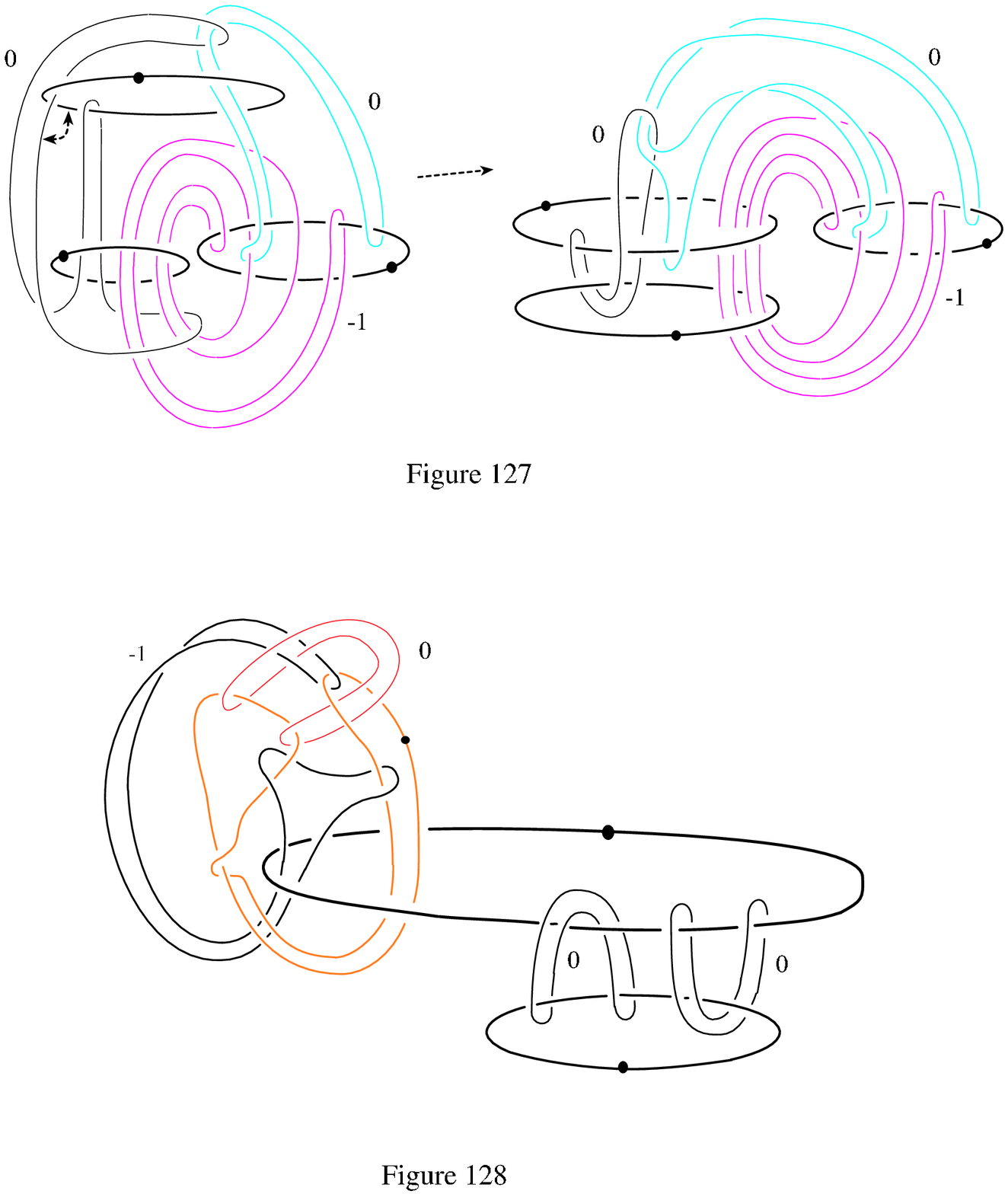}

\newpage

\section{ More on A }

In this section we will give even a simpler and more concrete description of the $2$--knot $A$. 
As a corollary we will show that $W_{+}$ is obtained from $W$ by performing the Fintushel--Stern knot surgery
operation by using the trefoil knot $K$ (\cite{fs}). 
It is easy to check that Figure 129, which describes $A$  (recall Figure 82), 
is isotopic to Figure 130. To reduce the clutter, starting with Figure 130 we will denote the small
$0$--framed circle, which links $A$ twice, by a single thick circle $\delta$. 

By performing the indicated handle slide to the handlebody of
Figure 130 we arrive to Figure 131, which can be drawn as Figure 132. By introducing a cancelling pair of $1$
and $2$--handles to the handlebody of Figure 132, then sliding them over other handles, and again cancelling
that handle pair we get Figure 133 (this move is self explanatory from the figures). Then by an isotopy we get
Figure 134, by the indicated handle slide we arrive to Figure 135. By drawing the ``slice $1$--handle'' 
(see \cite{a2}) as a $1$--handle and a pair of $2$--handles we get the diffeomorphism Figure 135 $\leadsto$
Figure 136, and a further handle slide gives Figure 137, which is an alternative picture of the $2$--knot
complement $A$. The reader can check that the boundary of Figure 137 is standard by the boundary diffeomorphism
Figure 137 $\leadsto$ Figure 138, which consists of a blowing-up + handle sliding + blowing-down operations
(done three times).

A close inspection reveals that the handlebody of the $2$--knot complement $A$ in Figure 135 is 
the same as Figure 139. Figure 139 gives another convenient way of checking that the boundary of this
handlebody is standard (eg, remove the dot from the slice $1$--handle and perform blowing up and sliding
and blowing down operations, three times, as indicated by the dotted lines of Figure 140). Now we can also
trace the loop $\delta$ into Figure 140, so Figure 140 becomes handlebody of $W_{+}$. By drawing the slice
$1$--handle as a
$1$--handle and a pair of $2$--handles  we get a diffeomorphism Figure 140 $\leadsto$ Figure 141. Clearly Figure
142 is diffeomorphic to Figure 141. Now by introducing a cancelling pair of $2$ and $3$--handles we obtain the
diffeomorphism Figure 142
$\leadsto$ Figure 143 (it is easy to check that the new $2$--handle of Figure 143 is attached along the
unknot on the boundary).

Now, let us recall the Fintushel--Stern knot surgery operation \cite{fs}: 
Let $X$ be a smooth $4$--manifold containing an imbedded torus $ T^{2} $ with trivial normal bundle, and
$ K\subset S^{3} $ be a knot. The operation $ {\bf X\leadsto X_{K}} $ of replacing a tubular
neighborhood of $T^{2}$ in $X$ by $ (S^{3} - K)\times S^{1}$ is the so called Fintushel--Stern
knot surgery operation. In
\cite{a1} and \cite{a3} an algorithm of describing the handlebody of $X_{K}$ in terms of the handlebody of $X$ is
given. From this algorithm we see that Figure 144 $\leadsto$ Figure 143 is exactly the operation 
$W\leadsto W_{K}$ where $K$ is the trefoil knot. And also it is easy to check that Figure 144 $\leadsto$
Figure 145 describes a diffemorphism to $W$. Hence we have proved:

{\Prop $W_{+}$ is obtained from $W$ by the Fintushel--Stern knot surgery operation along an imbedded torus by
using the trefoil knot K}

\Rm Note that we in fact proved that the knot complement $S^{4}-A$ is obtained by from 
$S^{1}\times B^{3}$ by the Fintushel--Stern knot surgery operation along an imbedded torus by
using the trefoil knot K. Unfortunately this torus is homologically trivial; if it wasn't,
 from \cite{fs}, we could have  concluded that $W_{+}$ (hence $H$) is exotic.

\Rm Now it is evident from From Figure 139 that $A$ is the
 $3$--twist spun of the trefoil knot (\cite{z}). This explains why $A$ is the fibered knot with
fibers $Q$ (which is the $3$--fold branched cover of the trefoil knot). After this paper was
written, we were pointed out that in \cite{go2} it  had proven that the Gluck
construction to a twist-spun knot gives back $S^{4}$. So in hind-sight we could have
delayed the Proposition 7 until this point and deduce its proof from \cite{go2}, but this
would have altered the natural evolution of the paper. Our hands-on proof of Proposition 7
should be seen as a part a  general technique which had been previously
utilized in \cite{a2}, \cite{ak2}.

Finally, note that if $A_{n}$ is the $n$--twist spun of the trefoil
knot (Figure 146), then one can check that its fundamental group generalizes 
the presentation of $A$:
$$ G_{n}= \< t,a \;| a=t^{-n}at^{n}, ata=tat\>=\< t,a \;| t^n=a^n, ata=tat\>$$

\newpage
\cl{\includegraphics[width=4in]{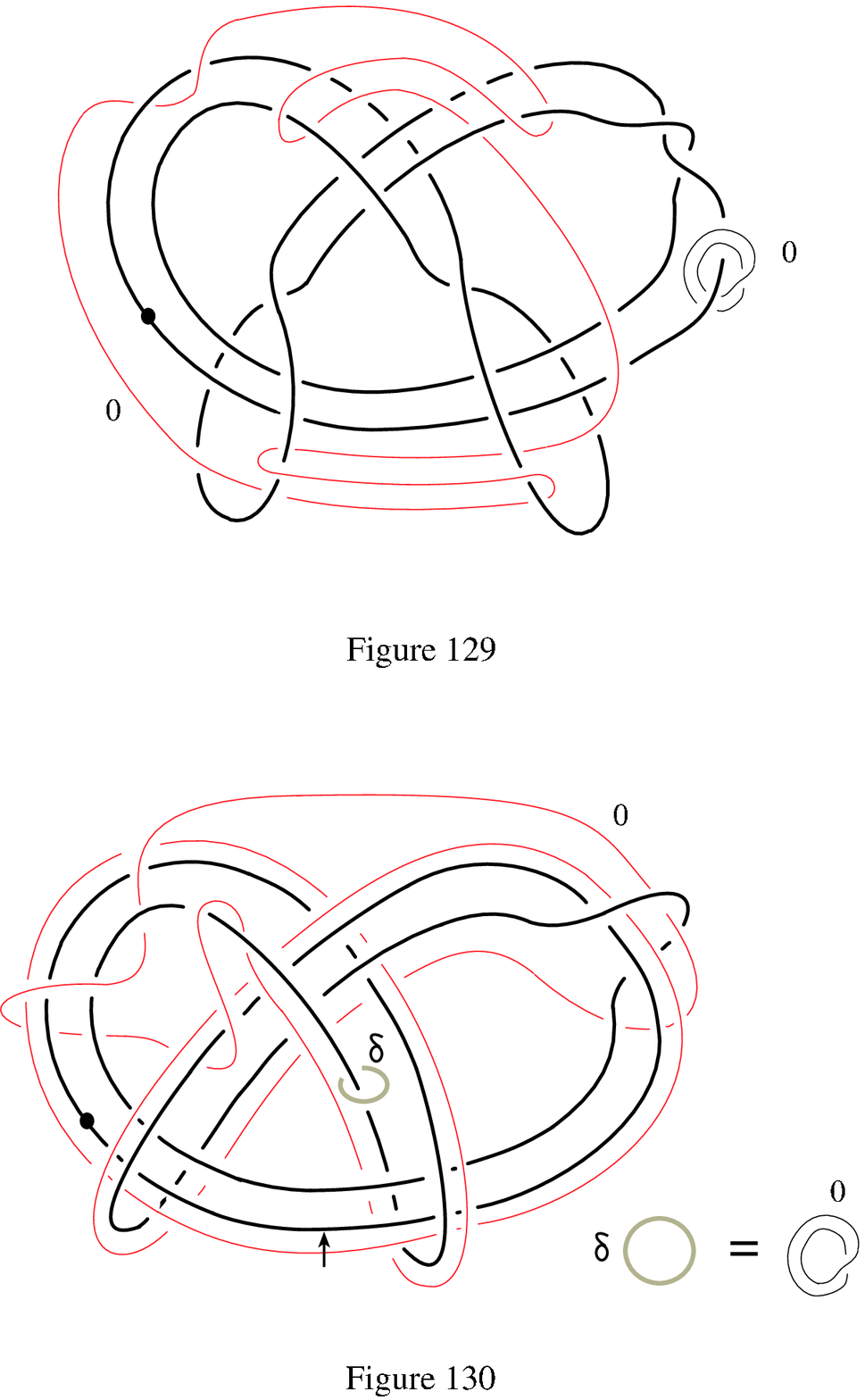}}

\newpage
\cl{\includegraphics[width=4in]{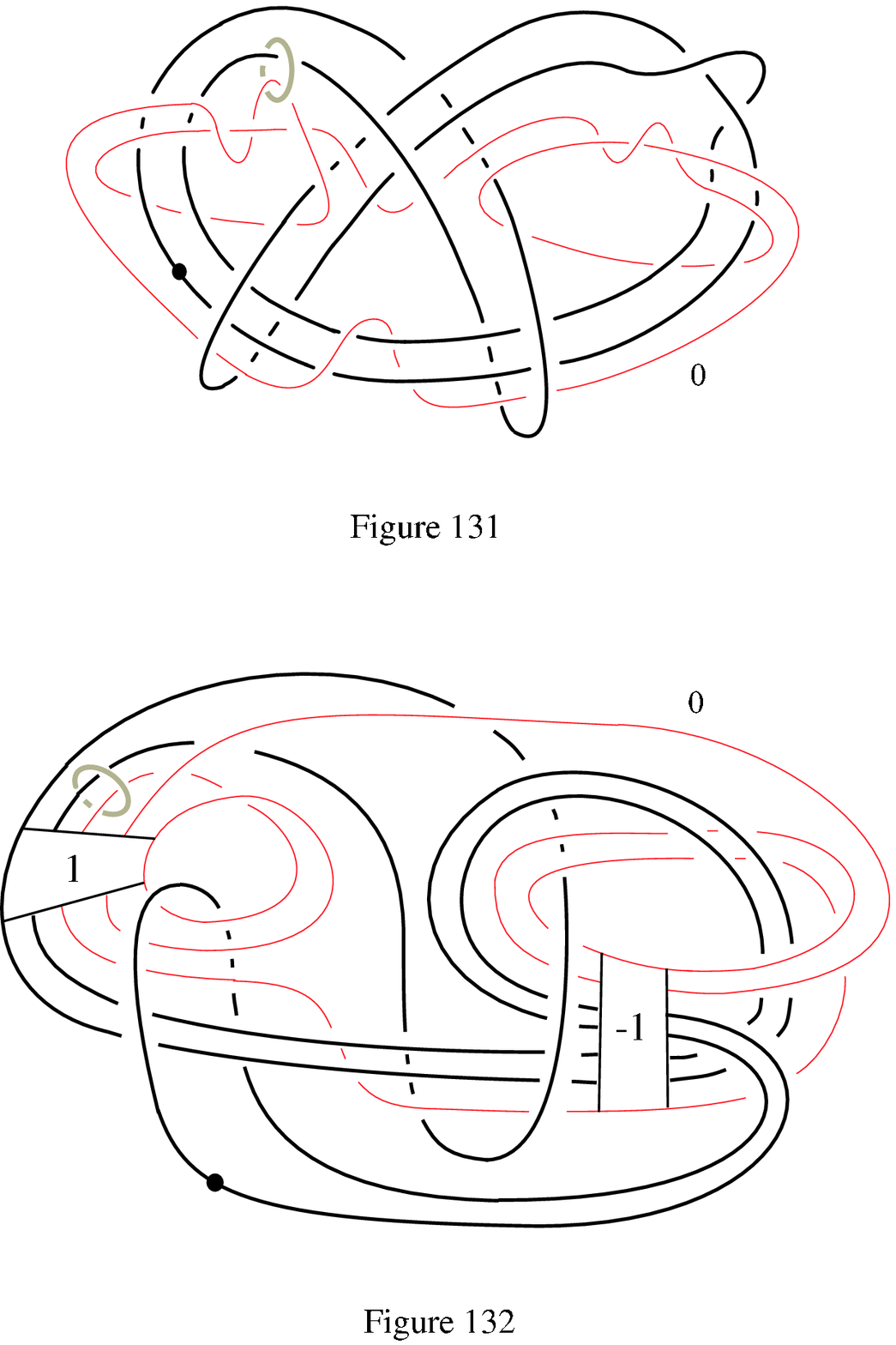}}

\newpage
\cl{\includegraphics[width=4in]{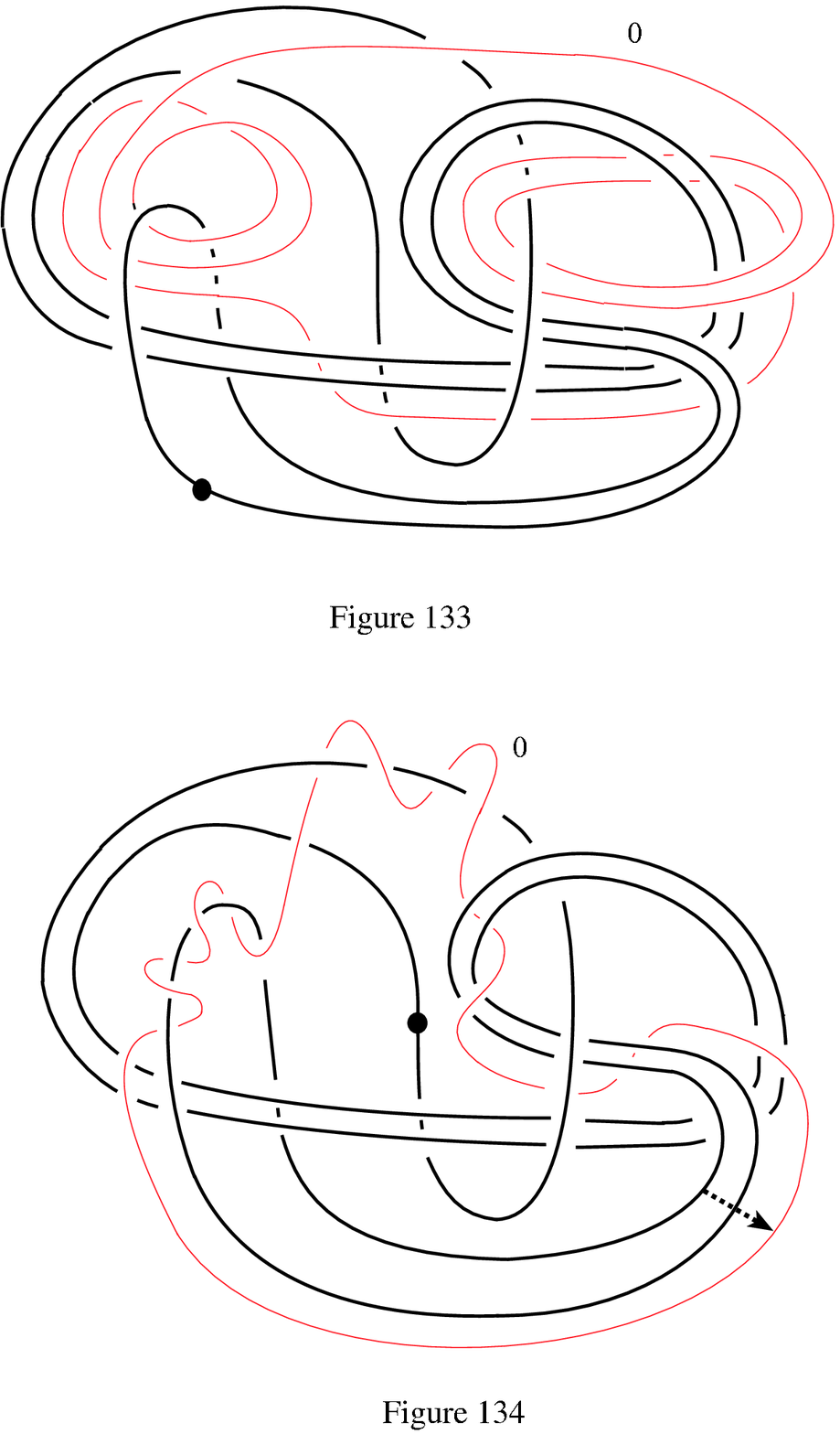}}

\newpage
\cl{\includegraphics[width=4in]{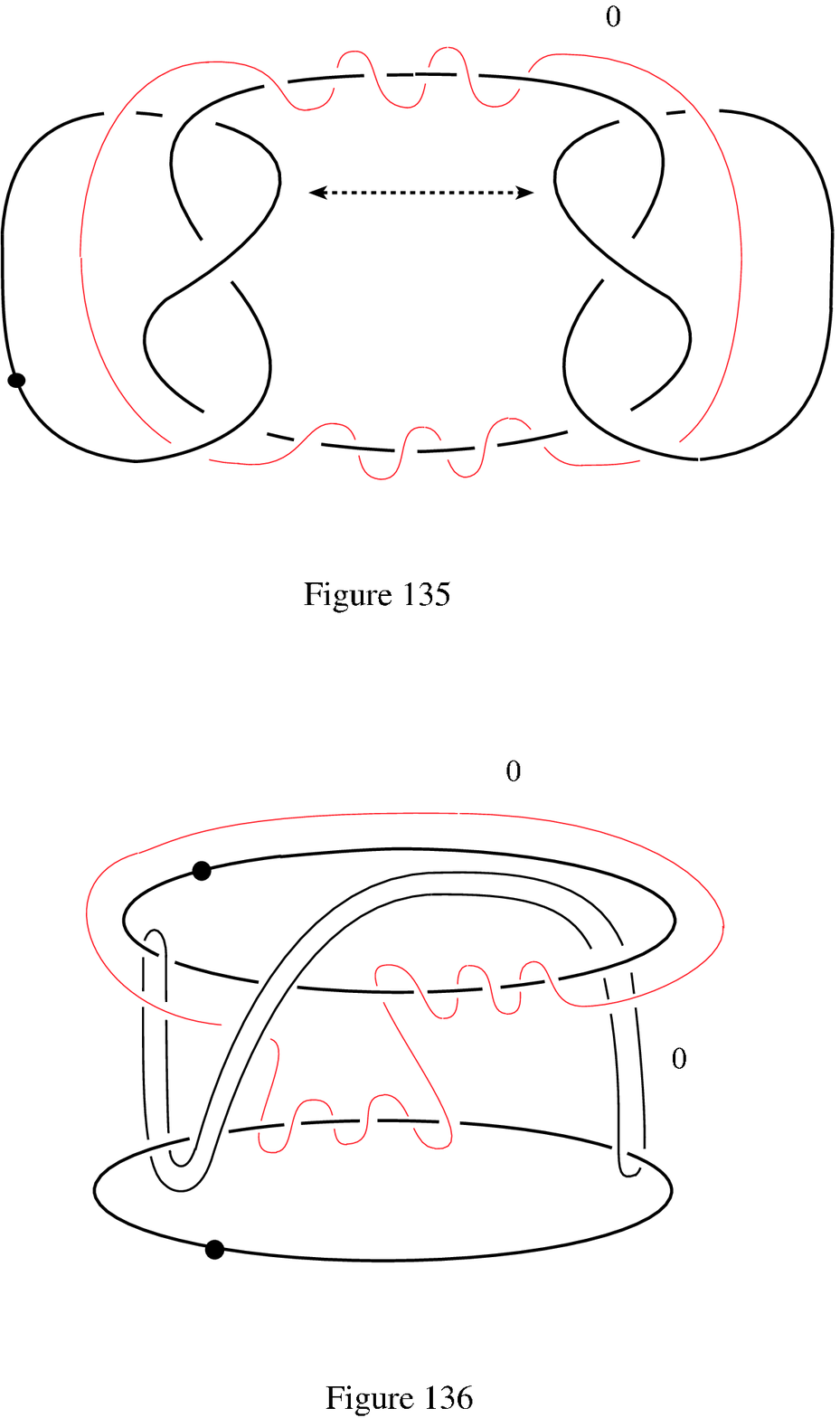}}

\newpage
\cl{\includegraphics[width=3.5in]{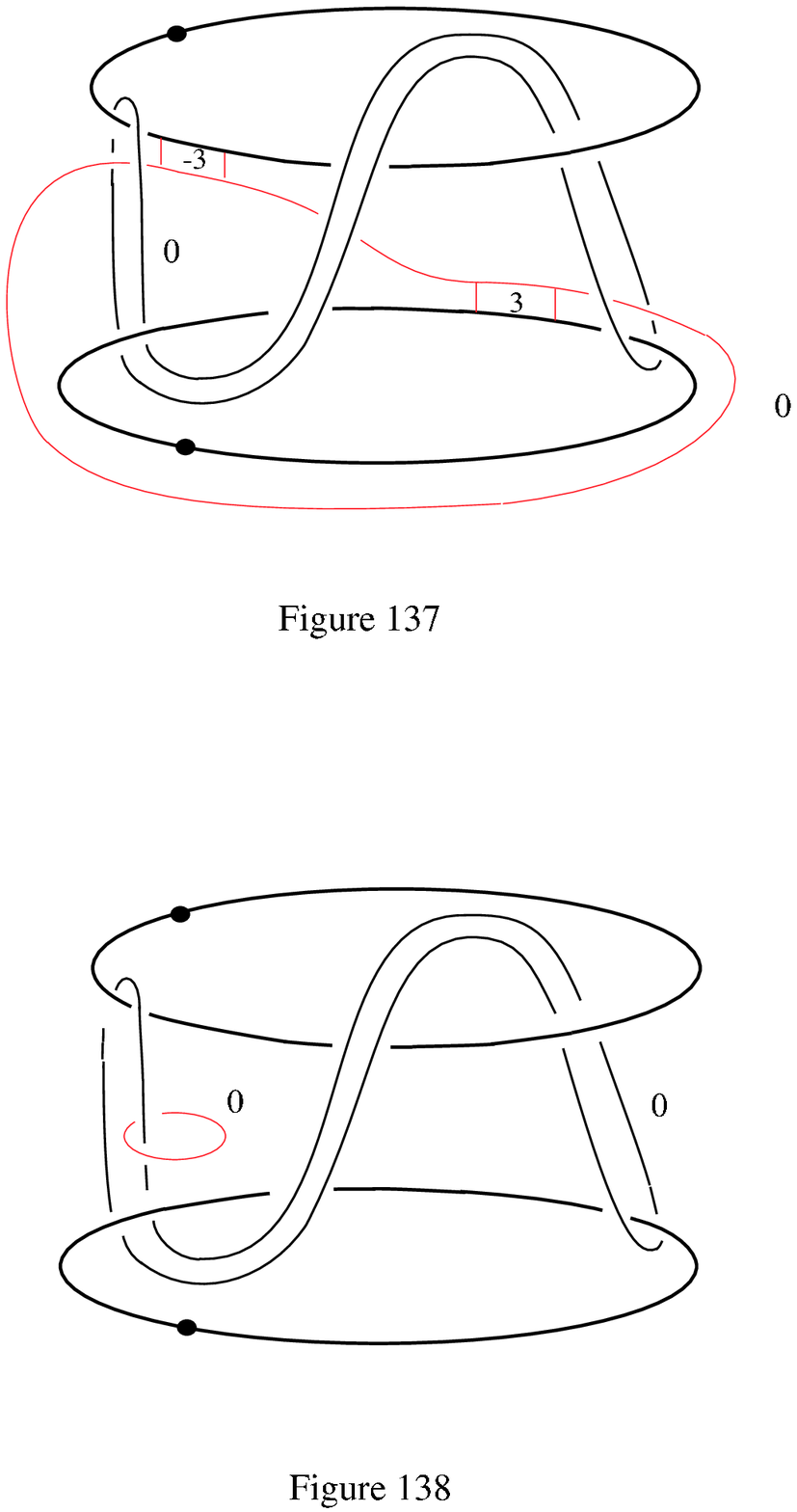}}

\newpage
\cl{\includegraphics[width=4in]{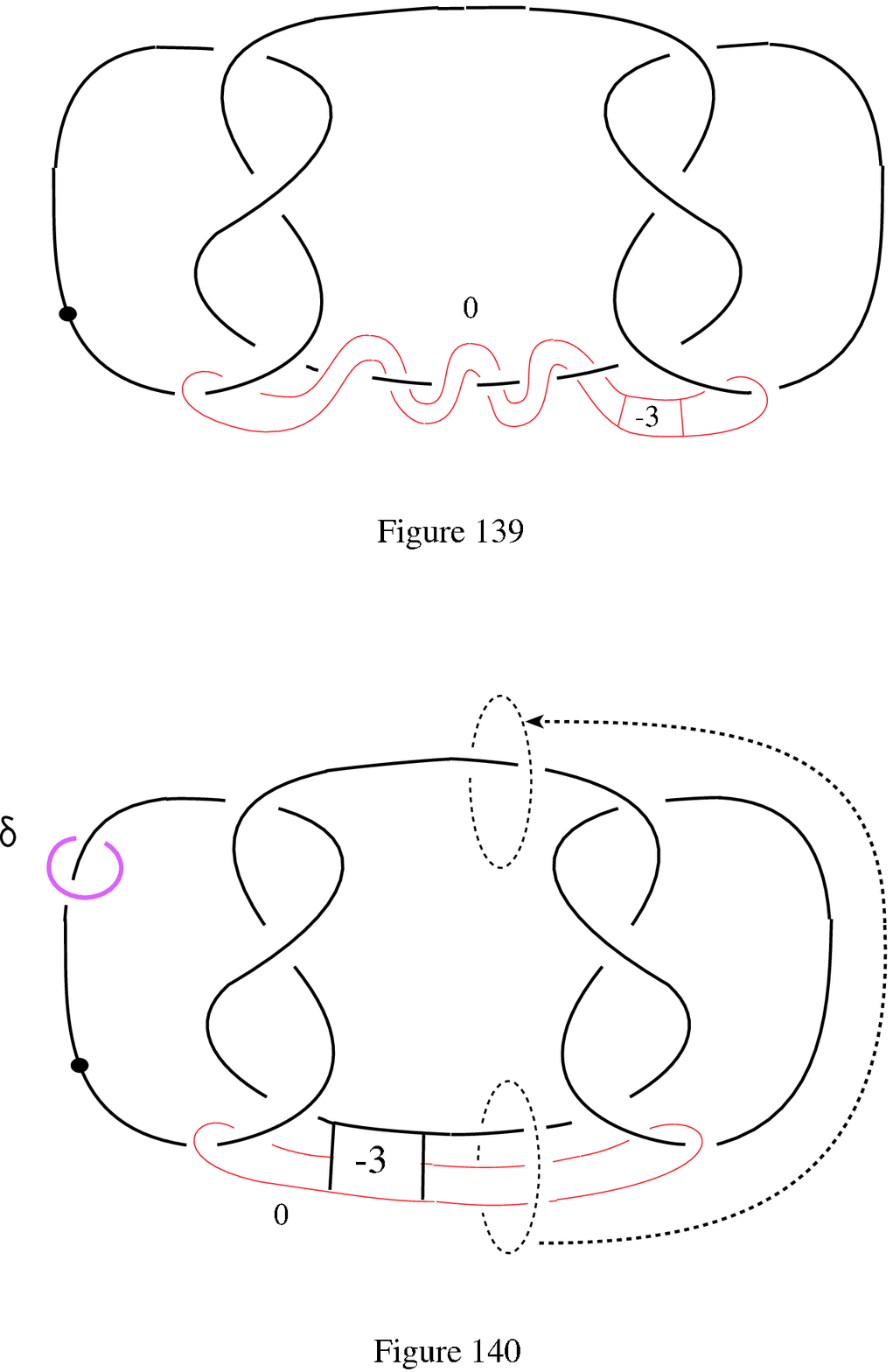}}

\newpage
\cl{\includegraphics[width=2.8in]{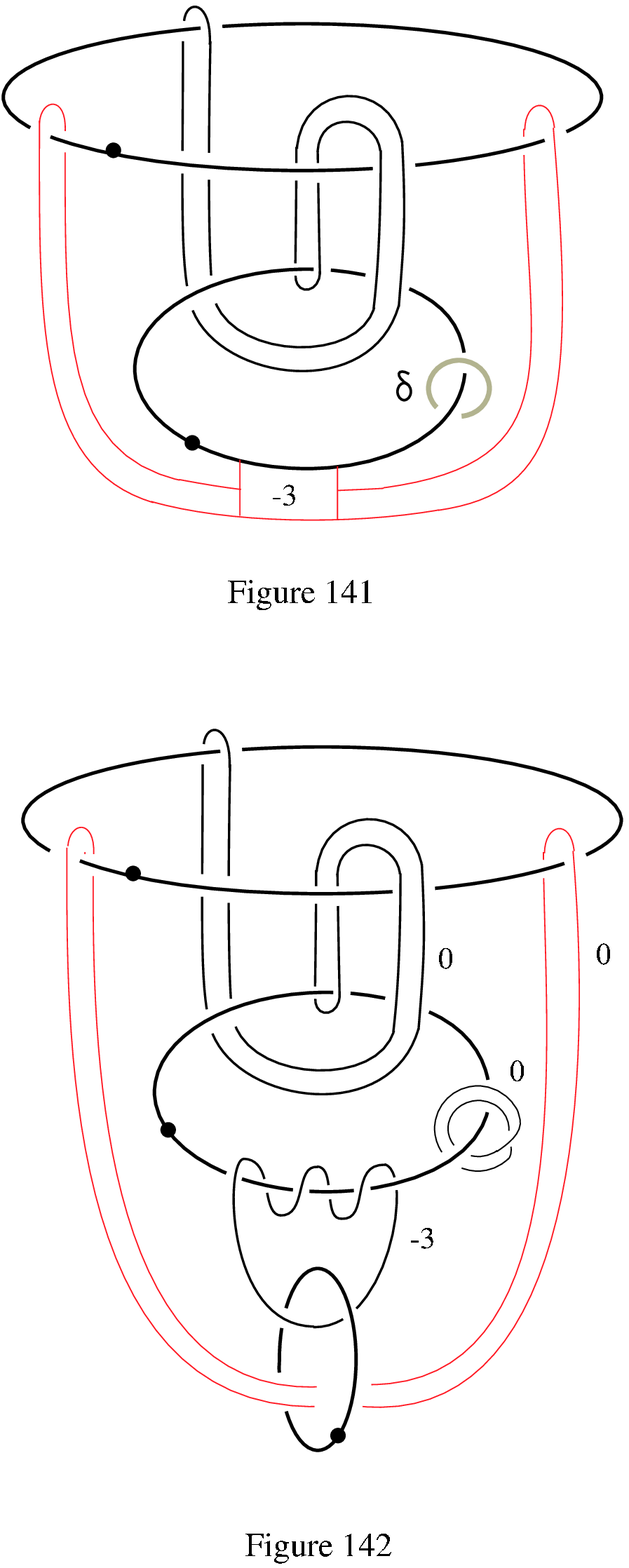}}

\newpage
\cl{\includegraphics[width=4.2in]{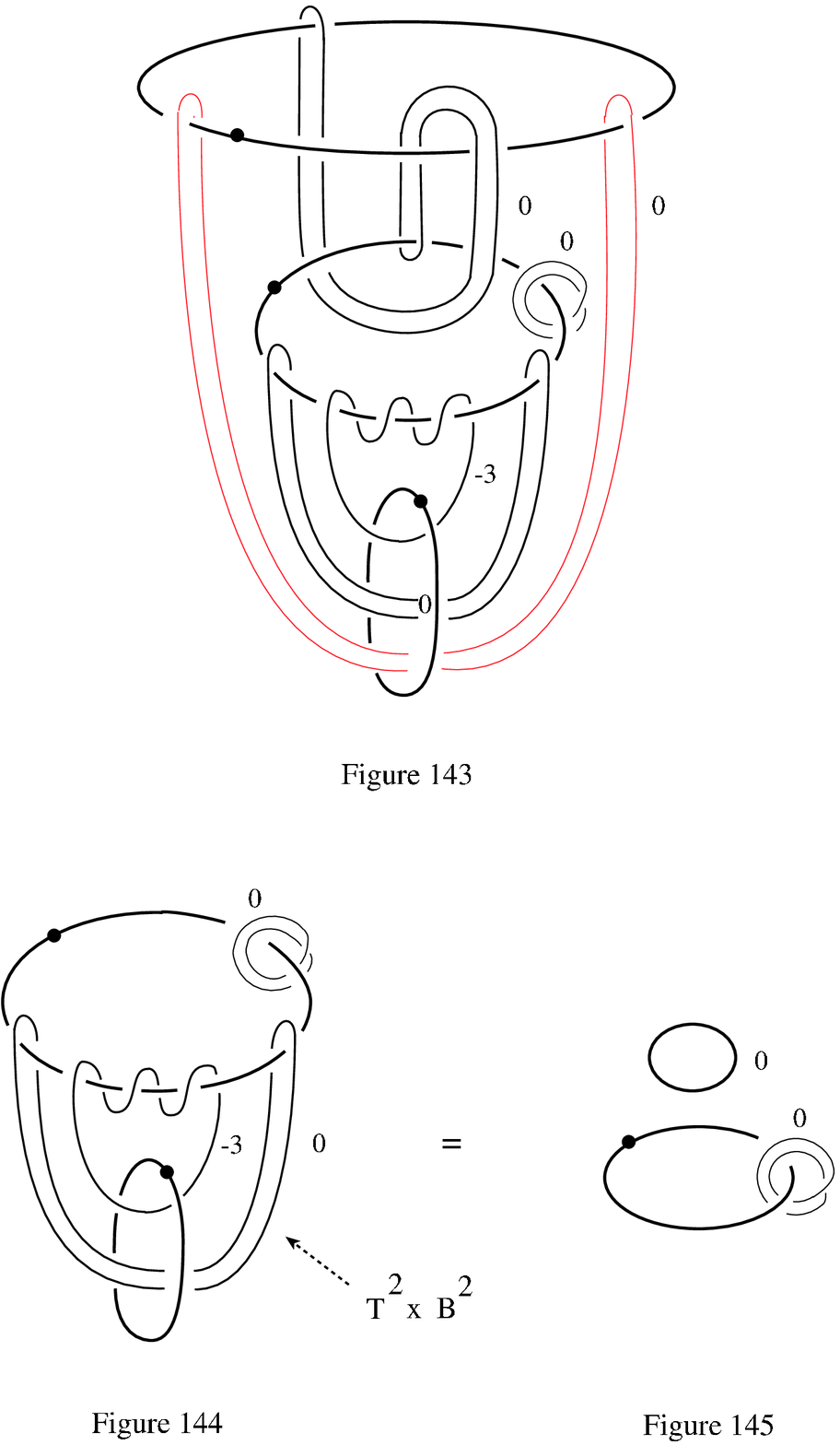}}

\newpage
\hbox{}

\cl{\includegraphics[width=3.5in]{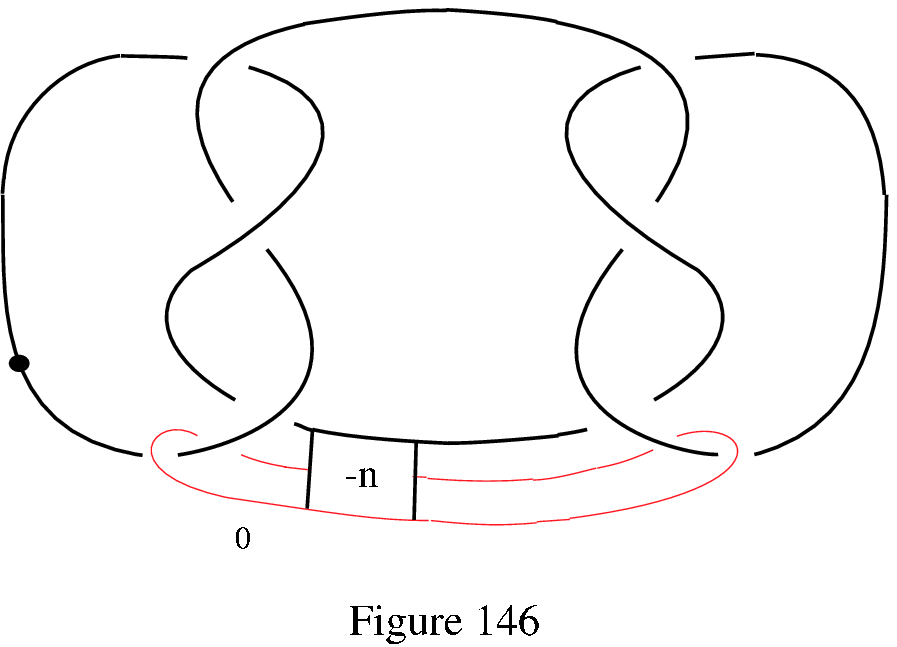}}

\vglue 0.5in

\end{document}